\overfullrule 0pt

\magnification=1200
\hsize=29pc
\frenchspacing
\advance\vsize-14pt 
\def\u{\underline}
\def\d{d\llap{\raise 1.9mm\hbox{\vrule height 0.16pt depth 0pt width 0.11cm}}}

\def\limind{\setbox1=\hbox{\oalign{\vadjust{\vskip -2pt}%
  \rm lim\cr \vadjust{\vskip -2pt}
    \hidewidth$-\mkern -12mu\rightarrow$\hidewidth\cr}}
    ÊÊ\mathop{\box1}}
\def\limproj{\setbox1=\hbox{\oalign{\vadjust{\vskip -2pt}%
  \rm lim\cr \vadjust{\vskip -2pt}
    \hidewidth$\leftarrow\mkern -12mu-$\hidewidth\cr}}
    ÊÊ\mathop{\box1}}
\font\smgoth=eufm10 at7pt
\font\goth=eufm10 at10pt
\font\Goth=eufm10 at13pt
 at14pt
\font\sm=cmr10 at8pt
\font\smint=cmr10 at9pt
\font\smitint=cmmi10 at9pt
\font\big=cmbx10 at12pt
 at12pt
\font\smit=cmmi10 at7pt

\centerline {\bf DIFF\'ERENTIELLES NON-COMMUTATIVES}
\medskip\centerline {\bf ET TH\'EORIE DE GALOIS DIFF\'ERENTIELLE}
\medskip\centerline {\bf OU AUX DIFF\'ERENCES.}
\bigskip
\centerline { par}
\bigskip
\centerline { \bf Yves Andr\'e} 
\bigskip
\bigskip
\centerline {\sm{Institut de Math\'ematiques, 175 rue du Chevaleret}}
\centerline {\sm{7\`eme \'etage}}
\centerline {\sm{F-75013 Paris}}
\centerline {\sm{andre@math.jussieu.fr}}
\bigskip
\bigskip
\bigskip\bigskip
\bigskip
Table des mati\`eres:
\bigskip Introduction.
\medskip {\bf I. Calcul diff\'erentiel, connexions et groupes de Galois.
 Cinq situations concr\`etes.} 
\par 1. En g\'eom\'etrie diff\'erentielle classique.
\par 2. En caract\'eristique $p$.
\par 3. Th\'eorie de Picard-Vessiot.
\par 4. Calcul aux diff\'erences.
\par 5. Aper\c cu sur le calcul diff\'erentiel quantique d'A. Connes.
\medskip {\bf II. Calcul diff\'erentiel non-commutatif et connexions.}     
\par 1. Alg\`ebres diff\'erentielles gradu\'ees et anneaux diff\'erentiels g\'en\'eralis\'es.
\par 2. Connexions.
\par 3. Produit tensoriel.
\par 4. La situation semi-classique.
\par 5. Localisation et rigidit\'e.
\medskip {\bf III. Groupes de Galois diff\'erentiels et extensions de
Picard-Vessiot en situation semi-classique.}     
\par 1. Solubilit\'e.
\par 2. Groupes de Galois diff\'erentiels.
\par 3. Le th\'eor\`eme de sp\'ecialisation.
\par 4. Extensions de Picard-Vessiot.
\par 5. Correspondance galoisienne.
\bigskip
\bigskip \medskip \rightline{\it 20/4/2000, \`a para\^{\i}tre dans Annales ENS.}

\vfill\eject 
\centerline{Introduction.} 
\bigskip La th\'eorie de Galois diff\'erentielle conna\^{\i}t un nouvel essor, devenant
source d'inspiration et d'applications dans des domaines de plus en plus vari\'es: nombres transcendants,
m\'ecanique c\'eleste, calcul formel, rigidit\'e de rev\^etements, sommes exponentielles... Par
ailleurs, la th\'eorie des groupes quantiques, la $q$-combinatoire et le d\'eveloppement de l'analyse
des \'equations aux diff\'erences suscitent un renouveau d'int\'er\^et pour la th\'eorie de Galois
aux ($q$-)diff\'erences. En amont, ce foisonnement d'ouvertures fait na\^{\i}tre un besoin de r\'enovation
des fondements de la th\'eorie, comme l'illustre le probl\`eme concret suivant.
\medskip Consid\'erons la s\'erie $q$-hyperg\'eom\'etrique de Heine 
\medskip \centerline{$\sum_{n\geq 0}\;{(q^{\alpha_1};q)_n \ldots
(q^{\alpha_r};q)_n\over (q^{\beta_1};q)_n\ldots (q^{\beta_{r-1}};q)_n.(q;q)_n}\;z^n\;.$}
\medskip \noindent Elle satisfait une \'equation aux $q$-diff\'erences d'ordre $r$. Lorsque $q$ tend vers
$1$, cette \'equation ``conflue" vers l'\'equation diff\'erentielle ordinaire satisfaite par la s\'erie
hyperg\'eom\'etrique de Gauss-Barnes $\;_rF_{r-1}(\matrix{\alpha_1&\ldots, &\alpha_r \cr \beta_1 &\ldots &
\beta_{r-1}}; z).$ On peut se poser les questions suivantes:
\par \noindent $i)$ comment comprendre la confluence du point de vue galoisien?  
\par \noindent $ii)$ A cette d\'eformation d'\'equation diff\'erentielle correspond-il de
mani\`ere naturelle des groupes quantiques comme $q$-d\'eformation du groupe de
Galois diff\'erentiel hyperg\'eom\'etrique?
\par Aborder ce type de questions requiert clairement une description uniforme des \'equations aux
($q$-)diff\'erences et des \'equations diff\'erentielles.
\medskip Notre point de vue sera celui des connexions. Rappelons qu'en alg\`ebre, une connexion sur un
$A$-module est une application additive $\nabla: M\rightarrow \Omega^1\otimes_A M$ v\'erifiant la r\`egle
de Leibniz $\nabla (am)= da\otimes m+a\nabla (m)$, o\`u $d: A \rightarrow \Omega^1$ est une d\'erivation
convenable. La recette pour traduire \'equations diff\'erentielles lin\'eaires en termes de
connexions est bien connue. Il s'av\`ere tout aussi possible d'interpr\'eter \'equations aux diff\'erences
lin\'eaires en termes de connexions, en consid\'erant des $A$-$A$-bimodules de diff\'erentielles
$\Omega^1$ non-commutatifs, i.e. o\`u l'action de $A$ n'est pas
la m\^eme \`a gauche et \`a droite. Cela fournit un cadre commode pour \'etudier des probl\`emes
mixtes diff\'erentiels-aux diff\'erences comme ci-dessus. Il y a donc lieu de se placer d'embl\'ee dans le
cadre du calcul diff\'erentiel non-commutatif, m\^eme si nos applications se limitent au cas
``semi-classique" o\`u l'anneau de fonctions $A$ est commutatif.  
\par Un autre avantage du point de vue des connexions est de faire le lien avec la g\'eom\'etrie
diff\'erentielle, et de prendre en compte certaines analogies frappantes entre la th\'eorie de l'holonomie
et la th\'eorie de Picard-Vessiot. En particulier, la possibilit\'e de consid\'erer des connexions \`a
courbure quelconque sugg\`ere la question:
\par \noindent $iii)$ peut-on consid\'erer un groupe d'holonomie comme groupe de Galois
diff\'erentiel?  
\medskip Pour b\^atir une th\'eorie du groupe de Galois
diff\'erentiel suffisamment souple, nous aurons recours \`a la th\'eorie tannakienne sur les anneaux (et \`a
un usage massif de la notion de platitude). 
\par La premi\`ere t\^ache est de d\'efinir le produit tensoriel de deux connexions, l'anneau de
base $A$ \'etant suppos\'e commutatif. Lorsque le bimodule $\Omega^1$ est commutatif, la r\`egle est bien
connue: $\nabla(m_1\otimes m_2)=\nabla(m_1)\otimes m_2+ m_1\otimes\nabla(m_2)$, o\`u l'on permute
tacitement $\Omega^1$ et $M_1$ dans le second terme. Pour \'etendre cette r\`egle au cas d'un bimodule
$\Omega^1$ non commutatif, on a donc besoin d'une application ``volte" $M_1\otimes \Omega^1\rightarrow
\Omega^1\otimes M_1$. Notre point de d\'epart est l'existence d'une telle {\it volte bilin\'eaire
canonique} (II.4.1.1). Elle permet de munir la cat\'egorie des connexions d'une structure
mono\"{\i}dale. En outre, l'\'echange des facteurs $M_1\otimes M_2\rightarrow
M_2\otimes M_1$ est une contrainte de commutativit\'e sym\'etrique. Cela donne une r\'eponse
n\'egative \`a la question $ii)$ ci-dessus (si, en modifiant la contrainte de commutativit\'e, on
trouvait un groupe quantique, il y aurait lieu d'apr\`es [Bru94] de le consid\'erer comme classique
puisque sa cat\'egorie de repr\'esentations le serait). 
\par On peut alors donner des crit\`eres g\'en\'eraux qui garantissent que cette cat\'egorie est
tannakienne (II.5.3.2, III.2.1.3). Un des corollaires est que {\it si $X$ est une vari\'et\'e alg\'ebrique
lisse sur un corps de caract\'eristique nulle ou une vari\'et\'e analytique lisse (r\'eelle, complexe 
ou $p$-adique rigide), alors tout ${\cal O}_X$-module coh\'erent muni d'une connexion} - non
n\'ecessairement int\'egrable - {\it est localement libre} (II.5.2.2). On donne aussi une
r\'eponse essentiellement positive \`a  la question $iii)$. 
\medskip Dans une situation \`a param\`etres comme ci-dessus (param\`etre $q$), le groupe de Galois
diff\'erentiel obtenu est un sch\'ema en groupes affine plat sur l'anneau des
constantes, qui commute au passage aux fibres (mais g\'en\'eralement pas de type fini). Ceci fournit
une r\'eponse, dans un cadre g\'en\'eral, \`a la question $i)$. On peut alors d\'eduire de
la connaissance des groupes de Galois diff\'erentiels hyperg\'eom\'etriques des informations sur les
groupes de Galois aux $q$-diff\'erences hyperg\'eom\'etriques pour $q$ g\'en\'erique (III.3.3).
\medskip On d\'eveloppe ensuite la ``th\'eorie de Picard-Vessiot" sans aucune hypoth\`ese
d'int\'egrabilit\'e (i.e. d'annulation de courbure). Elle montre, en termes figur\'es, qu'on peut toujours
int\'egrer {\it symboliquement} un syst\`eme diff\'erentiel ou aux diff\'erences lin\'eaire, m\^eme non
int\'egrable.     

\bigskip Passons en revue les diff\'erentes parties de ce travail. 
\medskip \noindent Le premier chapitre expose cinq situations concr\`etes qui ont motiv\'e
la th\'eorie, et l'illustrent. Les seuls r\'esultats originaux se trouvent en I.4. On y explique comment
interpr\'eter un endomorphisme semi-lin\'eaire comme connexion non-commutative, et comment la recherche de
solutions d'\'equations aux diff\'erences se ram\`ene \`a celle de sections horizontales de
telles connexions. Cette interpr\'etation a \'et\'e sugg\'er\'ee par un aspect du calcul
diff\'erentiel quantique d'A. Connes esquiss\'e en I.5.             
\par \noindent Dans I.1, on rappelle quelques points de g\'eom\'etrie diff\'erentielle classique:
th\'eor\`eme de Frobenius, connexions sur des fibr\'es principaux ou vectoriels et groupes d'holonomie, en
insistant sur l'alg\`ebre du calcul diff\'erentiel ext\'erieur. Outre le r\^ole d'exemples pour les
groupes de Galois diff\'erentiels d\'efinis ult\'erieurement, ces groupes d'holonomie fourniront une
interpr\'etation g\'eom\'etrique suggestive des constructions alg\'ebriques abstraites de la th\'eorie
(III.4.2.2).     
\par \noindent En I.2, on traite dans le m\^eme esprit la situation alg\'ebrique en caract\'eristique non
nulle: th\'eor\`eme de Cartier-Ekedahl, calcul diff\'erentiel \`a puissances divis\'ees...    
\par \noindent Le paragraphe I.3 rappelle bri\`evement la th\'eorie de Picard-Vessiot-Kolchin et sa
reformulation tannakienne. 
\medskip \noindent Le second chapitre traite du calcul diff\'erentiel non-commutatif et des connexions, et
r\'esout le probl\`eme du produit tensoriel. L'alg\`ebre diff\'erentielle ext\'erieure est remplac\'ee par
une alg\`ebre diff\'erentielle gradu\'ee $\Omega^{\ast}$ arbitraire. En tronquant en degr\'e $\leq 1$, on
obtient une d\'erivation $d$ de $A=\Omega^0$ \`a valeurs dans le $A$-$A$-bimodule $\Omega^1$. Ceci fournit
une vaste g\'en\'eralisation de la notion classique d'anneau diff\'erentiel, et un cadre commun aux
situations de I. On d\'efinit dans ce contexte les notions d'id\'eal diff\'erentiel, d'extension
diff\'erentielle, d'anneau diff\'erentiel simple...         
\par \noindent En II.2, on traite des connexions non-commutatives. On y examine particuli\`erement le cas
o\`u le module sous-jacent $M$ est un $A$-$A$-bimodule et o\`u la connexion $\nabla$ v\'erifie la r\`egle
de Leibniz \`a gauche. Motiv\'e par les travaux de J. Mourad [Mou95], M. Dubois-Violette et T. Masson
[DubM96] en g\'eom\'etrie non-commutative, on s'int\'eresse au cas o\`u le d\'efaut de suj\'etion
\`a la r\`egle de Leibniz \`a droite est donn\'e par un homomorphisme de bimodule $\phi(\nabla):M\otimes
\Omega^1 \rightarrow \Omega^1\otimes M$. On dira dans ce cas que $\nabla$ est une biconnexion, de
volte $\phi(\nabla)$. On calcule cette volte dans le cas particulier non trivial des \'equations aux
diff\'erences.   
\par \noindent En II.3, on d\'efinit le produit tensoriel de deux biconnexions: $\nabla =
\nabla_1\otimes id_2+({\phi}(\nabla_1)\otimes id_2)\circ(id_1 \otimes \nabla_2)$. On obtient ainsi une
structure mono\"{\i}dale sur la cat\'egorie des biconnexions, et on discute la dualit\'e. Dans ce
contexte g\'en\'eral, on peut s'attendre \`a ce que certaines sous-cat\'egories soient
$\otimes$-\'equivalentes \`a des cat\'egories de repr\'esentations de groupes quantiques non
classiques.       
\par \noindent A partir de II.4, on se place dans la situation ``semi-classique" o\`u $A$ est
commutatif mais o\`u le bimodule $\Omega^1$ n'est pas n\'ecessairement commutatif. On prouve qu'alors
toute connexion est une biconnexion, et que l'\'echange na\"{\i}f des facteurs est une contrainte de
commutativit\'e (sym\'etrique). Le paragraphe le plus substantiel de ce chapitre est II.5: on y
donne un crit\`ere sur $(A,d)$ (dont l'hypoth\`ese principale est la simplicit\'e) pour que la
cat\'egorie des modules \`a connexion de type fini sur $A$ \`a volte inversible soit tannakienne sur
le corps des constantes.      

\medskip Le troisi\`eme chapitre traite des aspects galoisiens dans la situation semi-classique. Il
commence par une \'etude d\'etaill\'ee de la notion de solubilit\'e d'une connexion dans une extension
diff\'erentielle, nourrie de nombreux exemples. On y montre en particulier que solubilit\'e et
int\'egrabilit\'e (i.e. nullit\'e de la courbure) sont des notions ind\'ependantes. 
\par \noindent En II.2, on consid\`ere, pour un module \`a connexion ${\cal M}=(M,\nabla)$ projectif de
type fini sur $A$ \`a volte inversible, la cat\'egorie ab\'elienne mono\"{\i}dale ${\scriptstyle <}{\cal
M}{\scriptstyle >}^{\scriptstyle\otimes}$ form\'ee des sous-quotients des sommes de tenseurs mixtes sur
${\cal M}$; puis on introduit de mani\`ere semi-axiomatique le groupe de Galois diff\'erentiel $Gal({\cal
M},\omega)$ de ${\cal M}$ attach\'e \`a un ``foncteur fibre" $\omega$. C'est un sch\'ema en groupes affine
et plat sur l'anneau des constantes $k$, de type fini si $k$ est un corps, mais pas en g\'en\'eral.
L'exemple principal est celui o\`u $\omega$ est le foncteur solution dans une extension diff\'erentielle
convenable d'anneau de constantes $k$. On introduit aussi le torseur des solutions, et le groupe de Galois
diff\'erentiel intrins\`eque (correspondant au foncteur oubli de la connexion), qu'on d\'ecrit en
caract\'eristique $p$ lorsque $\Omega^1$ est commutatif. On \'etablit ensuite diverses fonctorialit\'es
des groupes de Galois diff\'erentiels; on en d\'eduit, en II.3, le th\'eor\`eme de sp\'ecialisation, qui
\'eclaire les ph\'enom\`enes de confluence, en particulier la confluence d'\'equations aux
($q$-)diff\'erences vers une \'equation diff\'erentielle.      
\par \noindent Les deux derniers paragraphes g\'en\'eralisent la th\'eorie de Picard-Vessiot. 
Une extension diff\'erentielle $(A',d')$ de $(A,d)$ est dite de Picard-Vessiot pour $\cal
M$ si $A'$ est fid\`element plate sur $A$, $(A',d')$ est simple de corps de constantes $k$, ${\cal M}$ est
soluble dans $(A',d')$, et $A'$ est engendr\'ee comme $A$-alg\`ebre par $\langle M, (\check{M}\otimes
A')^{\nabla}\rangle$ et $\langle \check{M},(M\otimes A')^{\nabla}\rangle$. Le th\'eor\`eme
III.4.2.3 \'etablit une \'equivalence entre extensions de Picard-Vessiot et foncteurs fibres sur
${\scriptstyle <}{\cal M}{\scriptstyle >}^{\scriptstyle\otimes}$. On donne ensuite des crit\`eres
d'existence et d'unicit\'e pour les extensions de Picard-Vessiot (elles apparaissent comme alg\`ebres de
fonctions des torseurs de solutions).  
\par \noindent Etant donn\'ee une extension de Picard-Vessiot $(A',d')/(A,d)$ pour ${\cal M}$, on
identifie en III.5.1 le groupe des automorphismes de $(A',d')/(A,d)$ au groupe de Galois diff\'erentiel
de ${\cal M}$ attach\'e au foncteur solution dans $(A',d')$. 
\par \noindent Pour formuler la correspondance galoisienne, il y a lieu de supposer $A$
semi-simple (i.e. produit fini de corps) et d'introduire la notion voisine et ``birationnelle" d'extension
de Picard-Vessiot fractionnaire; il s'agit essentiellement de l'anneau total des fractions d'une extension
de Picard-Vessiot lorsque cet anneau est semi-simple. Si le corps des constantes est alg\'ebriquement clos
de caract\'eristique nulle, la correspondance galoisienne pour les extensions
de Picard-Vessiot fractionnaires a lieu sans surprise (III.5.2). Elle unifie la correspondance de
Picard-Vessiot diff\'erentielle classique et son analogue aux diff\'erences, et englobe le cas de
syst\`emes mixtes \`a plusieurs variables, \'eventuellement non int\'egrables.       

\bigskip \noindent {\smitint{Remerciements}}. {\smint{Cet article est l'aboutissement d'un projet
commenc\'e avec [An87],[An89], et dont diverses phases ont \'et\'e expos\'ees au S\'eminaire d'Alg\`ebre de
Paris 86, au Colloque diff\'erentiel de Plovdiv 93, au Colloque Galois diff\'erentiel de Luminy
99 et \`a l'Atelier diff\'erentiel de Strasbourg 99; je remercie les organisateurs de m'en
avoir donn\'e l'occasion. Je remercie M. Dubois-Violette et C. Kassel de m'avoir fait
conna\^{\i}tre les articles [DubM96] et [N97] respectivement.}}
 \vfill \eject

\centerline {{\big{ {\bf \S { I}. Calcul diff\'erentiel, connexions et groupes de Galois.}}}}
\centerline {{\big{ {\bf Cinq situations concr\`etes.}}}}

\bigskip {\bf 1. En g\'eom\'etrie diff\'erentielle classique.} 
\medskip {\bf 1.1. Calcul diff\'erentiel ext\'erieur de
Cartan-De Rham-K\"ahler et th\'eor\`eme de Frobenius.} 
\medskip \noindent Soient $X$ une vari\'et\'e diff\'erentiable $C^{\infty}$ de dimension $m$,
$T(X)$ son fibr\'e tangent, $\Lambda^{\ast}(X)$ l'alg\`ebre ext\'erieure sur le fibr\'e cotangent
$T^{\vee}(X)$. L'existence du crochet de Lie $[\;,\;]$ sur $T(X)$ permet de d\'efinir la
diff\'erentiation ext\'erieure $d: \Lambda^p(X)\rightarrow \Lambda^{p+1}(X)$ par la formule bien
connue:
$$d\omega(D_1,..., D_{p+1})=\sum_1^{p+1}(-)^iD_i(\omega(...,\hat{D_i},...))+ 
\sum_{i< j}(-)^{i+j} \omega([D_i,D_j],...\hat{D_i},...\hat{D_j},...D_{p+1}). $$
Gr\^ace \`a l'identit\'e de Jacobi satisfaite par $[\;,\;]$, on a $d^2=0$, ce qui fait de
$\Lambda^{\ast}(X)$ un fibr\'e en alg\`ebres diff\'erentielles gradu\'ees (commutatives au sens
gradu\'e). Le lemme de Poincar\'e dit qu'au niveau des faisceaux de sections, le complexe obtenu
est acyclique en degr\'e $>0$.
\par Un champ de $p$-plans $\Sigma$ sur $X$ est dit {\it involutif} si chaque
fibre $\Sigma_x$ est stable sous $[\;,\;]$. Sous cette condition, le th\'eor\`eme de
Frobenius affirme que $\Sigma$ est le fibr\'e tangent d'un feuilletage: par tout point $x\in X$
passe une (unique) vari\'et\'e int\'egrale connexe maximale $X_{\Sigma, x}$ pour $\Sigma$ (de
dimension \'egale \`a $p$).   
\par \noindent Traduction en calcul diff\'erentiel ext\'erieur: l'application qui associe
\`a $\Sigma$ son orthogonal $\Sigma^{\bot}$ dans $\Lambda^{\ast}(X)$ est une bijection entre
champs de $p$-plans et id\'eaux de $\Lambda^{\ast}(X)$ localement engendr\'es par
$m-p\;1$-formes ind\'ependantes. De plus, $\Sigma^{\bot}$ est stable sous $\;d\;$ si et seulement
si $\Sigma$ est involutif (cf. [W71] p.73). La condition que $X_{\Sigma, x}$ est une vari\'et\'e
int\'egrale se traduit par $(\Sigma^{\bot})_{\mid X_{\Sigma, x}}=0$.

\medskip {\bf 1.2. Connexions et holonomie.} 
\medskip \noindent Soient $G$ un groupe de Lie et $\pi : P\rightarrow X$ un fibr\'e principal
\`a droite sous $G$. On note $t_g$ la translation par $g\in G$ sur $P$. Pour tout $p\in P$, on
note $V_p(P)$ l'espace (``vertical") des vecteurs tangents \`a la fibre en $p$, c'est-\`a-dire
$Ker\;\pi_{\ast p}$. Rappelons qu'une connexion $\aleph$ sur $P$ est la donn\'ee pour tout
$p\in P$, d'un suppl\'ementaire (``horizontal") $H_p(P)$ de $V_p(P)$ dans $T_p(P)$, variant de
mani\`ere $C^{\infty}$ avec $p$, tel que ${t_g}_{\ast}(H_p(P))=H_{pg}(P)$. Elle induit un
isomorphisme $H_p(P)\cong T_{\pi(p)}(X)$. Si $C$ est une courbe $C^{\infty}$ par morceaux sur
$X$, et $x$ un point de cette courbe, la connexion permet de relever $C$ en une courbe
``horizontale" ${\tilde C}_p$ sur $P$ passant par un point $p \in P$ fix\'e au-dessus de $x$: tout
vecteur tangent \`a
${\tilde C}_p$ est horizontal (on dit aussi ``parall\`ele"). On a ${\tilde
C}_{pg}=t_g({\tilde C}_p)$.
\par \noindent 
Si $C$ est un {\it lacet} bas\'e en $x = \pi(p)$, l'extr\'emit\'e de ${\tilde C}_p$ est un point
de la fibre $\pi^{-1}(x)$; c'est donc le translat\'e $t_{g(C)}(p)$ de $p$ par un \'el\'ement bien
d\'efini $g(C)$ de $G$. Lorsque $C$ varie, les
$g(C)$ forment un sous-groupe de Lie $Hol_p(\aleph)$ de $G$, le groupe d'{\it holonomie} point\'e
en $p$. Sa composante neutre $Hol^0_p(\aleph)$ est form\'e des $g(C)$ avec $C$ contractile.
Lorsque $p$ varie, les $Hol_p(\aleph)$ sont conjugu\'es dans $G$, si $X$ est connexe. On montre
que le fibr\'e principal $\pi : P\rightarrow X$ provient en fait d'un fibr\'e principal $\pi' :
P'\rightarrow X$ de groupe $Hol_p(\aleph)$ (cf. [Li55], [KN63]). 
\par Soit {\goth g} l'alg\`ebre de Lie de $G$. Pour tout $a\in
\;${\goth g}, notons $D_{\scriptstyle a}$ le champ de vecteurs de $P$ g\'en\'erateur infinit\'esimal du
groupe
\`a un param\`etre $s \mapsto t_{exp\;sa}$. C'est un champ de vecteur vertical, et lorsque $a$
varie, $(D_{\scriptstyle a})_p$ engendre $V_p(P)$. On d\'efinit une $1$-forme $\omega$ \`a
valeurs dans {\goth g} par: 
\par  pour tout $a \in \;$ {\goth g},
$\;\omega_p((D_{\scriptstyle a})_p)=a\;;\;\; \omega_p(H_p(P))=0.$ 
\par \noindent Elle v\'erifie $t_g^{\ast}(\omega)=ad(g^{-1})\omega$. Dans tout ouvert $U$
trivialisant $P$ ($P_{\mid U}\cong U\times G$), elle s'\'ecrit sous la forme $\omega_{\mid U}=
g^{-1}dg - ad(g^{-1})\pi^{\ast}\theta_U$, o\`u $\theta_U$ est une $1$-forme \`a valeurs dans
{\goth g} sur $U$. Il est clair que la donn\'ee de $\omega$ \'equivaut \`a celle de $\aleph$. 
\par \noindent La $2$-forme de {\it courbure} $\Omega$ (\`a
valeurs dans {\goth g}) est caract\'eris\'ee par 
\par$\Omega(D, D') = 0\;$  si $D$ ou $D'$ est vertical,
\par $\Omega( D, D')=\omega([ D, D'])\;$ si $D$ et $D'$ sont horizontaux. 
\par \noindent On voit donc que la courbure est identiquement nulle si et seulement si
pour tout $p$, $H_p(P)$ est un champ de $m$-plans involutif sur $P$, ou encore, d'apr\`es le
th\'eor\`eme de Frobenius, si et seulement si $\aleph$ est {\it plate}, i.e. $(P,\aleph)$ est
localement du type produit $(X \times G, \;pr_1^{\ast}T(X))$. Ceci \'equivaut aussi \`a la
trivialit\'e de $Hol^0_p(\aleph)$. 
\par \noindent Plus g\'en\'eralement, le th\'eor\`eme d'Ambrose-Singer affirme que
$Lie\;Hol^0_p(\aleph)$ est le sous-espace de {\goth g} engendr\'e par les $\Omega_{q}(D,D')$,
o\`u $q$ d\'ecrit tous les points de $P$ qu'on peut joindre \`a $p$ par une courbe horizontale
({\it loc. cit.}). 
\medskip {\bf 1.3. Connexions sur un fibr\'e vectoriel.}
\medskip \noindent Lorsque $G= GL_n$, une connexion
sur $P$ induit une ``connexion" sur le fibr\'e vectoriel standard $E$ sur $X$ associ\'e \`a $P$,
c'est-\`a-dire un op\'erateur additif 
$$\nabla: \Gamma(E) \rightarrow \Gamma(T^{\vee}(X)\otimes E)$$ v\'erifiant la r\`egle de Leibniz.
Dans tout ouvert $U$ trivialisant $P$ (donc aussi $E$), il est donn\'e par  $\nabla(e_{\mid U})=
de_{\mid U} - \theta_U \otimes e_{\mid U}$. On peut voir $\nabla$ comme une r\`egle associant \`a
tout champ de vecteur $D$ sur $X$ un endomorphisme $\nabla_{\scriptstyle D}$ de $E$. Il est clair que la
donn\'ee de $\nabla$ \'equivaut \`a celle de $\aleph$. 
\par Le {\it principe d'holonomie} affirme qu'un champ d'objets parall\`eles d\'efinit en chaque
point un objet invariant par le groupe d'holonomie et r\'eciproquement. Par exemple si
$X$ poss\`ede une m\'etrique riemannienne $g$, et $E= T(X)$ muni de la connexion de Levi-Civita, 
$Hol_p(\aleph)$ est un sous-groupe du groupe orthogonal $O(T_{\pi(p)}(X),g)$ ($Hol^0_p(\aleph)$
est un sous-groupe compact de $SO(T_{\pi(p)}(X),g)\;$) (cf. [Bry99] pour un r\'ecent survol).  
\medskip La {\it courbure} de $\nabla$ est la $2$-forme sur $X$ \`a valeurs dans $End\;E$
d\'efinie par $R_{{\scriptstyle \nabla}}(D,D') = [\nabla_{\scriptstyle D} ,\nabla_{D'}]-\nabla_{[D,D']}$. Elle s'annule si et
seulement si $\aleph$ est plate; on dit alors que $\nabla$ est {\it plate} ou {\it
int\'egrable}. Il revient au m\^eme de demander que l'action de $Hol_p(\nabla)$ sur $E_x$
se factorise \`a travers $\pi_1(X,x)$ (repr\'esentation de {\it monodromie} en $x=\pi(p)$).

\bigskip {\bf 2. En caract\'eristique $p$.} 
\medskip {\bf 2.1. Calcul diff\'erentiel ext\'erieur, puissance $p$-i\`eme, et variantes du
th\'eor\`eme de Frobenius.} 
\medskip \noindent Soient $X$ une vari\'et\'e affine lisse de dimension $m$ sur un corps $k$
parfait de caract\'eristique $p>0$, $T(X)$ l'alg\`ebre de Lie des champs de vecteurs,
$\Omega^{\ast}(X)$ l'alg\`ebre ext\'erieure sur $T^{\vee}(X)$ (formes de K\"ahler). $T(X)$
est munie d'une $p$-structure (l'application puissance $p$-i\`eme). Cette structure
suppl\'ementaire se refl\`ete en calcul diff\'erentiel ext\'erieur par l'existence de
l'{\it op\'eration de Cartier}: un isomorphisme d'anneaux gradu\'es $C:H^{\ast}_{DR}(X):=
H^{\ast}(\Omega^{\ast}(X)) \rightarrow \Omega^{\ast}(X)$, d'inverse d\'efini par $C^{-1}(f) =
[f^p], \; C^{-1}(df)=[f^{p-1}df]$.    
\par L'analogue du th\'eor\`eme de Frobenius en caract\'eristique $p$ est probl\'ematique (les
feuilletages ne sont pas d\'etermin\'es par leur partie d'ordre un...) On peut toutefois obtenir
un \'enonc\'e net en faisant intervenir la $p$-structure ([E87] 2.4):  

\proclaim Proposition 2.1.1. Il y a une correspondance bijective entre sous-fibr\'es
involutifs $\Sigma$ de $T(X)$ stables par puissance $p$-i\`eme, et morphismes finis plats $h:
X\rightarrow X'$ de hauteur $1$, donn\'ee par $\Sigma = T(X/X')$ (fibr\'e tangent relatif).
\par \noindent 
Soit $E$ un fibr\'e vectoriel sur $X$ muni d'une connexion alg\'ebrique $\nabla$.
Outre la courbure $R_{{\scriptstyle \nabla}}$, il y a une autre obstruction \`a l'existence de sections
horizontales: la $p$-{\it courbure} $R_{{\scriptstyle \nabla},p}$, d\'efinie par $R_{{\scriptstyle \nabla}, p}(D) = (\nabla_{\scriptstyle D})^p
-\nabla_{(D^p)}$. Ce sont les seules obstructions, d'apr\`es le r\'esultat suivant de Cartier:

\proclaim Proposition 2.1.2. Soit $F: X \rightarrow X' = X\times_{k, x\mapsto x^p}k\;$ le
Frobenius relatif (\'el\'evation des coordonn\'ees de $X$ \`a la puissance $p$-i\`eme). Il y a une
\'equivalence de cat\'egories entre fibr\'es $E$ sur $X$ munis d'une connexion int\'egrable
$\nabla$ \`a $p$-courbure nulle, et fibr\'es $E'$ sur $X'$, donn\'ee par $(E,\nabla)\mapsto E'=
F_{\ast}(E^{\nabla}),\;E=F^{\ast}(E')$.
\par Voir [Kat70] 5.1.1. Le lien avec 2.1.1 est le suivant: soient $P'$ un
fibr\'e principal sous $GL_n$ associ\'e \`a $E'$, $P$ son image inverse sur $X$ via $F$, et $h:
P \rightarrow P'$ le morphisme fini plat induit par $F$. La connexion $\nabla$ provient
d'un sous-fibr\'e de $T(P)$, involutif et stable par puissance $p$-i\`eme.   

\medskip {\bf 2.2. Puissances divis\'ees et lemme de Poincar\'e.} 
\medskip \noindent Pour disposer d'un analogue du lemme de Poincar\'e formel en
caract\'eristique $p$, on est conduit \`a introduire des anneaux de s\'eries de puissances
divis\'ees (cf. [Bou81] IV 5, [BeO78] app.). Soient $x$ un point $k$-rationnel de $X$, et
$t_1,\ldots,t_m$ des coordonn\'ees locales en $x$. Le compl\'et\'e formel $\hat{\cal O}_x$ de
l'anneau local en $x$ s'identifie \`a
$k[[t_1,\ldots,t_m]]\cong {S(T_x(X))}^{\wedge}$, o\`u $S$ d\'esigne l'alg\`ebre sym\'etrique. Le
compl\'et\'e \`a puissances divis\'ees ${\cal O}^{pd}_x$ de ${\cal O}_x$ est l'alg\`ebre
topologique $(S(T_x(X)^{\vee}))^{\vee}$.         
\par On lui associe un complexe de De Rham ${\cal O}^{pd}_x \;{\buildrel{d}\over
\longrightarrow}\; ({\Omega}^1)^{pd}_x \;{\buildrel{d}\over
\longrightarrow}\; ({\Omega}^2)^{pd}_x \;{\buildrel{d}\over
\longrightarrow}\; \ldots  $ qui est acyclique en degr\'e $> 0$. 
L'homomorphisme naturel ${\cal O}_x \rightarrow {\cal O}^{pd}_x$ a pour noyau l'id\'eal
engendr\'e par les $t_i^p$. On note $\bar{\cal O}_x \hookrightarrow {\cal O}^{pd}_x$ l'anneau
quotient de ${\cal O}_x$, et $\bar X_x = Spec\;\bar{\cal O}_x$. Bien que $\bar X_x$ soit
non-r\'eduit, le calcul diff\'erentiel ext\'erieur se comporte bien: $\Omega^{\ast}(\bar X_x)\cong
\bar{\cal O}_x\otimes_{\scriptstyle k} \Lambda^{\ast}(T^{\vee}_x(X)),\; H^{\ast}_{DR}(\bar X_x)\cong
\Lambda^{\ast}(\oplus k.t_i^{p-1}dt_i)$, et $\Omega^{pair}(\bar X_x)$ admet des puissances
divis\'ees. 
\par $T(\bar X_x)$ est une $p$-alg\`ebre de Lie simple sur $k$ non classique (analogue en
caract\'eristique $p$ de l'alg\`ebre des champs de vecteurs sur un tore); la sous-alg\`ebre de
Lie form\'ee des champs de vecteurs qui respectent l'id\'eal maximal s'identifie \`a
$Lie\;Aut(\bar X_x)$, cf. [Mat99].     
\proclaim Proposition 2.2.1. Toute connexion int\'egrable $\nabla$ sur $\bar X_x$ est soluble dans
${\cal O}^{pd}_x$. Elle est soluble dans $\bar{\cal O}_x$ si et seulement si sa $p$-courbure est
nulle.
\par \noindent Ceci se d\'emontre par le ``t\'elescopage" habituel: pour toute section $e$ de $E$,
\medskip \centerline{$\Sigma_{(n_1,\ldots,n_m)}\;\Pi (-t_i)^{[n_i]}\;\Pi (\nabla({d\over
dt_i}))^{n_i}(e)$}
\par \noindent est une section horizontale de $E\otimes{\cal O}^{pd}_x\;$ ($\;^{[n]}$
d\'esigne une puissance divis\'ee); pour la seconde assertion, on note que $R_{{\scriptstyle \nabla},p}=0$ si et
seulement si tous les $(\nabla({d\over dt_i}))^{p}$ s'annulent.

\bigskip {\bf 3. Th\'eorie de Picard-Vessiot.} 
\medskip {\bf 3.1. R\'esum\'e.} 
\medskip \noindent Soit $K$ un corps diff\'erentiel de caract\'eristique $0$, muni d'une
d\'erivation $\partial$. On suppose le corps des constantes $k:=K^{\partial}$ alg\'ebriquement
clos (e.g. $K={\bf C}(z), \partial = {d\over dz}$). On consid\`ere une \'equation diff\'erentielle
lin\'eaire
\medskip \centerline{${(\ast)}\;\;\;\;\;\;\;\partial^{\mu}y + a_1 \partial^{\mu -1}y + \ldots +
a_{\mu}y = 0$}
\par \noindent \`a coefficients $a_i$ dans $K$.
\par \noindent Une extension diff\'erentielle $(L, \partial)$ de
$(K,\partial)$ est dite de {\it Picard-Vessiot} si 
\par $i)$ $\;k=L^{\partial}$, 
\par $ii)$ $L$ est engendr\'e par $\mu$ solutions $k$-ind\'ependantes de $(\ast)$ et leurs
d\'eriv\'ees.  
\par \noindent D'apr\`es Kolchin, il existe une extension de Picard-Vessiot, unique \`a
isomorphisme pr\`es. Le {\it groupe de Galois diff\'erentiel} $Gal(\ast):= Aut_{\partial}(L/K)$
est un sous-groupe alg\'ebrique du groupe des automorphismes du $k$-espace $Sol(\ast)$ des
solutions de $(\ast)$ dans $L$. On a une correspondance galoisienne entre sous-groupes ferm\'es
de $Gal(\ast)$ et extensions diff\'erentielles interm\'ediaires $K\subset F\subset L$. Voir
[Le90] pour un expos\'e concis et d\'ecant\'e de la th\'eorie.  

\medskip {\bf 3.2. Point de vue tannakien.} 
\medskip \noindent Consid\'erons l'anneau de polyn\^omes tordu $K[\partial]$, et posons $\Omega^1
= (K \partial)^{\vee}$. L'\'equation diff\'erentielle $(\ast)$ donne lieu \`a un
$K[\partial]$-module:
\par \centerline{ $Hom_{\scriptstyle k}(K[\partial]/K[\partial](\partial^{\mu} + a_1 \partial^{\mu -1} + \ldots +
a_{\mu}),\;K)$.} 
\par \noindent On note $M$ le $K$-espace de dimension $\mu$ sous-jacent. L'action de $\partial$
sur $M$ se d\'ecrit encore comme connexion: application additive $\nabla: M \rightarrow \Omega
\otimes_{\scriptstyle k} M$ v\'erifiant la r\`egle de Leibniz. 
\par \noindent L'espace des solutions de $(\ast)$ dans toute
extension diff\'erentielle $F$ de $K$ s'identifie \`a l'espace $(M\otimes F)^{\nabla}$ des
sections horizontales de $M\otimes F$.    
\par \noindent Le produit tensoriel d'espaces vectoriels \`a connexion est d\'efini par la r\`egle
habituelle $$\nabla(m_1\otimes m_2) = (\nabla(m_1))\otimes m_2 + m_1\otimes (\nabla(m_2)).$$ 
Le foncteur ``solutions dans $L$" \'etablit alors une \'equivalence de cat\'egories tensorielles
entre la cat\'egorie tensorielle engendr\'ee par $(M,\nabla)$ et la cat\'egorie des
repr\'esentations de dimension finie de $Gal(\ast)$ sur $k$ (cf. [De90]9, [An89[). Nous
renvoyons \`a [Se93] pour une introduction \'el\'ementaire \`a la dualit\'e de Tannaka
alg\'ebrique, et \`a
[Ber92] pour un expos\'e de techniques de calcul du groupe de Galois diff\'erentiel.

\bigskip
{\bf 4. Calcul aux diff\'erences.} 
\medskip {\bf 4.1. Diff\'erences finies et polyn\^omes tordus d'Ore.} 
\medskip \noindent Le calcul aux diff\'erences au sens large fait intervenir un endomorphisme
$\sigma$ d'un anneau de fonctions $A$, et s'int\'eresse \`a l'op\'erateur de diff\'erence
\medskip\centerline{$\delta_{\sigma}: a\mapsto \gamma(a^{\sigma}-a),$}
\medskip\noindent o\`u $\gamma$ est un
\'element convenable de $A$. La propri\'et\'e fondamentale d'un tel op\'erateur est d'\^etre une
{\it $\sigma$-d\'erivation}, i.e. d'\^etre additif et de v\'erifier la r\`egle 
\medskip\centerline{$\delta_{\sigma}(ab)=
\delta_{\sigma}(a)b+a^{\sigma}\delta_{\sigma}(b).$}   
\medskip\noindent{\it Exemples de $\sigma$-d\'erivations:} $\bullet$ $\sigma=$ translation de pas $h$,
 $(\delta_{\sigma}(f))(z)={f(z+h)-f(z)\over h}$,  
\par \noindent $\bullet$ $\sigma=$ dilatation de module $q$,
 $(\delta_{\sigma}(f))(z)={f(qz)-f(z)\over (q-1)z}$,
\par \noindent $\bullet$ $\sigma=$ identit\'e,  $\delta_{\sigma}(f)={df\over dz}$ (correspondant
au passage \`a la limite $h\rightarrow 0$ ou $q\rightarrow 1$).
\par La th\'eorie des $\sigma$-d\'erivations a \'et\'e d\'evelopp\'ee dans cette optique par Ore
et d'autres (cf. [Coh77]; citons aussi le calcul diff\'erentiel galoisien [He92] qui s'y
rattache). Dans cette approche, on consid\`ere l'anneau de polyn\^omes tordu $A[X]_{\sigma,
\delta}$ (avec la loi de commutation $Xa=a^{\sigma}X+\delta(a)$, o\`u $\delta$ est une
$\sigma$-d\'erivation). On associe alors \`a tout syst\`eme lin\'eaire aux diff\'erences
attach\'e \`a $\sigma$ un $A$-module (de type fini) $M$ muni d'endomorphisme $\theta$
pseudo-lin\'eaire (i.e. v\'erifiant la r\`egle $\theta(am)=
\delta(a)m+a^{\sigma}\theta(m)$), qu'on interpr\`ete comme un $(A[X]_{\sigma,
\delta})$-module; on tire alors parti, lorsque $A$ est un corps, du fait que $A[X]_{\sigma,
\delta}$ est {\it euclidien}, et que tout $(A[X]_{\sigma,\delta})$-{\it module de type fini est
somme directe de sous-modules cycliques} [Coh71]. C'est ainsi, par exemple, que Praagman
\'etablit, en prouvant un lemme de Hensel pour les polyn\^omes tordus, la classification formelle
des syst\`emes lin\'eaires aux diff\'erences [Pr83] (voir aussi [Duv83]).  
\par C'est un point de vue dual, en un certain sens, que nous adopterons: dans le cas
$\sigma=$ identit\'e, il s'agirait du calcul diff\'erentiel ext\'erieur de Cartan-
K\"ahler-Koszul; dans le cas g\'en\'eral, il s'agira de ses g\'en\'eralisations non-commutatives.
L'id\'ee de base est d'interpr\'eter les $\sigma$-d\'erivations comme de vraies d\'erivations \`a
valeurs dans des bimodules non-commutatifs.   

\medskip {\bf 4.2. $\sigma$-d\'erivations et sesquimodules.}
\medskip \noindent Soient $k$ un anneau commutatif unitaire et $A$ une $k$-alg\`ebre commutative
unitaire munie d'un $k$-endomorphisme $\sigma$. Rappelons qu'une $k$-d\'erivation $\delta$ de
$A$ \`a valeurs dans un $A$-$A$-bimodule $U$ est une application $k$-lin\'eaire v\'erifiant la
r\`egle de Leibniz 
\medskip\centerline{$\delta(ab)= \delta(a)b+a\delta(b).$}
\medskip\noindent Comme $A$ est suppos\'e commutatif, ces applications forment un $A$-$A$-bimodule
not\'e $Der_{\scriptstyle k}(A,U)$.
\par \noindent Soit $T$ un $A$-module \`a droite. Consid\'erons-le comme un $A$-$A$-bimodule
avec la r\`egle
\medskip\centerline{$\;a.t=t.a^{\sigma}{\rm \;\;pour\; tout \;\;}a\in A {\rm \;\;et\; tout}\;\;
t\in T.$}
\medskip\noindent Nous appellerons {\it sesquimodule} un tel bimodule. Les homomorphismes de
sesquimodules co\"{\i}ncident bien entendu avec les homomorphismes de $A$-modules \`a droite
sous-jacents, et forment eux-m\^emes des sesquimodules. 
\par \noindent En particulier, soit $A_{sesq}$ le sesquimodule provenant
du $A$-module \`a droite $A$. Alors $Der_{\scriptstyle k}(A,A_{sesq})$ n'est autre que le sesquimodule des
$\sigma$-d\'erivations $k$-lin\'eaires de $A$.  

\proclaim Proposition 4.2.1. Il existe un sesquimodule $\Omega^1_{\sigma}$ et une d\'erivation
$d\in Der_{\scriptstyle k}(A,\Omega^1_{\sigma})$ tels que pour tout sesquimodule $T$, l'application
$$Hom(\Omega^1_{\sigma},T)\rightarrow Der_{\scriptstyle k}(A,T): h\mapsto h\circ d$$
est un isomorphisme (de sesquimodules). En particulier, les $\sigma$-d\'erivations $k$-lin\'eaires
de $A$ s'identifient aux formes $A$-lin\'eaires (\`a droite) sur $\Omega^1_{\sigma}$.   
\par {\it Preuve}. Soit $I=Ker(A\otimes_{\scriptstyle k} A \rightarrow A)$ le noyau de la multiplication. On
sait que l'application $d: x\mapsto 1\otimes x -  x\otimes 1$ est une $k$-d\'erivation
de $A$ dans $I$, et que $I=A.dA$; de plus $f\mapsto f\circ d$ d\'efinit un isomorphisme
$Hom_{(A,A)}(I,T)\cong Der_{\scriptstyle k}(A,T)$ (cf.[Bou81] III,10, p.132). Comme
$A$ est commutatif, $I$ est un id\'eal de $A\otimes_{\scriptstyle k} A$. Utilisons maintenant le fait
que $T$ est un sesquimodule. Soit $J$ l'id\'eal de $A\otimes_{\scriptstyle k} A$
engendr\'e par les \'el\'ements $1\otimes x^{\sigma}-x\otimes 1$. Modulo
l'identification entre $A$-$A$-bimodules et $A\otimes_{\scriptstyle k} A$-modules, on a donc $JT=0$. Compte tenu
de l'identification de $I\otimes_{\scriptstyle k} (A\otimes_{\scriptstyle k} A/J)$ et de $I/JI$, on obtient
$Hom(I/JI,T)\cong Der_{\scriptstyle k}(A,T)$, et on voit que $\Omega^1_{\sigma}=I/IJ$ convient.       

\medskip \noindent {\bf 4.2.2. Exemples.} $i)\;$ $A= k[z], k(z), k[[z]]$ ou $k((z))$. Si
$z^{\sigma}\neq z$, $A$ est stable par $f\mapsto \delta_{\sigma}(f)={f^{\sigma}-f \over
z^{\sigma}-z}\;$; $\delta_{\sigma}$ est une base du
$A$-module des $\sigma$-d\'erivations $k$-lin\'eaires de $A$ (\`a remplacer par $d\over dz$ si
$z^{\sigma}=z.$). $\Omega^1_{\sigma}$ est le $A$-module libre de rang un $dz.A$, et la
d\'erivation $d: A
\rightarrow dz.A$ est donn\'ee par $d(f(z))=dz.\delta_{\sigma}(f)$.       
\medskip $ii)\;$ $k$ est un corps, $A=k[z_1,z_2], z_i^{\sigma}=q_iz_i,$ avec $q_i\in
k.$  En d\'eveloppant l'\'egalit\'e $d(z_1z_2)=d(z_2z_1)$, on trouve la relation
$\;dz_2(q_1-1)z_1=dz_1(q_2-1)z_2.$
\par \noindent Il y a lieu de distinguer trois cas:
\par $a)$ $q_1=q_2=1:\Omega^1_{\sigma}=dz_1A\oplus dz_2A$,  
\par $b)$ $q_1\neq 1, q_2=1:\Omega^1_{\sigma}=dz_1A$, 
\par $c)$ $q_1\neq 1, q_2\neq 1: \Omega^1_{\sigma}$ est libre de rang un sur $A$. Une base
$\delta_{\sigma}$ du $A$-module des $\sigma$-d\'erivations
$k$-lin\'eaires de $A$ est donn\'ee par
$\delta_{\sigma}(z_1)=(q_1-1)z_1,\delta_{\sigma}(z_2)=(q_2-1)z_2.$

\medskip {\bf 4.3. Syst\`emes lin\'eaires aux diff\'erences et connexions.}
\medskip \noindent {\bf 4.3.1.} Dans la situation de l'exemple 4.2.2 $i)$, consid\'erons un
syst\`eme lin\'eaire aux diff\'erences
\medskip \centerline{${(\ast\ast)}\;\;\;\;\;\;\;\sigma(Y)={\cal A}Y$}
\par \noindent o\`u ${\cal A} \in GL_{\mu}(A)$. 
\par \noindent Pour une pr\'esentation plus intrins\`eque, [vdPS97] introduit
l'anneau de polyn\^omes de Laurent tordu $A[X,X^{-1}]_{\sigma}\;$:
$Xa=a^{\sigma}X$ (anneau des op\'erateurs aux $\sigma$-diff\'erences), et associe
\`a ${(\ast\ast)}$ le $A[X,X^{-1}]_{\sigma}$-module $M=A^{\mu}$, o\`u $X$ agit via un
endomorphisme $\sigma$-lin\'eaire $\Phi$ repr\'esent\'e par ${\cal A}^{-1}$ dans la base
canonique $(m_1,\ldots,m_{\mu})$. Le syst\`eme ${(\ast\ast)}$ \'equivaut \`a dire que les $\sum_i
Y_{ij}m_i$ constituent une base form\'ee de points fixes de $\Phi$.
\par \noindent Supposons que ${{\cal A}-1\over
z^{\sigma}-z}$ soit \`a coefficients dans $A$. Le syst\`eme
$(\ast\ast)$ s'\'ecrit encore sous la forme
$$\delta_{\sigma}(Y)=-{{\cal A}^{-1}-1\over
z^{\sigma}-z}\;Y^{\sigma}.$$ En termes plus intrins\`eques, on munit $M$ de la connexion
$\nabla: M \rightarrow \Omega^1_{\sigma}\otimes_{\scriptstyle A} M$ d\'efinie par  
$$\nabla(m)=dz\otimes {\Phi-id\over z^{\sigma}-z}\;(m).$$ En vertu de la structure de
sesquimodule de $\Omega^1_{\sigma}$, elle v\'erifie la r\`egle de Leibniz
$$\nabla(am)=a\nabla(m)+da\otimes m.$$
R\'eciproquement, la donn\'ee de $\nabla$ permet de retrouver $\Phi$; le noyau de $\nabla$ est
constitu\'e des points fixes de $\Phi$.
\medskip Ce formalisme des connexions fournit non seulement un cadre alg\'ebrique unifi\'e
pour l'\'etude des \'equations diff\'erentielles et des \'equations aux diff\'erences, mais encore
un cadre alg\'ebrique commode pour \'etudier la ``confluence" des \'equations aux
$q$-diff\'erences vers une \'equation diff\'erentielle lorsque $q$ tend vers $1$: prendre par
exemple $k = {\bf C}[[q-1]], A = k[[z]][{1\over z}]$.   
\medskip \noindent {\bf 4.3.2.} Dans un cadre plus g\'en\'eral, soient $A$ une $k$-alg\`ebre commutative
unitaire munie d'un $k$-endomorphisme $\sigma$, $\gamma$ un \'el\'ement de $A$, $\delta_{\sigma,\gamma}$ la
$\sigma$-d\'erivation $a\mapsto \gamma(a^{\sigma}-a),$ vue comme un d\'erivation $A\rightarrow A_{sesq}$. 
\par \noindent Soit $M$ un $A$-module muni d'un endomorphisme $\sigma$-lin\'eaire $\Phi$. On lui
associe la connexion $\nabla: M \rightarrow A_{sesq}\otimes_{\scriptstyle A} M$ d\'efinie par  
$$\nabla(m)=1\otimes \gamma(\Phi-id )(m).$$ Elle v\'erifie la r\`egle de Leibniz
$$\nabla(am)=a\nabla(m)+\delta_{\sigma,\gamma}(a)\otimes m.$$
\medskip \noindent {\bf 4.3.3. Exemple:
$F$-isocristaux.} Soient $k_0$ un corps parfait de caract\'eristique $p$, $A$ le corps des fractions de
l'anneau des vecteurs de Witt de $k_0$,
$\sigma$ l'endomorphisme de Frobenius de $A$. Un $F$-isocristal sur $k_0$ est un $A$-espace
vectoriel de dimension finie $M$ muni d'un automorphisme $\sigma$-lin\'eaire $\Phi$. Ils
forment une cat\'egorie tannakienne sur ${\bf Q}_p$. On peut les interpr\'eter comme connexions
non-commutatives, et consid\'erer les r\'esultats de [Saa72]VI, 3.2, 3.4, comme des calculs de
groupes de Galois diff\'erentiels.      

\medskip {\bf 4.4. Equations lin\'eaires mixtes (diff\'erentielles-aux diff\'erences).} 
 \medskip \noindent {\bf 4.4.1.}
On consid\`ere une $k$-alg\`ebre $A$ (commutative unitaire), munie d'endomorphismes
$\sigma_1,\ldots,\sigma_m$. Pour tout $i\in \lbrace1,\ldots,m\rbrace$, on se donne un
$A$-$A$-bimodule quotient de $\Omega^1_{\sigma_i}$. On a donc un homomorphisme canonique de
$A$-$A$-bimodules $$I=Ker(A\otimes_{\scriptstyle k} A \rightarrow A) \longrightarrow \oplus
\;\Omega^1_{\sigma_i}$$ 
dont on notera simplement $\Omega^1$ l'image, et une d\'erivation
canonique 
\medskip \centerline{$d: A \longrightarrow \Omega^1$}
\par \noindent telle que $\Omega^1 = dA.A \;(= A.dA)$.   
 \medskip Cette situation correspond \`a la situation \'etudi\'ee par Bialynicki-Birula
[Bi62], qui se limite au cas o\`u $A$ est un corps de caract\'eristique $0$ (outre des
endomorphismes de $A$, cet auteur se donne aussi une famille de d\'erivations; elles sont prises
en compte ici dans $\Omega^1_{id}$). Elle recouvre, semble-t-il, toutes
celles o\`u interviennent des syst\`emes lin\'eaires mixtes diff\'erentiels-aux
diff\'erences, qu'on interpr\'etera en termes de connexions $\nabla: M
\rightarrow \Omega^1\otimes_{\scriptstyle A} M$ en g\'en\'eralisant la construction ci-dessus. On a donc cod\'e
une famille d'endomorphismes ou d\'erivations de $A$ en une seule d\'erivation $d$ (mais \`a
valeurs dans un bimodule qui n'est pas forc\'ement $A$). 
\medskip \noindent {\bf 4.4.2. Exemples.} $i)$ $A= k[z_1,z_2,\ldots ,z_m]$ (ou un
localis\'e/compl\'et\'e convenable), $\sigma_i(z_j)=z_j$ si $i\neq j$, $\sigma_i(z_i)$ est
fonction de la variable $z_i$ seule. Alors les $\delta_i: f \mapsto {f^{\sigma_i}-f \over
z^{\sigma_i}-z}\;$ (resp. ${\partial f\over {\partial z_i}}$ si $\sigma_i=id$) commutent deux \`a
deux. On prend $\Omega^1 = \oplus \;\Omega^1_{\sigma_i}\cong \oplus \; dz_i.A, \; $ et $d$ est
donn\'ee par $ df=
\sum dz_i.\delta_i(f)$.   
\par \noindent Comme cas particulier de cet exemple, en choisissant $A={\bf
Q}(a,b,c,z),\;\sigma_1:a\mapsto a+1,\sigma_2:b\mapsto b+1,\sigma_3:c\mapsto c+1,$ $\sigma_4=id,
\delta_i=\sigma_i - id \;(i\leq 3),\;\delta_4= {d\over dz}$, on traitera de mani\`ere uniforme
les relations de contigu\"{\i}t\'e et l'\'equation diff\'erentielle pour la fonction
hyperg\'eom\'etrique de Gauss $_2F_1(a,b;c;z)$. 
\par \noindent $ii)$ G\'en\'eralisant 4.3.2, les $F$-isocristaux sur une base affine lisse, ou
(dans un langage diff\'erent mais \`a peu pr\`es \'equivalent) les
\'equations diff\'erentielles lin\'eaires $p$-adiques munies d'une structure de Frobenius
forte, s'interpr\`etent comme des exemples de telles connexions ``mixtes".   

\bigskip {\bf 5. Aper\c cu sur le calcul diff\'erentiel quantique d'A. Connes.} 

\medskip {\bf 5.1.} L'id\'ee g\'en\'erale ([Con94] IV) est de ``quantifier" le calcul
diff\'erentiel en rempla\c cant la diff\'erentielle classique $df$ par une diff\'erentielle
op\'eratorielle $\d f=i[F,f]$. Ici, $f$ est un \'el\'ement arbitraire d'une $\bf C$-alg\`ebre
involutive $(A, \ast)$ d'op\'erateurs d'un espace de Hilbert {\Goth H}, et $F$ est un
op\'erateur auto-adjoint de carr\'e nul de {\Goth H}. On impose aux $\d f$ d'\^etre compacts
(ou seulement, suivant la situation, $p$-sommables). Variante ${\bf Z}/2$-gradu\'ee: {\Goth
H}={\Goth H}$^{\pm},\;F$ est impair et $[\;,\;]\;$ est un supercommutateur.   
\par L'alg\`ebre diff\'erentielle gradu\'ee $\Omega^{\ast}$ du calcul diff\'erentiel
quantique s'obtient en posant $\Omega^n= \lbrace \sum a^0\d a^1 \ldots \d a^n, \;a^i \in
A\rbrace$. La trace (r\'egularis\'ee) des op\'erateurs remplace l'int\'egration des formes. 
\par L'exemple de base est celui o\`u $A=C^{\infty}(X)$ est l'alg\`ebre des fonctions sur une
vari\'et\'e spinorielle $X$, {\Goth H} est l'espace $L^2(X,{\cal S})$ des sections de carr\'e
int\'egrable du fibr\'e des spineurs sur $X$, et $F$ est le signe $D\vert D\vert^{-1}$ de
l'op\'erateur de Dirac
$D$ (auto-adjoint et non born\'e). On rappelle que la distance g\'eod\'esique riemannienne sur
$X$ s'obtient par la formule $d(x,x')= \sup \lbrace \vert x(f)-x'(f)\vert ,\; f\in
A,\;\vert\vert [D,f]\vert\vert\leq 1\rbrace
,\; x $ et $x'$ \'etant vus comme caract\`eres de $A$. 
\medskip {\bf 5.2.} Dans cet ordre d'id\'ees, Connes a donn\'e une
pr\'esentation g\'eom\'etrique du mod\`ele de Weinberg-Salam des interactions \'electro-faibles
([Con90] V). L'espace $X$ choisi est somme disjointe de deux copies $M$ et $M'$ d'une vari\'et\'e
spinorielle compacte de dimension $4$, plac\'ees \`a tr\`es courte distance $\ell$.
L'op\'erateur de Dirac est modifi\'e en sorte que  
$$i[D, (a,a')]=(\matrix{d(a,a') & i(a'-a)/\ell \cr i(a-a')/\ell & d(a,a')});$$ ainsi la condition
$\vert\vert [D,(a,a')]\vert\vert\leq 1$ \'equivaut, \`a une constante pr\`es, \`a $\vert
a(x)-a(x')\vert\leq d(x,x'), \;\vert
a'(x)-a'(x')\vert\leq d(x,x'), \;\vert
a(x)-a'(x)\vert\leq \ell$ pour tous $x,x' \in M$.  
\par \noindent Nous avons l\`a une situation ``semi-classique" o\`u $A =
C^{\infty}(X) = C^{\infty}(M)^2$ est commutative, mais o\`u le bimodule
$\Omega^1$ ne l'est pas; l'apparition de diff\'erentielles discr\`etes $i(a'-a)/\ell$ a inspir\'e
l'interpr\'etation des diff\'erences finies propos\'ee ci-dessus (\S 4) comme diff\'erentielles
non-commutatives. 
\par Dans {\it loc. cit.}, les objets g\'eom\'etriques primordiaux sont le fibr\'e $\cal E$ non
trivial de fibre $\bf C$ sur $M$ et ${\bf C}^2$ sur $M'$, l'espace des connexions unitaires sur
$\cal E$ de ``trace" nulle sur $M'$ (le groupe $U(1)\times SU(2)$ de Weinberg-Salam
appara\^{\i}t alors comme groupe de sym\'etrie de ces connexions), et une fonctionnelle d'action
sur ces connexions construite \`a partir de leur courbure. En collaboration avec Lott
([Con94] VI), Connes est ensuite parvenu \`a une pr\'esentation g\'eom\'etrique du mod\`ele
standard (groupe de sym\'etrie $U(1)\times SU(2)\times SU(3)$), o\`u l'alg\`ebre $A$ n'est
toutefois plus commutative. 

 \bigskip \bigskip 
  \vfill \eject

\centerline{\big{ {\bf \S { II}} {\bf Calcul diff\'erentiel non-commutatif et connexions.}}}

\bigskip {\bf 1. Alg\`ebres diff\'erentielles gradu\'ees et anneaux diff\'erentiels
g\'en\'eralis\'es.} 

\medskip {\bf 1.1. A.d.g.} Le calcul diff\'erentiel non-commutatif s'est beaucoup
d\'evelopp\'e r\'ecemment \`a partir de trois sources: topologie alg\'ebrique [Kar87][Kar95],
g\'eom\'etrie non-commutative [Con90][Con94] - en particulier \`a propos du mod\`ele de
Connes-Lott du mod\`ele standard (Dubois-Violette [Dub97], Kastler, Kerner, Madore, Masson,
Michor, Mourad [Mou95], Testard...) - et groupes quantiques (Cartier, Klimyk, Maltsiniotis,
Schm\"udgen...).
\par Le cadre g\'en\'eral est celui des alg\`ebres diff\'erentielles gradu\'ees (a.d.g.). On part
d'une $k$-alg\`ebre associative unitaire $A$, non n\'ecessairement commutative, et on consid\`ere
une a.d.g.
$\Omega^{\ast}=\oplus_{n\geq 0}\Omega^n$, avec $\Omega^0=A$. Il est sous-entendu que la
diff\'erentielle $d$ est de degr\'e $1$ (et bien entendu $d(k.1)=0, d^2=0$ et
$d(\omega.\omega')=d\omega.\omega'+(-1)^{deg\;\omega}\omega.d\omega'$). Le noyau $C$ de $d$
dans $A$ est une sous-$k$-alg\`ebre, l'alg\`ebre des {\it constantes}.

\medskip \noindent {\bf 1.1.1. Exemples.} $i)$ Il existe une a.d.g. {\it universelle}
$\Omega^{\ast}_{\scriptscriptstyle univ}
=\Omega^{\ast}_{\scriptscriptstyle univ}(A/C)\;$ avec
$\Omega^0=A, \;Ker_{\Omega^0}(d) = C\;$: c'est l'alg\`ebre tensorielle
sur le $A$-$A$-bimodule
$\Omega^1_{\scriptscriptstyle univ}
 = I = Ker(A\otimes_{\scriptstyle C} A
\rightarrow A)$ noyau de la multiplication (avec $d: x \in A \mapsto 1\otimes
x - x\otimes 1$). Comme $A$-module \`a gauche, $\Omega^n_{\scriptscriptstyle univ}
$ s'identifie \`a
$A\otimes_{\scriptstyle C} (A/C)\otimes_{\scriptstyle C} \ldots \otimes_{\scriptstyle C} (A/C)$ par l'application $a_0 da_1 \ldots da_n
\mapsto a_0\otimes(a_1\; mod\; C)\otimes\ldots\otimes(a_n\; mod\; C)$ (cf. [Kar87]I). 
\par \noindent Par exemple, si $A$ est l'anneau des fonctions sur une $k$-vari\'et\'e alg\'ebrique
affine $X$, et $C=k$, alors les \'el\'ements de ${\Omega^1}_{\scriptscriptstyle univ}
$ s'identifient
aux fonctions $f(x,x')$ sur $X\times X$ qui s'annulent sur la diagonale: $\sum adb$ correspond
\`a la fonction $f(x,x')= \sum a(x)(b(x')-b(x))$; on peut ainsi voir
$\Omega^{\ast}_{\scriptscriptstyle univ}
$ comme un ``calcul aux diff\'erences universel".    

\medskip \noindent $ii)$ Si $A$ est commutative, il existe une a.d.g. {\it commutative} (au
sens gradu\'e) universelle avec $\Omega^0=A,\;Ker_{\Omega^0}(d) = C\;$: c'est le quotient de
$\Omega^{\ast}_{\scriptscriptstyle univ}
$ par l'id\'eal bilat\`ere diff\'erentiel engendr\'e par
$I^2\subset I$, qui n'est autre que l'alg\`ebre antisym\'etrique sur le module des
diff\'erentielles de K\"ahler $\Omega^1_{A/C}=I/I^2$ (au niveau de $\Omega^1_{A/C}$,
l'antisym\'etrie r\'esulte de l'application de $d$ \`a la formule
$a.db=db.a$). En caract\'eristique $\neq 2$, c'est donc l'a.d.g. de De Rham
$\Omega^{\ast}_{DR}(A/C)=\Lambda^{\ast}\Omega^1_{A/C}.$  

\medskip \noindent $iii)$ Dans la situation de I.4.2, il est naturel d'introduire
l'a.d.g. $\Omega^{\ast}_{\sigma}$ quotient de $\Omega^{\ast}_{\scriptscriptstyle univ}
$ par 
l'id\'eal bilat\`ere diff\'erentiel engendr\'e par $JI$; son terme de degr\'e $1$ est le
sesquimodule $\Omega^1_{\sigma}$ consid\'er\'e en I.4.2.1. Remarquons que $\Omega^n_{\sigma}$ est
un sesquimodule relativement \`a l'endomorphisme $\sigma^n$ de $A$ (point de vue
proche de celui des $\epsilon$-diff\'erentielles de [LT99]; voir aussi [Su99]). 
\par \noindent Dans le cas ``\`a une variable" I.4.2.2.$i)$, le calcul
$dz.dz=d(z.dz)=d(dz.z^{\sigma})= -dz.dz.\delta_{\sigma}(z^{\sigma})$ montre qu'en
caract\'eristique $\neq 2$, {\it et si $\sigma(z)\neq -z$}, alors $\Omega^{>1}=0$.

\medskip \noindent $iv)$ De m\^eme, dans le cas ``\`a plusieurs variables" I.4.4.2, on prend pour
$\Omega^{\ast}$ le quotient de $\Omega^{\ast}_{\scriptscriptstyle univ}
$ par 
l'id\'eal bilat\`ere diff\'erentiel engendr\'e par $Ker(\Omega^1_{\scriptscriptstyle univ}
\rightarrow
\oplus \Omega^1_{\sigma_i})$. Dans $\Omega^2$, on observe alors que
$dz_i.dz_j=d(z_i.dz_j)=d(dz_j.z_i^{\sigma_j})= -dz_j.dz_i$ si
$i\neq j$, et $2dz_i.dz_i=0$ si $\sigma_i(z_i)\neq -z_i$; on en tire alors ais\'ement, en
caract\'eristique $\neq 2$, que si $\sigma_i(z_i)\neq -z_i$ pour tout $i$, $\Omega^2 \cong
\oplus_{i<j} \; dz_i.dz_j.A$ (comme $A$-module \`a droite). 

\medskip \noindent {\bf 1.1.2.} On dit qu'une a.d.g. $\Omega^{\ast}$ est {\it r\'eduite} ou
engendr\'ee par $\Omega^0$ si c'est la plus petite sous-a.d.g. de cran $0$ \'egal \`a
$\Omega^0$. Il revient au m\^eme de dire que $\Omega^{\ast}$ est quotient de
$\Omega^{\ast}_{\hbox{\smit{\hbox{\smit{univ}}}}}$. 

\medskip {\bf 1.2. Anneaux diff\'erentiels (g\'en\'eralis\'es).} 
\medskip \noindent {\bf 1.2.1.} La recherche d'un cadre unifi\'e o\`u inscrire les exemples de I
incite \`a g\'en\'eraliser la notion d'anneau diff\'erentiel. Nous appellerons {\it anneau
diff\'erentiel} (g\'en\'eralis\'e) la donn\'ee d'une d\'erivation
$$d: A\longrightarrow \Omega^1$$
d'une $k$-alg\`ebre associative unitaire $A$ \`a valeurs dans un $A$-$A$-bimodule $\Omega^1$ (il
est sous-entendu que le $k$-$k$-bimodule sous-jacent \`a $\Omega^1$ est
commutatif; autrement dit, $\Omega^1$ est un $A\otimes_k A^{op}$-module \`a gauche). La notion usuelle
d'anneau diff\'erentiel correspond au cas $\Omega^1=A$.   
\par \noindent Par abr\'eviation abusive, nous parlerons quelquefois de l'anneau
diff\'erentiel $(A,d)$, ou m\^eme $A$. Nous noterons couramment $C$ la $k$-alg\`ebre des
constantes $Ker(d)$. 
\medskip \noindent {\bf 1.2.2.} On dit que l'anneau diff\'erentiel $(A{\buildrel{d}\over
\rightarrow}\Omega^1)$ est {\it r\'eduit} si l'image de $\;d\;$
engendre $\Omega^1$ comme $A$-module \`a droite. Par la r\`egle de Leibniz $d(ab)=d(a)b+ad(b)$,
on voit qu'elle l'engendre aussi $\Omega^1$ comme $A$-module \`a  gauche, i.e. $$\Omega^1= dA.A =
A.dA = A.dA.A.$$ 
\par \noindent {\bf 1.2.3.} Un {\it morphisme} d'anneaux diff\'erentiels 
$(A {\buildrel{d}\over
\rightarrow}\Omega^1) \rightarrow (A'{\buildrel{d'}\over
\rightarrow}\Omega'^1)$ est un couple
$u=(u^0, u^1)$ form\'e d'un
morphisme de $k$-alg\`ebres unitaires $A {\buildrel{u^0}\over
\rightarrow}A'$ et d'une application $\Omega^1{\buildrel{u^1}\over
\rightarrow}\Omega'^1$ v\'erifiant
$$u^1\circ d=d' \circ u^0,\;u^1(a\omega b)=u^0(a)u^1(\omega)u^0(b)$$ pour tous $a,b \in A, \omega
\in \Omega^1$. 
\par \noindent L'un des objets de la th\'eorie de Galois diff\'erentielle (g\'en\'eralis\'ee) est
l'\'etude des automorphismes de certains anneaux diff\'erentiels, cf. infra III.3.2.  
\medskip \noindent Toute a.d.g. fournit, par restriction \`a ses termes de degr\'e $0$ et $1$, un
anneau diff\'erentiel au sens ci-dessus. Ceci d\'efinit un foncteur ``d'oubli".  

\proclaim Lemme 1.2.4. Le foncteur d'oubli 
$\;\lbrace {\hbox{a.d.g. r\'eduites}}\rbrace \rightarrow \lbrace {\hbox{anneaux
diff\'erentiels r\'eduits}}\rbrace\;$ admet un adjoint \`a gauche ``a.d.g", qui associe \`a
l'anneau diff\'erentiel $\;A {\buildrel{d}\over \rightarrow}\Omega^1\;$ l'alg\`ebre
diff\'erentielle gradu\'ee
$\;\Omega^{\ast}= ``a.d.g"(A{\buildrel{d}\over \rightarrow}\Omega^1)\;$ quotient de
$\;\Omega^{\ast}_{\scriptscriptstyle univ}
(A/k)\;$ par l'id\'eal bilat\`ere diff\'erentiel engendr\'e
par $\;Ker(\Omega^1_{\scriptscriptstyle univ}
\rightarrow \Omega^1)$. 

\par \noindent En effet, pour toute a.d.g. $\Omega'^{\ast}$, la surjectivit\'e de
$$Hom_{a.d.g.}(``a.d.g"(A{\buildrel{d}\over
\rightarrow}\Omega^1),\Omega'^{\ast})\; \longrightarrow \; Mor(A {\buildrel{d}\over
\rightarrow}\Omega^1, \;A'{\buildrel{d'}\over
\rightarrow}\Omega'^1)$$ r\'esulte ais\'ement de la propri\'et\'e universelle de
$\Omega^{\ast}_{\scriptscriptstyle univ}
(A/k)$, et l'injectivit\'e de ce qu'un morphisme $``a.d.g"(A{\buildrel{d}\over
\rightarrow}\Omega^1)\rightarrow \Omega'^{\ast}$ est d\'etermin\'e par sa valeur sur $\Omega^0 =
A$ et sur $\Omega^1= A.dA.A$. Noter qu'on aurait pu remplacer
$\Omega^{\ast}_{\scriptscriptstyle univ}
(A/k)$ par $\Omega^{\ast}_{\scriptscriptstyle univ}
(A/C)$.  

\medskip \noindent {\bf 1.2.5.} Une {\it extension diff\'erentielle} est un morphisme
d'anneaux diff\'erentiels $u=(u^0, u^1)$ avec $u^0$ injectif, et o\`u $\Omega'^1\cong \Omega^1
\otimes_{\scriptstyle A} A'$ comme $A'$-module \`a droite (la dissym\'etrie est li\'ee \`a 
notre choix de privil\'egier, ult\'erieurement, les connexions \`a gauche).

 \medskip {\bf 1.3. Id\'eaux diff\'erentiels.} 
\medskip \noindent {\bf 1.3.1.} Le lemme pr\'ec\'edent sugg\`ere de d\'efinir la notion d'{\it
id\'eal diff\'erentiel g\'en\'eralis\'e} de  
$(A{\buildrel{d}\over
\rightarrow}\Omega^1)$ comme trace en degr\'e $\leq 1$ d'un id\'eal diff\'erentiel (bilat\`ere)
de $\break ``a.d.g"(A{\buildrel{d}\over \rightarrow}\Omega^1),\Omega'^{\ast})$, ou ce qui
revient au m\^eme, comme noyau d'un morphisme d'anneaux diff\'erentiels $(A {\buildrel{d}\over
\rightarrow}\Omega^1) \rightarrow (A'{\buildrel{d'}\over
\rightarrow}\Omega'^1)$. Explicitement, c'est donc la donn\'ee d'un id\'eal bilat\`ere $I$
de $A$ et d'un sous-$A$-$A$-bimodule $I^1$ de $\Omega^1$ tel que $dI\subset I^1$.        

\medskip \noindent {\bf 1.3.2.} Cette notion naturelle du point de vue cat\'egorique semble peu
utile en pratique. Aussi r\'eserverons-nous le nom d'{\it id\'eal diff\'erentiel} (au sens strict)
aux id\'eaux diff\'erentiels g\'en\'eralis\'es
$(I{\buildrel{d}\over \rightarrow}I^1)$ qui v\'erifient $I^1= \Omega^1.I,$ l'image de
$\Omega^1 \otimes I$ dans $\Omega^1$ (noter la dissym\'etrie droite-gauche). C'est donc le noyau
d'un morphisme d'anneaux diff\'erentiels $u=(u^0, u^1)$ avec $u^1\otimes 1: \Omega^1\otimes
A' \rightarrow \Omega'^1$ bijectif. Comme l'id\'eal diff\'erentiel est alors d\'etermin\'e par
$I$, nous dirons aussi, par abus, que $I$ est un id\'eal diff\'erentiel de $A$.
\par Consid\'erons le cas particulier o\`u le $A$-module \`a droite $\Omega^1$ est projectif de
type fini (c'est le cas dans tous les exemples int\'eressants). Il admet donc un ``dual \`a
gauche"
${}^\vee\Omega^1$, qui est un
$A$-$A$-bimodule. On dispose des homomorphismes de $A$-$A$-bimodules habituels (cf.
[Bru94]1) 
$$\epsilon :\;\; {}^\vee\Omega^1 \otimes_{\scriptstyle A} \Omega^1 \rightarrow A \;\;\;\hbox{(\'evaluation,
not\'ee aussi }\;\langle \;,\;\rangle )$$
$$\eta :\;\; A\rightarrow\Omega^1 \otimes_{\scriptstyle A} {}^\vee\Omega^1 
\;\;\;\hbox{(co\'evaluation)}$$
$$\hbox{v\'erifiant}\;\;\;\;(id\otimes \epsilon)(\eta\otimes id)=id_{\Omega^1},\;\;(
\epsilon \otimes id)(id\otimes \eta)=id_{({}^\vee\Omega^1)}. $$

\proclaim Lemme 1.3.3. Supposons $\Omega^1$ fid\`ele et projectif de
type fini \`a droite sur $A$. Un id\'eal bilat\`ere $I$ de $A$ est un id\'eal diff\'erentiel si
et seulement si $\;\langle {}^\vee\Omega^1 , dI\rangle \; \subset I$, ou encore si et seulement si
$\;\langle {}^\vee\Omega^1 , dI.A\rangle \; = I$.
\par \noindent Il suffit de prouver la seconde assertion. Remarquons que $J=\;\langle {}^\vee\Omega^1
, dI.A\rangle \;= \epsilon({}^\vee\Omega^1 \otimes A.dI.A)$ est un id\'eal bilat\`ere de $A$. On a
$AdI.A= (\epsilon \otimes id)(id\otimes \eta)(AdI.A) \subset \Omega^1\otimes J\cong \Omega^1 J$.
Il reste \`a d\'emontrer que $J$ contient $I$. Puisque
$\Omega^1$ est fid\`ele et projectif de type fini \`a droite sur $A$, il suffit de faire voir que
$\Omega^1\otimes I \subset \Omega^1\otimes J$. Or pour tous $a,b \in A,\;i\in I,\; adb.i =
ad(bi)-ab.di\in \Omega^1\otimes J$.   
\medskip Dans le cas classique o\`u $A$ et $\Omega^1$ sont commutatifs, les
\'el\'ements de ${}^\vee\Omega^1$ s'interpr\`etent comme des d\'erivations de $A$, et on retrouve
une notion famili\`ere d'id\'eal diff\'erentiel: id\'eal stable sous ${}^\vee\Omega^1 \subset
Der(A).$    
\medskip \noindent {\bf 1.3.4.} Un anneau diff\'erentiel $(A,d)$ est dit {\it simple} si ses seuls
id\'eaux diff\'erentiels sont $0$ et $A$. 
\par  \noindent {\it Exemples.} $\bullet$ Les anneaux
diff\'erentiels de fonctions analytiques sur une vari\'et\'e analytique lisse (r\'eelle ou complexe,
voire $p$-adique rigide) sont simples, mais pas les anneaux diff\'erentiels de fonctions $C^{\infty}$ (les
fonctions infiniment plates en un point forment un id\'eal diff\'erentiel).  
\par \noindent $\bullet$ De m\^eme, les anneaux locaux d'une vari\'et\'e alg\'ebrique $X$ sur
un corps $k$ de caract\'eristique nulle sont diff\'erentiellement simples ($\Omega^1$ \'etant pris
\'egal au module des diff\'erentielles de K\"ahler) si $X$ est lisse, mais pas en g\'en\'eral.    
\par \noindent $\bullet$ Les anneaux diff\'erentiels (de caract\'eristique $p$) $\bar{\cal
O}_x$ et ${\cal O}^{pd}_x$ consid\'er\'es en I.2.2 sont simples. 
\par \noindent $\bullet$ Dans la situation d'un anneau aux diff\'erences $(A,\sigma)$ (I.4.2),
on prendra garde que la simplicit\'e de $A\rightarrow \Omega^1_{\sigma}$ n'\'equivaut pas \`a ce
que les seuls id\'eaux de $A$ stables par $\sigma$ sont $0$ et $A$. Par exemple, cet
anneau diff\'erentiel est simple dans le cas de $A=k[z]$ muni de la dilatation $\sigma z = qz$,
bien que $zA$ soit un id\'eal stable sous $\sigma$.
\par \noindent $\bullet$ Si $(A',d')$ est extension d'un anneau
diff\'erentiel simple $(A,d)$ (cf. 1.2.5), tout id\'eal diff\'erentiel de $(A',d')$ distinct de
$A'$ \'evite $A$ priv\'e de $0$; en particulier $(A',d')$ est simple si $A'$ est une localisation
de $A$. 
\par \noindent $\bullet$ Soit $k={\bf C}[[x,y]],\;A=k[z_1,z_2]/(xz_1+yz_2-1),\;\Omega^1=\Omega^1_{A/k}$.
Alors $C=k$ et l'id\'eal maximal {\goth p} de $k$ v\'erifie  {\goth p}$A=A$. 
On a toutefois le r\'esultat suivant:
\par 
\proclaim Lemme 1.3.5. Supposons $\;A \;$ est commutatif. 
\hfill \break $i)$ Si $\;(A,d)\;$ est simple, alors $\;C \;$ est un corps. 
\hfill \break $ii)$ Soit $\;Q(A)\;$ l'anneau total des fractions de $\;A \;$ (i.e. le localis\'e de
$\;A \;$ par le mono\"{\i}de de tous ses \'el\'ements non-diviseurs de $0$). Alors $\;d \;$ admet un unique
prolongement en une d\'erivation $\;d: Q(A)\rightarrow Q(A)\otimes\Omega^1\otimes
Q(A)$. Supposons en outre que $\Omega^1$ soit un $A$-module \`a droite sans torsion, que $\Omega^1
\otimes_{\scriptstyle A}Q(A) \rightarrow Q(A)\otimes_{\scriptstyle A}
\Omega^1\otimes_{\scriptstyle A}Q(A)$ soit surjectif, et que $\;(A,d)\;$ soit simple. Alors
$\;(Q(A),d)\;$ est simple, et son corps de constantes est $C$. 
\par {\it Preuve.} $i)$ Si $c$ est un \'el\'ement non nul de $C$, la simplicit\'e de $(A,d)$ entra\^{\i}ne
que $cA=A$, donc $c$ a un inverse, qui est n\'ecessairement une constante.  
\par \noindent $ii)$ On peut voir $A$ comme sous-$k$-alg\`ebre de $Q(A)$. L'existence et l'unicit\'e
du prolongement de $d$ \`a $Q(A)$ sont claires: pour $a$ inversible, on a la formule
$d(a^{-1})=-a^{-1}\otimes da \otimes a^{-1}$. Sous les hypoth\`eses que $\Omega^1$ soit un $A$-module \`a
droite sans torsion et que $\Omega^1 \otimes_{\scriptstyle A}Q(A) \rightarrow
Q(A)\otimes_{\scriptstyle A} \Omega^1\otimes_{\scriptstyle A}Q(A)$ soit surjectif, on a $d:
Q(A)\rightarrow \Omega^1\otimes_{\scriptstyle A}Q(A)$, et $\Omega^1 \hookrightarrow
\Omega^1\otimes_{\scriptstyle A}Q(A)$. On voit alors que l'intersection de tout id\'eal
diff\'erentiel de $Q(A)$ avec $A$ est un id\'eal diff\'erentiel de $A$. Si $(A,d)$ est simple, on
en d\'eduit qu'il en est de m\^eme de $(Q(A),d)$. Par ailleurs, soit $c$ une constante non nulle de
$Q(A)$. Alors $cA \cap A$ est un id\'eal diff\'erentiel non nul de $A$, donc \'egal \`a $A$; on en
d\'eduit $1/c \in A$, d'o\`u $1/c \in C$ et $c \in C$.     
  
\medskip \noindent {\bf 1.3.6.} Pour traiter du cas o\`u $A$ est commutatif mais o\`u
$C$ n'est pas un corps (c'est-\`a-dire, en pratique, du cas de familles d'\'equations diff\'erentielles ou
aux diff\'erences), il y a lieu de remplacer la notion d'anneau diff\'erentiel simple par celle d'anneau
diff\'erentiel ``simple par couches". Pour tout id\'eal premier {\goth p} de
$C$, notons $\kappa$({\goth p}) le corps de fractions de $C/${\goth p}. 
\par \noindent Nous dirons que $(A,d)$ est {\it simple par couches} si $A$ est fid\`element plat
sur $C$, et si pour tout id\'eal premier {\goth p} de $C$, l'anneau diff\'erentiel
$\;\kappa$({\goth p})$\otimes_{\scriptstyle C} A \rightarrow \kappa$({\goth p})$\otimes_{\scriptstyle C} A.dA\;$ induit est
simple, et l'anneau des constantes de $\kappa$({\goth p})$\otimes_{\scriptstyle C} A$ est $\kappa$({\goth p}). 
\par \noindent {\it Exemples.} $\bullet$ Tout anneau diff\'erentiel simple est simple par couches.
\par \noindent $\bullet$ Si $(A,d)$ est simple par couches, il en est de m\^eme de l'anneau
diff\'erentiel $(A,d)\otimes C'$ obtenu par changement d'anneau de constantes $C\rightarrow C'$,
d\`es lors que $C'$ est une localisation ou un quotient de $C$; en effet tout id\'eal premier
de $C'$ provient alors d'un id\'eal premier de $C$.   
\par \noindent $\bullet$ Soit $A=({\bf C}[[q-1]])[z]$ (ou un localis\'e),
$\sigma=$ dilatation de module $q$, $\Omega^1 = \Omega^1_{\sigma}$ (cf. I.4.2.2 $i)$, 4.3;
situation correspondant \`a la confluence des \'equations aux $q$-diff\'erences vers une
\'equation diff\'erentielle). C'est un anneau diff\'erentiel simple par couches: le
point est que $\delta_{\sigma}(z^n)=n_q.z^{n-1}$ et $n_q=1+q+\ldots +q^{n-1}$ est une unit\'e
dans ${\bf C}[[q-1]])$.
\par \noindent $\bullet$ Les anneaux diff\'erentiels dont l'anneau des constantes est de
caract\'eristique mixte ne sont g\'en\'eralement pas simples par couches, e.g.
${\bf Z}[z]{\buildrel{d}\over \rightarrow}{\bf Z}[z]dz$ (l'anneau des constantes de $({\bf
Z}/p{\bf Z})[z]$ est $({\bf Z}/p{\bf Z})[z^p]$).  

\medskip \noindent {\bf 1.3.7. Remarque.} Toutes ces notions se faisceautisent sans difficult\'e.
En fait, la plupart des constructions qui suivent gardent un sens dans tout topos.

\bigskip {\bf 2. Connexions.}

\medskip {\bf 2.1. D\'efinition.} Soient $A {\buildrel{d}\over
\rightarrow}\Omega^1$ un anneau diff\'erentiel (g\'en\'eralis\'e), $C$ son anneau de constantes.
Soit $M$ un $A$-module \`a gauche. Une {\it connexion} sur $M$ est une application $k$-lin\'eaire 
\medskip \centerline{$\nabla: M \longrightarrow \Omega^1\otimes_{\scriptstyle A} \;M$}
\medskip \noindent v\'erifiant la r\`egle de Leibniz
\medskip \centerline{$\nabla(am)=a\nabla(m)+da\otimes m.$}
\medskip \noindent On peut aussi, selon Atiyah, d\'efinir une connexion en introduisant le module
${\cal P}^1(M)$ des {\it jets \`a gauche d'ordre un} [Kar87]1.8: comme $C$-module, ${\cal P}^1(M)$
n'est autre que $M\oplus (\Omega^1\otimes_{\scriptstyle A}
\;M)$; on le munit d'une structure de $A$-module \`a gauche en posant 
\medskip \centerline{$a(m,\omega\otimes n) =(am,
da\otimes m + a\omega \otimes n).$} 
\medskip \noindent
La donn\'ee d'une connexion $\nabla$ \'equivaut alors \`a celle d'un scindage $D=(id,\nabla)$ de
la suite naturelle de $A$-modules
\medskip \centerline{$0 \rightarrow \Omega^1\otimes_{\scriptstyle A} \;M \rightarrow {\cal P}^1(M)
\;{\buildrel{pr_1}\over
\rightarrow}\; M \rightarrow 0.$}
\medskip \noindent Un morphisme de $A$-modules \`a connexion est un homomorphisme de
$A$-modules {\it horizontal}, i.e. compatible aux connexions. Les
$A$-modules \`a connexion forment une cat\'egorie ab\'elienne
$k$-lin\'eaire (et m\^eme $C$-lin\'eaire \`a gauche).  
\medskip \noindent {\bf 2.1.1. Exemples.} $\bullet$ $A {\buildrel{d}\over
\rightarrow}\Omega^1$ d\'efinit une connexion sur $A$, appel\'ee connexion triviale. Plus
g\'en\'eralement, on qualifie de {\it triviale} toute connexion isomorphe \`a
une connexion de la forme $d\otimes id_N$ sur $A\otimes_{\scriptstyle C} N$, o\`u $N$ est un $C$-module \`a
gauche quelconque.
\par \noindent Pour tout $A$-module \`a connexion $(M,\nabla)$, on peut identifier ${\rm
Mor}((A,d),(M,\nabla))$ au $C$-module $M^{\nabla}:=Ker\nabla$. Les \'el\'ements de
$M^{\nabla}$ sont dits {\it horizontaux} (en souvenir de la situation g\'eom\'etrique
rappel\'ee en I.1.2).    
\par \noindent $\bullet$ Si $\Omega^1$ est projectif de type fini \`a droite, alors les id\'eaux
diff\'erentiels de $A {\buildrel{d}\over
\rightarrow}\Omega^1$ (1.3.2) correspondent exactement aux sous-modules \`a connexion de la
connexion triviale $(A,d)$.    
\par \noindent $\bullet$ Si $\Omega^1=A$, la notion de connexion se sp\'ecialise en la notion
usuelle de module diff\'erentiel: un endomorphisme $k$-lin\'eaire $\nabla$ de $M$ v\'erifiant
$\nabla(am)= da.m+a\nabla(m)$.   
\medskip \noindent {\bf 2.1.2. Remarque.} Lorsque $\Omega^1$ est projectif de type fini \`a
droite (de dual \`a gauche ${}^{\vee}\Omega^1$), la donn\'ee d'une connexion $\nabla$ \'equivaut
\`a celle d'une application $k$-lin\'eaire 
$\;\;{}^{\vee}\Omega^1\rightarrow End_{\scriptstyle k}\;M, \;\;\;D\mapsto
\nabla_{\scriptstyle D}\;\;$ v\'erifiant 
\medskip \centerline{$\nabla_{\scriptstyle D}(am)=\nabla_{D.a}(m)+\langle
D,da\rangle .m \;\; \hbox{(en posant }\;\nabla_{\scriptstyle D}(m)=\langle D,\nabla(m)\rangle \;
\in M).$}
\medskip \noindent En g\'en\'eral les op\'erateurs
$\nabla_{\scriptstyle D}$ ne v\'erifient pas la r\`egle de Leibniz, d'o\`u la sup\'eriorit\'e de la formulation
en termes de connexions. 
\medskip Rappelons, pour m\'emoire des signes, le lemme
bien connu suivant:
\proclaim Lemme 2.1.3. Pour toute a.d.g. $\Omega^{\ast}$ ``prolongeant" $A
{\buildrel{d}\over \rightarrow}\Omega^1$, $\nabla$
s'\'etend en une unique application $C$-lin\'eaire 
$\;\Omega^{\ast}\otimes_{\scriptstyle A} M \longrightarrow \Omega^{\ast
+1}\otimes_{\scriptstyle A} M,\;$
 encore not\'ee $\nabla$, v\'erifiant pour tout $\omega, \omega' \in
\Omega^{\ast}$, et tout $m\in M$ l'identit\'e 
\medskip $\nabla(\omega.(\omega'\otimes m))=d\omega.(\omega'\otimes
m)+(-1)^{deg\;\omega}\omega.\nabla(\omega'\otimes m).$
\par L'existence se d\'eduit facilement du calcul $\nabla(\omega.(\omega'\otimes
m))=\nabla((\omega.\omega')\otimes m)=d(\omega.\omega')\otimes m+(-1)^{deg\;\omega+deg\;
\omega'}\omega.\omega'\nabla(m)=d\omega.\omega'\otimes m +(-1)^{deg\;\omega}d\omega'\otimes
m +(-1)^{deg\;\omega+deg\;
\omega'}\omega.\omega'\nabla(m) = d\omega.(\omega'\otimes
m)+(-1)^{deg\;\omega}\omega.\nabla(\omega'\otimes m).$
\medskip \noindent {\bf 2.1.4. Remarque.} Dans le cas des connexions sur les modules
\`a droite ([Kar87]1.10, [N97].2), la r\`egle des signes d\'efinissant l'extension $M\otimes_{\scriptstyle A}
\Omega^{\ast} \longrightarrow  M\otimes_{\scriptstyle A} \Omega^{\ast +1}$ est  
\medskip \centerline{$\nabla((m\otimes\omega).\omega')=\nabla(m\otimes
\omega).\omega'+(-1)^{deg\;\omega}(m\otimes
\omega).d\omega'.$}
\medskip \noindent Il y a donc lieu de prendre garde, m\^eme dans le cas classique o\`u
$\Omega^{\ast}$ est gradu\'ee commutative, aux changements de signe lorsqu'on \'ecrit les
connexions comme applications \`a valeurs dans $\Omega^1\otimes M$ ou dans $M\otimes \Omega^1$.

\medskip {\bf 2.2. Fonctorialit\'e.} 
\medskip \noindent Consid\'erons un morphisme d'anneaux diff\'erentiels $u=(u^0, u^1): \;(A
{\buildrel{d}\over
\rightarrow}\Omega^1) \rightarrow (A'{\buildrel{d'}\over
\rightarrow}\Omega'^1)$ (cf. 1.2.2). A tout $A$-module \`a connexion $(M,\nabla)$, on associe un
$A'$-module \`a connexion $u^{\ast}(M,\nabla)= (A'\otimes_{\scriptstyle A} M, \nabla')$, o\`u $\nabla'$ est
d\'efinie par 
\medskip \centerline{$\nabla'(a'\otimes m)=a'.(u^1\otimes id_M)(\nabla(m))+da'\otimes m.$}
\medskip  \noindent On le note aussi $(M_{A'},\nabla_{A'})$.
\par \noindent Le scindage correspondant $(id,\nabla')$ de
la suite naturelle de $A$-modules
\par \centerline{$0 \rightarrow \Omega'^1\otimes_{\scriptstyle A} \;M \rightarrow {\cal P}^1(A'\otimes_{\scriptstyle A} M)
\;{\buildrel{pr_1}\over
\rightarrow}\; A'\otimes_{\scriptstyle A} M \rightarrow 0$}
\par \noindent n'est autre que le prolongement $A'$-lin\'eaire de l'application $A$-lin\'eaire
compos\'ee 
\par \centerline{$M \;{\buildrel{(id,\nabla)}\over
\longrightarrow} \;{\cal P}^1( M) \rightarrow {\cal P}^1(A'\otimes_{\scriptstyle A} M).$}

\bigskip {\bf 2.3. Courbure, int\'egrabilit\'e, descente.} 
\medskip \noindent Le lemme suivant ([Kar87]1.13, [N97]2.5) est bien connu, surtout dans le cas commutatif
\proclaim Lemme 2.3.1. La composition $\nabla^2$ de $\nabla$ par elle-m\^eme
\medskip \centerline{$\Omega^{\ast}\otimes_{\scriptstyle A} M \longrightarrow \Omega^{\ast
+2}\otimes_{\scriptstyle A} M$} 
\medskip \noindent est un homomorphisme de $\Omega^{\ast}$-module \`a gauche, uniquement
d\'etermin\'e par sa restriction $\;\;M\longrightarrow \Omega^2\otimes_{\scriptstyle A} M$. 
\par Lorsque $\Omega^{\ast}$ est l'a.d.g.
canonique attach\'ee \`a $A {\buildrel{d}\over \rightarrow}\Omega^1\;$ (1.2.3.), on appelle cet
homomorphisme $M\rightarrow \Omega^2\otimes_{\scriptstyle A} M$ (la {\it courbure} de $\nabla$). 
 \par Le lemme suivant est imm\'ediat  
\proclaim Lemme 2.3.2. Soit $u=(u^0, u^1): \;(A
{\buildrel{d}\over
\rightarrow}\Omega^1) \rightarrow (A'{\buildrel{d'}\over
\rightarrow}\Omega'^1)$ un morphisme d'anneaux diff\'erentiels. La courbure de
$u^{\ast}(M,\nabla)$ est l'homomorphisme compos\'e de $id_{A'}\otimes \nabla^2:
A'\otimes M\rightarrow
A'\otimes \Omega^2\otimes_{\scriptstyle A} M$ et de l'homomorphisme naturel $A'\otimes \Omega^2\otimes_{\scriptstyle A} M
\rightarrow \Omega'^2\otimes_{\scriptstyle A} M$.
\par On dit que
$\nabla$ est {\it plate} ou {\it int\'egrable} si sa courbure est nulle. 
\medskip \noindent {\bf 2.3.3. Exemples.} $\bullet$ Une connexion {\it triviale}, i.e. de la forme
$d\otimes id_N$ sur $A\otimes_{\scriptstyle C} N$, o\`u $N$ est un
$C$-module \`a gauche quelconque, est toujours int\'egrable. 
\par \noindent $\bullet$ Si $A$ est commutatif de caract\'eristique $\neq 2$ (cf. 1.1.1), et si
$\Omega^1$ est un bimodule commutatif projectif de type fini, l'int\'egrabilit\'e d'une connexion
$\nabla$ \'equivaut \`a ce que l'application $k$-lin\'eaire ${}^{\vee}\Omega^1\rightarrow
End_{\scriptstyle k}\;M, \;\;D\mapsto \nabla_{\scriptstyle D} $ commute au crochet de Lie.  
\par \noindent $\bullet$ Examinons la situation de 1.1.1.$iv)$: calcul aux diff\'erences
\`a plusieurs variables, en caract\'eristique $\neq 2$, avec $\sigma_i(z_i)\neq \pm z_i$, de sorte
que $\Omega^2
\cong \oplus_{i<j} \; dz_i.dz_j.A$. Consid\'erons un syst\`eme
lin\'eaire aux diff\'erences
\medskip \centerline{${(\ast\ast\ast)}\;\;\;\;\;\;\;\sigma_i(Y)={\cal A}_{(i)}.Y$}
\par \noindent o\`u ${\cal A}_{(i)} \in GL_{\mu}(A)$. Suivant II.4.3, 4.4, associons-lui le module
libre
$M=\oplus_{k=1}^{k=\mu} A.m_{\scriptstyle k}$ muni de la connexion suivante (cf. I.4.3): 
\medskip\centerline{$\nabla(m_{\scriptstyle k})=\oplus_{i,l}\; dz_i\otimes {{\cal
A}^{-1}_{(i){lk}}-\delta_{lk}\over z_i^{\sigma_i}-z_i}\;m_l.$}
\medskip \noindent Un calcul direct montre que sa courbure est donn\'ee par 
\medskip\centerline{$\nabla^2(m_{\scriptstyle k})=\oplus_{i<j,\;l}\; dz_i.dz_j\otimes {{\cal
R}_{(ij){lk}}\over (z_i^{\sigma_i}-z_i)(z_j^{\sigma_j}-z_j)}\;m_l,$}
\medskip\centerline{${\hbox{o\`u}}\;\;{\cal R}_{(ij)}={\cal A}^{-1}_{(i)}{{\cal
A}^{-1}_{(j)}}^{\sigma_i}-{\cal A}^{-1}_{(j)}{{\cal A}^{-1}_{(i)}}^{\sigma_j},\;\;i<j.$}
\medskip \noindent La connexion est donc int\'egrable si et seulement si $\;{{\cal
A}_{(j)}}^{\sigma_i}{\cal A}_{(i)}={{\cal A}_{(i)}}^{\sigma_j}{\cal
A}_{(j)},\;\forall i,j.$
\medskip \noindent {\bf 2.3.4.} Dans le cas int\'egrable, on dispose d'un
complexe de $C$-modules $\Omega^{\ast}\otimes_{\scriptstyle A} M$, appel\'e {\it complexe de De Rham
non-commutatif de $(M,\nabla)$}. Sa cohomologie appara\^{\i}t par exemple dans les travaux
d'Aomoto sur les syst\`emes hyperg\'eom\'etriques ou $q$-hyperg\'eom\'etriques, cf. e.g. [Ao90]
(voir aussi [TV97]).  
\medskip \noindent {\bf 2.3.5.} Dans [N97]3.12, on trouve un dictionnaire entre donn\'ees de
recollement pour un
$A$-module $M$ relativement \`a $A/C$, et connexions $\nabla$ sur $M$ (avec $\Omega^1=
\Omega^1_{\scriptscriptstyle univ}
(A/C)$). La condition de cocycle pour la donn\'ee de
recollement correspond \`a l'int\'egrabilit\'e de
$\nabla$; on obtient alors une donn\'ee de descente effective: $M$ provient d'un $C$-module. 

\medskip {\bf 2.4. Biconnexions et voltes.} Supposons maintenant que $M$ soit un
$A$-$A$-bimodule, tel que le $k$-$k$-bimodule sous-jacent soit
commutatif. Rappelons qu'un {\it op\'erateur d'ordre un} de $M$ vers un autre
$A$-$A$-bimodule $N$ est une application $k$-lin\'eaire $D: M\rightarrow N$ telle que pour tout
$a\in A$, $m\mapsto D(am)-aD(m)$ est un homomorphisme de $A$-modules \`a droite; il revient au
m\^eme de dire que pour tout $b\in A$, $m\mapsto D(mb)-D(m)b$ est un homomorphisme de $A$-modules
\`a gauche ([DubM96], lemma 1; le point est qu'une translation \`a droite arbitraire
${\bf{.}}b$ commute \`a toute translation \`a gauche $a{\bf{.}}$, de sorte  que
$[[D,a{\bf{.}}],{\bf{.}}b]=[[D,{\bf{.}}b],a{\bf{.}}]$). Le r\'esultat suivant est un cas
particulier de [DubM96], theorem 1:
\proclaim Proposition 2.4.1. Il existe un unique homomorphisme de bimodules
$\varphi_s(D):\Omega^1_{\scriptscriptstyle univ}
\otimes_{\scriptstyle A} M \rightarrow N$ (resp.
$\varphi_{\scriptstyle d}(D): M\otimes_{\scriptstyle A}\Omega^1_{\scriptscriptstyle univ}
\rightarrow N$) tel que l'on ait    
\medskip $D(amb)=aD(m)b+\varphi_s(D)(d_{\scriptscriptstyle univ}
a\otimes m)b+a\varphi_{\scriptstyle d}(D)(m\otimes
d_{\scriptscriptstyle univ}
b).$ 
\par \noindent Par exemple, $\varphi_{\scriptstyle d}(D)$ est induit par l'homomorphisme de
bimodules
$M\otimes_{\scriptstyle k}\Omega^1_{\scriptscriptstyle univ}
\rightarrow N$ d\'efini par $m\otimes db.a \mapsto
D(mb)a-D(m)ba$.  
\medskip \noindent Il est clair que toute connexion (\`a gauche) sur $M$ est un
op\'erateur d'ordre $1$ de $M$ vers $N=\Omega^1\otimes_{\scriptstyle A} \;M$, avec $\varphi_s(\nabla)=\pi
\otimes id_M$, o\`u $\pi$ est l'homomorphisme canonique $\Omega^1_{\scriptscriptstyle univ}
\rightarrow \Omega^1$. 
\medskip {\it Preuve} de 2.4.1: abr\'egeons provisoirement l'expression $D(am)-aD(m)$ (resp.
$D(mb)-D(m)b$) en $a\bullet m$ (resp. $m \bullet b$). Puisque $D$ est un op\'erateur d'ordre
un, on a 
\medskip\centerline{$D(amb)-D(am)b-aD(mb)+aD(m)b=0$, c'est-\`a-dire } 
\medskip\centerline{$D(amb)=D(am)b+(a\bullet m)b+a(m \bullet b)$.}
\medskip Rappelons que $\Omega^1_{\scriptscriptstyle univ}\cong A/k\otimes_k A$ comme $A$-module \`a droite.
 Puisque l'application $m\mapsto (m \bullet b)$ est $A$-lin\'eaire \`a gauche et nulle pour $b\in k$,
on peut d\'efinir un homomorphisme de bimodules $\psi_d(D): M\otimes_k 
\Omega^1_{\scriptscriptstyle univ} \rightarrow N$ en posant $\psi_d(D)(m\otimes d_{\scriptscriptstyle univ}b\otimes
c) = (m \bullet b)c$. Pour tout $a\in A$, on a $\psi_d(D)(m\otimes a\; d_{\scriptscriptstyle univ}b\otimes c) =
\psi_d(D)(m\otimes d_{\scriptscriptstyle univ}(ab)\otimes c - d_{\scriptscriptstyle univ}a\otimes bc) = (m \bullet
ab)c -
 (m \bullet a)bc = (ma \bullet b)c = \psi_d(D)(ma\otimes d_{\scriptscriptstyle univ}b\otimes c)$. Donc $\psi_d(D)$ 
 d\'efinit bien un homomorphisme de bimodules $\varphi_d(D): M\otimes_A \Omega^1_{\scriptscriptstyle univ}
\rightarrow N$. On proc\`ede de m\^eme de l'autre c\^ot\'e pour d\'efinir $\varphi_s(D)$, et on a
bien $D(amb)=D(am)b+\varphi_s(D)( d_{\scriptscriptstyle univ}a\otimes m)b+a\varphi_d(D)(m\otimes
d_{\scriptscriptstyle univ}b).$ L'unicit\'e de $\varphi_s(D)$ (resp. $\varphi_d(D)$) se voit en posant $b=1$
(resp. $a=1$) dans cette formule).

\bigskip A la suite du travail de J. Mourad [Mou95] sur l'adaptation en g\'eom\'etrie
non-commutative du concept de connexion lin\'eaire en g\'eom\'etrie diff\'erentielle (connexion
sur le fibr\'e tangent), M. Dubois-Violette et T. Masson [DubM96] ont d\'egag\'e la notion de
``bimodule $\Omega$-connection". Nous reprendrons cette notion dans une tout autre perspective,
en la rebaptisant simplement biconnexion: 
 \medskip \noindent {\bf Definition 2.4.2.} On suppose que l'anneau diff\'erentiel $(A
{\buildrel{d}\over \rightarrow}\Omega^1)$ est r\'eduit, i.e. $\Omega^1=dA.A$ (ainsi $\pi$ est
surjective). Une {\it biconnexion} (\`a gauche) est une connexion $\nabla$ pour laquelle
$\varphi_{\scriptstyle d}(\nabla)$ {\it se factorise \`a travers $id_M \otimes \pi$}. Il donne
alors lieu \`a un homomorphisme canonique de bimodules 
\medskip\centerline{${\phi}(\nabla)={\phi}_{\scriptstyle
D}(\nabla): M\otimes_{\scriptstyle A}\Omega^1\;\longrightarrow\;\Omega^1\otimes_{\scriptstyle A}
\;M$} 
\medskip \noindent que nous appellerons la {\it volte} (droite-gauche) de la biconnexion $\nabla$.
\medskip \noindent {\bf 2.4.3. Exemples.} $i)$ La connexion triviale $d: A \rightarrow
\Omega^1$: la volte est l'unique prolongement bilin\'eaire de l'application $1_{\scriptstyle A}\otimes
\omega \rightarrow \omega \otimes 1_{\scriptstyle A}$ (en particulier, c'est un isomorphisme). 
\par \noindent $ii)$ Le cas classique o\`u $\Omega^{\ast}$ est gradu\'ee commutative, et $M$
est un bimodule commutatif; la volte est alors l'isomorphisme na\"{\i}f d'\'echange des
facteurs
 $M\otimes_{\scriptstyle A}\Omega^1\;\longrightarrow\;\Omega^1\otimes_{\scriptstyle A}
\;M$.
\par \noindent $iii)$ Si $(A',d')$ est une extension diff\'erentielle de $(A,d)$ (cf. 1.2.5),
on peut consid\'erer $(A',d')$ comme un $A$-$A$-bimodule \`a biconnexion; sa volte est l'unique
prolongement  bilin\'eaire de l'application $1_{A'}\otimes \omega
\rightarrow \omega \otimes 1_{A'}$.  
\par \noindent $iv)$ Syst\`emes aux diff\'erences (I.4.3): $A=k(z),\; \sigma=$
endomorphisme de $A$ avec $z^{\sigma}\neq z$. On peut \'ecrire toute connexion sous la forme 
$\nabla(m)=dz\otimes {\Phi-id\over z^{\sigma}-z}\;(m),$ o\`u $\Phi: M\rightarrow
M$ est $\sigma$-lin\'eaire. C'est toujours une biconnexion, dont la volte s'\'ecrit  
\medskip\centerline{$\phi(\nabla)(m\otimes da)=da\otimes \Phi(m).$}
\medskip \noindent  Dans le cas trivial $M=A,\Phi=\sigma$, on observe que c'est bien un isomorphisme,
m\^eme si $\;\sigma$ n'est pas injectif. 
 \par Revenons au cas g\'en\'eral.    
\proclaim Lemme 2.4.4. Soit $f: (M_1,\nabla_1) \rightarrow (M_2,\nabla_2)$ un morphisme de
biconnexions (i.e. homomorphisme de
$A$-$A$-bimodules horizontal, $\nabla_1$ et $\nabla_2$ \'etant des biconnexions). On a 
$\;{\phi}(\nabla_2)(f\otimes id_{\scriptstyle \Omega^1})=(id_{\scriptstyle
\Omega^1} \otimes f){\phi}(\nabla_1).$
\par \noindent Cela r\'esulte imm\'ediatement de la d\'efinition $\varphi_{\scriptstyle
d}(\nabla)(m\otimes db.a) = D(mb)a-D(m)ba$. 

\proclaim Lemme 2.4.5. Soit $(M',\nabla'){\buildrel{f}\over\rightarrow}
(M,\nabla){\buildrel{g}\over\rightarrow} (M",\nabla")$ une suite exacte de $A$-$A$-bimodules \`a
connexion. On suppose que $\nabla$ est une biconnexion. Alors
\hfill \break $i)$ si $g$ est surjective, $\nabla"$ est une
biconnexion. De plus si $\nabla'$ est une biconnexion, et si les voltes de $\nabla$ et $\nabla'$ sont
inversibles, il en est de m\^eme de $\nabla"$.
\hfill \break $ii)$ Si $f$ est injective, et fait du bimodule $M'$ un facteur direct de $M$, alors
$\nabla'$ est une biconnexion. De plus si $\nabla"$ est une biconnexion, et si les voltes de $\nabla$ et
$\nabla"$ sont inversibles, il en est de m\^eme de $\nabla'$.
\hfill \break $iii)$ Si $\Omega^1$ est plat \`a droite sur $A$, les $A$-$A$-bimodules \`a biconnexion
forment une cat\'egorie ab\'elienne $k$-lin\'eaire. 
\hfill \break $iv)$ Si $\Omega^1$ est plat aussi \`a gauche, alors la sous-cat\'egorie pleine form\'ee des
biconnexions \`a volte inversible est stable par sous-quotients.
\par {\it Preuve.} Via $\varphi_{\scriptstyle d}\;$, on obtient un morphisme de suites exactes de
bimodules
\medskip \centerline{$\matrix{M' \otimes\Omega^1_{\scriptscriptstyle univ}
& \rightarrow  &
M \otimes\Omega^1_{\scriptscriptstyle univ}
&
\rightarrow & M" \otimes\Omega^1_{\scriptscriptstyle univ}
\cr 
 \downarrow && \downarrow && \downarrow 
\cr  \Omega^1\otimes M' & \rightarrow  &
  \Omega^1 \otimes M & \rightarrow &  \Omega^1\otimes M"}$}
\medskip \noindent Si $g$ est surjective, il se compl\`ete en 
\medskip \centerline{$\matrix{M' \otimes\Omega^1_{\scriptscriptstyle univ}
& \rightarrow  & M \otimes\Omega^1_{\scriptscriptstyle univ}
& \rightarrow & M" \otimes\Omega^1_{\scriptscriptstyle univ}
& \rightarrow  & 0 \cr 
 \downarrow && \downarrow && \downarrow && \downarrow 
\cr  \Omega^1\otimes M' & \rightarrow  &
  \Omega^1 \otimes M & \rightarrow &  \Omega^1\otimes M" & \rightarrow  & 0}$}
\medskip \noindent Comme $\nabla$ est une biconnexion, ce morphisme se factorise par
$M' \otimes\Omega^1_{\scriptscriptstyle univ}
\rightarrow  
M \otimes\Omega^1 \rightarrow M" \otimes\Omega^1 \rightarrow 0$ (resp. par 
$\; M' \otimes\Omega^1 \rightarrow M \otimes\Omega^1 \rightarrow M" \otimes\Omega^1 \rightarrow
0\;$ si $\Omega^1$ est plat \`a droite et $f$ est injective, puisqu'alors $0\rightarrow  \Omega^1\otimes
M' \rightarrow  \Omega^1 \otimes M  \rightarrow   \Omega^1\otimes M" \rightarrow  0$ est exacte). Une
chasse au diagramme ais\'ee prouve alors $i)$ et $iii)$. 
\medskip\noindent Si $f$ est injective, et si $f$ fait du bimodule $M'$ un facteur direct de $M$ (ou si
$\Omega^1$ est plat \`a droite et \`a gauche), on a un morphisme de suites exactes  
\medskip \centerline{$\matrix{0 & \rightarrow & M' \otimes\Omega^1_{\scriptscriptstyle univ}
& \rightarrow  & M \otimes\Omega^1_{\scriptscriptstyle univ}
& \rightarrow & M" \otimes\Omega^1_{\scriptscriptstyle univ}
\cr  \downarrow && \downarrow && \downarrow && \downarrow 
\cr 0 & \rightarrow & \Omega^1\otimes M' & \rightarrow  &
  \Omega^1 \otimes M & \rightarrow &  \Omega^1\otimes M" }$}
\medskip \noindent Comme $\nabla$ est une biconnexion, on voit imm\'ediatement qu'il en est de m\^eme de 
$\nabla'$, et que si $\nabla"$ est une biconnexion, le morphisme de suites exactes pr\'ec\'edent se
factorise par 
$0 \rightarrow M' \otimes\Omega^1
\rightarrow  
M \otimes\Omega^1 \rightarrow M" \otimes\Omega^1 $. Une chasse au
diagramme ais\'ee prouve alors $ii)$ et $iv)$. 
\medskip\noindent {\bf 2.4.6. Remarque.} Voici un l\'eger raffinement, utile plus bas (4.5.5). Supposons
que $\nabla, \nabla', \nabla"$ sont des biconnexions. Supposons $\Omega^1$ plat \`a droite sur $A$,
$f$ injective, $g$ surjective. On alors un morphisme de suites exactes   
\medskip \centerline{$\matrix{&& M' \otimes\Omega^1
& \rightarrow  & M \otimes\Omega^1
& \rightarrow & M" \otimes\Omega^1 &\rightarrow &0
\cr  && \downarrow && \downarrow && \downarrow  && \downarrow
\cr 0 & \rightarrow & \Omega^1\otimes M' & \rightarrow  &
  \Omega^1 \otimes M & \rightarrow &  \Omega^1\otimes M" &\rightarrow &0}$}
\medskip \noindent Supposons $\phi(\nabla)$ inversible. Alors pour que  $\phi(\nabla")$ soit inversible,
il faut et il suffit que l'homomorphisme $Im(M' \otimes\Omega^1\rightarrow M\otimes\Omega^1)\rightarrow
\Omega^1\otimes M'$ induit par $\phi(\nabla')$ soit inversible.

\bigskip {\bf 3. Produit tensoriel.} 

\medskip {\bf 3.1. D\'efinition.} Dans le cas classique commutatif, la connexion canonique sur le
produit tensoriel de deux modules \`a connexion est donn\'ee par la r\`egle bien connue
$\nabla_1\otimes id + id \otimes \nabla_2$, l'isomorphisme na\"{\i}f d'\'echange
$M_1\otimes_{\scriptstyle A}\Omega^1\;\rightarrow\;\Omega^1\otimes_{\scriptstyle A} M_1$ \'etant sous-entendu dans le
second terme $id \otimes \nabla_2$. 
\par \noindent Dans le cas o\`u $\Omega^1$ est un bimodule non-commutatif, il n'y a plus
d'isomorphisme d'\'echange et la formule perd son sens. On peut toutefois la ``sauver", dans
le cas des biconnexions, en rempla\c cant l'isomorphisme na\"{\i}f d'\'echange par
${\phi}(\nabla_1)$.
\par Soient donc $(M_1,\nabla_1)$ et $(M_2,\nabla_2)$ deux $A$-$A$-bimodules \`a biconnexion. 

\proclaim Lemme 3.1.1. Il existe une (unique) connexion (\`a gauche) $\nabla$ sur $M_1\otimes_{\scriptstyle A}
M_2$ telle que $\;\;\;\nabla(m_1\otimes m_2)=\nabla_1(m_1)\otimes m_2 +
({\phi}(\nabla_1)\otimes id_2) (m_1\otimes \nabla_2(m_2))$
\medskip \noindent (en d'autres termes, $\nabla = \nabla_1\otimes id_2+({\phi}(\nabla_1)\otimes
id_2)\circ(id_1 \otimes \nabla_2)$). 
\medskip \noindent C'est une biconnexion, dont la volte est 
$\;\;\;{\phi}(\nabla)=({\phi}(\nabla_1)\otimes
id_2)\circ(id_1 \otimes {\phi}(\nabla_2)).$ 
\par {\it Preuve}. Pour prouver la premi\`ere assertion, il est commode d'utiliser la d\'efinition
des connexions en termes de jets: consid\'erons les homomorphismes de $A$-module \`a gauche  
\par \centerline{$D_1=(id,\nabla_1): M_1 \rightarrow {\cal P}^1(M_1)\;\;\;{\rm et}\;\;\;{\cal
P}^1(M_1)\otimes M_2\;{\buildrel{\cong}\over\rightarrow}\; {\cal P}^1(M_1\otimes M_2).$}
\medskip \noindent Alors $D=(id,\nabla)$ n'est autre que la composition 
\medskip \centerline{$M_1\otimes M_2 \;{\buildrel{(D_1,\nabla_2)}\over\rightarrow}\; ({\cal
P}^1(M_1)\otimes M_2)
\oplus  (M_1 \otimes \Omega^1 \otimes M_2)\;{\buildrel{(id,\;
{\phi}(\nabla_1)\otimes id_2)}\over\rightarrow}\;$}
\medskip \centerline{$\rightarrow({\cal P}^1(M_1)\otimes M_2)
\oplus(\Omega^1\otimes M_1\otimes M_2) \;{\buildrel{+}\over\rightarrow}\;{\cal
P}^1(M_1)\otimes M_2\;{\buildrel{\cong}\over\rightarrow}\; {\cal P}^1(M_1\otimes M_2)$}
\medskip \noindent qui est clairement une section $A$-lin\'eaire de ${\cal P}^1(M_1\otimes M_2)
\rightarrow M_1\otimes M_2$. On obtient donc une connexion bien d\'efinie $\nabla$. La seconde assertion
r\'esulte du calcul 
\medskip \centerline{$\nabla((m_1\otimes m_2)a)-\nabla(m_1\otimes m_2).a=
\nabla_1(m_1)\otimes m_2a + 
({\phi}(\nabla_1)\otimes id_2) (m_1\otimes \nabla_2(m_2a)) $}
\medskip \centerline{$- \nabla_1(m_1)\otimes m_2a - ({\phi}(\nabla_1)\otimes id_2) (m_1\otimes
\nabla_2(m_2)a)$}
\medskip \noindent o\`u deux termes se compensent.

\proclaim Lemme 3.1.2. Le produit tensoriel ainsi d\'efini est un bifoncteur sur les $A$-$A$-bimodules \`a
biconnexion.
\par Cela r\'esulte imm\'ediatement de la d\'efinition du produit tensoriel et de 2.4.4. 

\medskip \noindent {\bf 3.1.3. Exemples.} $\bullet$ Dans la situation de 2.4.3.$iii)$, et si
$(M,\nabla)$ est un $A$-$A$-bimodule \`a biconnexion, l'image inverse $u^{\ast}(M,\nabla)$
peut \^etre vue comme produit tensoriel de biconnexions
$(A',d')\otimes (M,\nabla)$.
\par \noindent $\bullet$ Syst\`emes aux diff\'erences (I.4.3): si
$\nabla_1$ et
$\nabla_2$ sont deux connexions correspondant aux endomorphismes $\sigma$-lin\'eaires $\Phi_1$ et
$\Phi_2$ respectivement, leur produit tensoriel est la connexion correspondant \`a $\Phi_1
\otimes \Phi_2$.
\par J'ignore si, en g\'en\'eral, le produit tensoriel de deux connexions int\'egrables
est int\'egrable. 

\medskip {\bf 3.2. Contraintes d'unit\'e et d'associativit\'e}. Ce sont celles induites
par celles des $A$-$A$-bimodules. 
\medskip $\bullet$ L'unit\'e est la connexion triviale ${\bf 1}=(A,d)$. On a
$End\;{\bf 1}\;= C \cap Z(A)\supset k,\;$ o\`u $Z(A)$ d\'esigne le centre de $A$. 
\medskip \noindent Le $End\;{\bf 1}$-$End\;{\bf 1}$-bimodule sous-jacent
\`a $\Omega^1$ est commutatif, en vertu du calcul $c.da=d(ca)=d(ac)=da.c, c \in C \cap Z(A)$.
\par $\bullet$ Pour l'associativit\'e, on remarque que $\nabla_1\otimes id_{23}+
({\phi}(\nabla_1)\otimes id_{23})(id_1\otimes (\nabla_2\otimes id_3+
({\phi}(\nabla_2)\otimes id_3)(id_2\otimes \nabla_3))=\nabla_1\otimes
id_{23}+({\phi}(\nabla_1)\otimes id_2)(id_1\otimes \nabla_2)\otimes
id_3)+({\phi}(\nabla_{12})\otimes id_3)(id_{12}\otimes \nabla_3),$ compte tenu de
${\phi}(\nabla_{12})=({\phi}(\nabla_1)\otimes id_2)\circ(id_1 \otimes
{\phi}(\nabla_2)).$
\medskip On r\'esume ces r\'esultats comme suit:
\proclaim Proposition 3.2.1. Les $A$-$A$-bimodules \`a biconnexion (\`a gauche)
forment une cat\'egorie $k$-lin\'eaire mono\"{\i}dale.

\medskip {\bf 3.3. Dualit\'e.} Ce point est un peu plus d\'elicat. Il y a lieu de distinguer
entre duaux \`a gauche et \`a droite (cf. [Bru94] pour une discussion g\'en\'erale des duaux). 
\medskip\noindent {\bf 3.3.1. Pr\'eduaux.} Soit $(M_{\scriptstyle d},\nabla_{\scriptstyle d})$ un
$A$-module \`a connexion {\it \`a droite} avec $M$ projectif de type fini. Le $A$-module \`a gauche
${}^{\ast}M_{\scriptstyle s}=Hom_A(M_{\scriptstyle d},A)$ est muni d'une connexion (\`a gauche)
${}^{\ast}\nabla_{\scriptstyle s}$ de la mani\`ere suivante. Identifiant les $A$-modules \`a gauche 
$\Omega^1\otimes{{}^{\ast}M}_{\scriptstyle s}$ et $Hom^A(M_{\scriptstyle d},\Omega^1_{\scriptstyle d})$
(homomorphismes de $A$-modules \`a droite), on d\'efinit ${}^{\ast}{\nabla}_{\scriptstyle
s}({}^{\ast}{m})\in
\Omega^1\otimes{{}^{\ast} M}_{\scriptstyle s}$ par   
\medskip\centerline{$\langle {}^{\ast}{\nabla}_{\scriptstyle s}({}^{\ast}{m}),m\rangle  =d\langle
{}^{\ast}{m},m\rangle  -\langle {}^{\ast}{m},\nabla_{\scriptstyle d}(m)\rangle   \; \in\Omega^1.$}     
\medskip\noindent Plus pr\'ecis\'ement, $(id_{\Omega^1}\otimes \epsilon )(
{}^{\ast}{\nabla}_{\scriptstyle s}({}^{\ast}{m})\otimes m ) = d(
\epsilon ({}^{\ast}{m}\otimes m)) -(\epsilon  \otimes id_{\Omega^1})
({}^{\ast}{m}\otimes \nabla_{\scriptstyle d}(m)),$ o\`u $\epsilon$ est l'homomorphisme d'\'evaluation.
\par \noindent On v\'erifie imm\'ediatement que c'est une connexion \`a gauche;
$({}^{\ast}M_{\scriptstyle s},{}^{\ast}\nabla_{\scriptstyle s})$ est le pr\'edual de $(M_{\scriptstyle
d},\nabla_{\scriptstyle d})$. 
\par \noindent De plus, si $M_{\scriptstyle d}$ est sous-jacent \`a un $A$-$A$-bimodule et si
$\nabla_{\scriptstyle d}$ est une biconnexion, alors $\nabla_{\scriptstyle s}$ est une biconnexion \`a
gauche, dont la volte v\'erifie $\langle \phi_{\scriptstyle d}(\nabla_{\scriptstyle
s})({}^{\ast}{m}\otimes da),m \rangle =
\langle {}^{\ast}{m},\phi_{\scriptstyle s}(\nabla_{\scriptstyle d})(da\otimes m)\rangle $.
\medskip Partant d'un $A$-module \`a connexion \`a gauche $(M_{\scriptstyle
s},\nabla_{\scriptstyle s})$, on d\'efinit de m\^eme son pr\'edual $(M^{\ast}_{\scriptstyle
d},\nabla^{\ast}_{\scriptstyle d})$; $\nabla^{\ast}_{\scriptstyle d}$ est une connexion sur le $A$-module
\`a droite $M{}^{\ast}_{\scriptstyle d}=Hom_A(M_{\scriptstyle s},A)$.    

\medskip\noindent {\bf 3.3.2. Des biconnexions \`a gauche aux biconnexions \`a droite.} Consid\'erons \`a
pr\'esent un bimodule \`a biconnexion (\`a gauche) $(M,\nabla_{\scriptstyle s})$ dont la volte
${\phi_{\scriptstyle d}}(\nabla_{\scriptstyle s})$ est {\it inversible}. On peut alors munir $M$ d'une
connexion \`a droite
$\nabla_{\scriptstyle d}$ d\'efinie par 
\par \centerline{$\nabla_{\scriptstyle d}= ({\phi}_{\scriptstyle d}(\nabla_{\scriptstyle s}))^{-1}\circ
\nabla_{\scriptstyle s}$}
\medskip\noindent C'est bien une connexion (\`a droite):
$\nabla_{\scriptstyle d}(ma)-\nabla_{\scriptstyle d}(m)a=({\phi}_{\scriptstyle
d}(\nabla_{\scriptstyle s}))^{-1}(\nabla_{\scriptstyle s}((ma)-\nabla_{\scriptstyle s}(m)a) = m\otimes
da,\;$ et m\^eme une biconnexion (\`a droite), de volte
${\phi}_{\scriptstyle s}(\nabla_{\scriptstyle d})=({\phi}_{\scriptstyle d}(\nabla_{\scriptstyle s}))^{-1}$.
\medskip Partant d'un bimodule \`a biconnexion \`a droite $(M,\nabla_{\scriptstyle d})$ dont la volte
${\phi_{\scriptstyle s}}(\nabla_{\scriptstyle d})$ est {\it inversible}, on construit de m\^eme une
biconnexion \`a gauche $\nabla_{\scriptstyle s}$ sur $M$. 

\medskip\noindent {\bf 3.3.3. Duaux.} Soit $(M,\nabla)$ un
$A$-$A$-bimodule \`a biconnexion {\it \`a gauche}. On suppose que le $A$-module {\it \`a droite}
$M_{\scriptstyle d}$ sous-jacent \`a $M$ est {\it projectif de type fini}. Soit ${}^{\vee} M$ son dual;
c'est un $A$-$A$-bimodule, projectif de type fini \`a gauche. 
\par Une connexion \`a gauche
${}^{\vee}\nabla={}^{\vee}\nabla_{\scriptstyle s}$ sur ${}^{\vee} M$ est dite {\it duale \`a gauche} de
$\nabla$ si c'est une biconnexion, et si les homomorphismes d'\'evaluation 
\medskip\centerline{$\epsilon :\;\; {}^{\vee}M \otimes_{\scriptstyle A} M \rightarrow A \;$} 
\par\noindent et de co\'evaluation
\par\centerline{$\eta :\;\; A\rightarrow M\otimes_{\scriptstyle A} {}^{\vee}M 
\;$} 
\medskip\noindent sont horizontaux, i.e. induisent des morphismes de connexion.
\proclaim Lemme 3.3.4. Une connexion duale \`a gauche ${}^{\vee}\nabla$ existe si et seulement si la volte
$\phi(\nabla)$ est inversible. Elle est donn\'ee par la formule  
 $\;\;\;\;\;{}^{\vee}\nabla= {}^{\ast}\nabla_{\scriptstyle s}$
\hfill \break o\`u ${}^{\ast}\nabla_{\scriptstyle s}$ est la pr\'eduale de la biconnexion
\`a droite $\nabla_{\scriptstyle d} $ sur
$M$ attach\'ee \`a la biconnexion \`a gauche $\nabla$. On a donc 
\hfill \break $\langle {}^{\vee}{\nabla}(\check{m}),m\rangle  =d\langle \check{m},m\rangle 
-\langle \check{m},\phi(\nabla)^{-1}\nabla(m)\rangle   \; \in\Omega^1.$
\hfill \break Sa volte v\'erifie 
\hfill \break $\langle {\phi}({}^{\vee}{\nabla})(\check{m}\otimes da),\;m\rangle  =
\langle   \check{m},{\phi}(\nabla)^{-1}(da\otimes m)\rangle $.
\par {\it Preuve.} Supposons d'abord ${\phi}(\nabla)$ inversible, et d\'efinissons ${}^{\vee}{\nabla}$ par
la formule ci-dessus. Un calcul direct montre que c'est une biconnexion, de volte donn\'ee par la formule
indiqu\'ee. Il s'agit de faire voir que $\epsilon$ et $\eta$ sont horizontaux. Posons $\nabla(m)=\sum
\omega_j\otimes m_j$. On a:
\par \noindent 
$(id\otimes \epsilon)(\nabla_{{}^{\vee}M \otimes M }(\check{m}\otimes m)) =
(id\otimes \epsilon)({}^{\vee}\nabla (\check{m})\otimes m) + (id\otimes
\epsilon)(({\phi}({}^{\vee}{\nabla})\otimes id_M )(\check{m}\otimes
\nabla(m))) =
d\langle  \check{m},m\rangle  -
\langle   \check{m},({\phi}(\nabla))^{-1}(\sum
\omega_j\otimes m_j)\rangle   +
\langle (\sum {\phi}({}^{\vee}{\nabla})(\check{m}\otimes \omega_j)
, m_j \rangle =
d\langle   \check{m},m\rangle  .$  
\par \noindent Il s'agit d'autre part de faire voir que $\eta(1) = \sum
m_i\otimes \check{m_i}$ est horizontal. Or $\nabla_{M\otimes {}^{\vee}M }(\sum
m_i\otimes \check{m_i})=\sum \nabla(m_i)\otimes \check{m_i} +({\phi}(\nabla)\otimes id)(\sum
m_i\otimes {}^{\vee}\nabla(\check{m_i})) $. En ``\'evaluant" contre un \'el\'ement
quelconque $m \in M$, on trouve 
\par \noindent
$\sum \nabla(m_i) \langle\check{m_i},m\rangle   +
{\phi}(\nabla)(\sum
m_i\otimes d\langle\check{m_i},m\rangle  ) - {\phi}(\nabla)(\sum
m_i\otimes \langle\check{m_i},{\phi}(\nabla)^{-1}\nabla(m)\rangle    = \sum
\nabla(m_i)
\langle\check{m_i},m\rangle   + \sum
\nabla(m_i.\langle\check{m_i},m\rangle  ) - \sum\nabla(m_i).\langle\check{m_i},m\rangle   -
\nabla(m) = 0$. On conclut que la connexion ${}^{\vee}\nabla$ est bien duale \`a gauche de $\nabla$.
\par \noindent R\'eciproquement, supposons que la biconnexion
duale \`a gauche ${}^{\vee}\nabla$ existe. Soit $\sum n_i\otimes \omega_i$ un \'el\'ement de $M\otimes
\Omega^1$ tel que $\phi(\nabla)(\sum n_i\otimes \omega_i)=0$, et soit $\check m$ un \'el\'ement quelconque
de ${}^{\vee}{M}$. On a      
\par \noindent $ \sum\langle \check{m}, n_i  \rangle \omega_i= (id_{\scriptstyle \Omega^1}\otimes
\epsilon )\circ  ((\phi({}^{\vee}\nabla)\otimes id_{\scriptstyle M})\circ (id_{\scriptstyle
 {}^{\vee}M}\otimes
\phi(\nabla))(\sum \check{m} \otimes n_i\otimes \omega_i) =0$; on en tire que  $\sum n_i\otimes
\omega_i=0$, et l'injectivit\'e de $\phi(\nabla)$. 
\par \noindent D'autre part, en \'ecrivant comme ci-dessus $\eta(1) = \sum
m_i\otimes \check{m_i}$, on a pour tout $m\in M$,      
 $\omega \otimes m = \sum \langle \omega \otimes m_i\otimes \check{m_i} , m \rangle = 
 \langle (id_{\scriptstyle \Omega^1}\otimes  \eta) (\omega ), m \rangle =
 \langle (id_{\scriptstyle
\Omega^1}\otimes  \eta) (\omega ), m \rangle
= \langle (\phi(\nabla)\otimes id_{\scriptstyle {}^{\vee}M}) \circ (id_{\scriptstyle M} \otimes
\phi({}^{\vee}\nabla)) (\eta(1) \otimes \omega), m \rangle \in  \phi(\nabla)(M\otimes \Omega^1),$ 
\par \noindent ce qui montre la surjectivit\'e de $\phi(\nabla)$. Ainsi $\phi(\nabla)$ est inversible. Or
la biconnexion duale \`a gauche ${}^{\vee}\nabla$ est unique, elle est donc donn\'ee par la formule du
lemme.         
\medskip Une construction sym\'etrique associe \`a tout bimodule \`a biconnexion \`a gauche $(M,\nabla)$
tel que $M$ soit projectif de type fini comme $A$-module {\it \`a gauche} et tel que la volte du pr\'edual
$\nabla^{\ast}$ soit inversible, une biconnexion \`a gauche duale \`a droite
$\nabla^{\vee}$ sur le dual \`a droite $M^{\vee}=Hom(M_{\scriptstyle s},A)$.

\medskip \noindent {\bf Exemple 3.3.5.} Syst\`emes aux diff\'erences (I.4.3): supposons pour
simplifier que $\sigma$ soit un {\it automorphisme} de $A$. Un module \`a connexion 
$(M,\nabla)$ admet un dual si et seulement si l'endomorphisme
$\sigma$-lin\'eaire $\Phi$ de $M$ correspondant est inversible. Dans ce cas,
${}^{\vee}\nabla=\nabla^{\vee}$ correspond \`a  ${}^t\Phi^{-1}$.

\bigskip {\bf 4. La situation semi-classique.}

\medskip {\bf 4.1.} Nous appellerons {\it situation semi-classique} celle o\`u l'on
consid\`ere un anneau diff\'erentiel r\'eduit $A\rightarrow \Omega^1=dA.A$ avec $A$ {\it
commutatif}, et o\`u les $A$-modules \`a gauche munis d'une connexion sont toujours
consid\'er\'es comme $A$-$A$-bimodules commutatifs; le $A$-$A$-bimodule $\Omega^1$ n'est
toutefois pas suppos\'e commutatif. En termes imag\'es, on a un espace de base classique, la
``quantification" ne portant que sur l'espace cotangent. 

\proclaim Proposition 4.1.1. Dans la situation semi-classique, toute connexion est une
biconnexion.  
\par {\it Preuve}. Par d\'efinition de $\phi(\nabla)$, il s'agit de montrer que si $\sum_i
m_i\otimes da_i =0$, alors $\sum_i \nabla(m_i.a_i)= \nabla(m_i).a_i$. Tout $A$-module \`a
connexion \'etant quotient d'un $A$-module {\it libre} \`a
connexion, on peut, d'apr\`es le lemme 2.4.5, se limiter au cas o\`u
$M$ est libre. On peut alors remplacer la condition $\sum_i
m_i\otimes da_i =0$ par $\sum_i
b_i.da_i =0,$ avec $b_i \in A$, et il s'agit de montrer que pour tout $m\in M, \;\sum_i
\nabla(m.b_i.a_i)= \sum_i \nabla(m.b_i).a_i$. Or $\sum_i
b_i.da_i =0$ entra\^{\i}ne que pour tout $c\in A$, $\sum_i b_ia_idc= \sum_i b_id(a_ic)=\sum_i
b_id(ca_i)=\sum_i b_idc.a_i$, d'o\`u $\sum_i b_ia_i\omega=\sum_i b_i\omega a_i$ pour tout $\omega
\in \Omega^1$, et en particulier $\sum_i b_ia_i \nabla(m)=\sum_i b_i \nabla(m)a_i$.
\par \noindent En combinant ceci aux \'egalit\'es $\sum_i
\nabla(m.b_i.a_i)=\sum_i
\nabla(b_i.a_i.m)=\sum_i db_i.a_i.m + \sum_i b_i.a_i\nabla(m)$ et $\sum_i \nabla(m.b_i).a_i =
 \sum_i \nabla(b_i.m).a_i=\sum_i db_i.a_i.m +  \sum_i b_i\nabla(m).a_i$, on obtient le r\'esultat
voulu. 
\medskip {\bf 4.2. Contrainte de commutativit\'e.}
C'est celle induite par {\it l'\'echange des facteurs}. Pour v\'erifier la coh\'erence,
\'ecrivons 
\par \noindent $\nabla_1(m_1)=\sum_i da_i\otimes m_1^i,\; \nabla_2(m_2)=\sum_j db_j\otimes
m_2^j$; alors 
\par \noindent $\nabla(m_1\otimes m_2) = \sum_i da_i\otimes m_1^i\otimes (m_2 -\sum_j
b_jm_2^j)+\sum_{ij} b_j da_i\otimes m_1^i\otimes m_2^j + \sum_j db_j \otimes m_1\otimes m_2^j =
\sum_{ij} [b_j, da_i]\otimes m_1^i\otimes m_2^j + \sum_i da_i \otimes m_1^i\otimes m_2 +\sum_j
db_j \otimes m_1\otimes m_2^j$.
\par \noindent La sym\'etrie voulue r\'esulte de la formule 
 $[b_j, da_i]= [a_i, db_j]$ (qu'on obtient en d\'eveloppant $d(a_ib_j)=d(b_ja_i)$).
 \medskip \noindent Notons d'autre part que $End\;{\bf 1}=C$.
\medskip Du fait de la contrainte de commutativit\'e, les duaux \`a gauche sont aussi des duaux \`a
droite (lorsqu'ils existent). Plus pr\'ecis\'ement, l'isomorphisme canonique de
$A$-modules ${}^{\vee}M \cong M^{\vee}$ est horizontal. Compte tenu de cette identification, nous noterons
$({\check M},\check{\nabla})$ plut\^ot que $({}^{\vee}M,{}^{\vee}\nabla)$
 ou $(M^{\vee},\nabla^{\vee})$, et nous l'appellerons alors simplement ``dual" de $(M,\nabla)$ (lorsqu'il existe, cf. 3.3.3). 
\par Un objet est dit {\it rigide} s'il admet un dual\footnote{$^{(1)}$}{il s'agit de la
notion usuelle d'objet rigide dans une cat\'egorie  mono\"{\i}dale sym\'etrique. Cette notion n'a bien
entendu rien
\`a voir une notion homonyme classique en th\'eorie des \'equations diff\'erentielles
lin\'eaires, li\'ee \`a l'absence de ``param\`etre  accessoire".}.
 Il d\'ecoule alors de 3.3.4 que 

\proclaim Lemme 4.2.1. Un $A$-module \`a connexion $(M,\nabla)$ est rigide si et seulement si $M$ est
projectif de type fini et la volte $\phi(\nabla)$ est inversible. 
\par \noindent Rappelons que la volte est automatiquement inversible si $\Omega^1$ est un
bimodule commutatif (cf. 2.4.3.$ii)$). 
\medskip Le th\'eor\`eme suivant r\'esume nos r\'esultats:

\proclaim Th\'eor\`eme 4.2.2. Dans la situation semi-classique, les $A$-modules \`a connexion
forment une cat\'egorie ab\'elienne $C$-lin\'eaire mono\"{\i}dale sym\'etrique ($C$-{\rm
tensorielle, dans la terminologie de [Bru94]}). La sous-cat\'egorie form\'ee des objets
rigides est autonome. 
\par
\medskip {\bf 4.3. Puissances sym\'etriques et altern\'ees.} Compte tenu de la
contrainte de commutativit\'e, elles s'obtiennent
\`a partir du produit tensoriel par les constructions quotients standard. Dans le cas d'un
syst\`eme aux diff\'erences d\'ecrit par une matrice $\sigma$-lin\'eaire $\Phi$ inversible (cf.
I.4.3), elles correspondent aux puissances sym\'etriques et altern\'ees de $\Phi$. 
\par \noindent Examinons le cas particulier de la puissance altern\'ee $\mu$-i\`eme d'une
connexion $\nabla$ sur $A^{\mu}$ correspondant \`a une \'equation aux diff\'erences 
\medskip\centerline{$\sigma^{\mu}y+a_{\mu -1}\sigma^{\mu -1}y+\ldots +a_0y=0,\;\;a_i \in A$}
\medskip\noindent mise sous forme de syst\`eme 
\medskip\centerline{$(\ast\ast)\;\;\;\;\sigma(Y)={\cal A}Y,\;\;\;\;\;{\cal
A}=\pmatrix{0&1&&\ldots&\cr
             0&0&1&&\cr &&\ddots &&\cr -a_0 &-a_1 &&\ldots &-a_{\mu}}$. }
\medskip\noindent Soit $\vec y=(y_1,\ldots,y_{\mu})$ la premi\`ere ligne de $\cal Y$. Le
d\'eterminant de Casorati de $\vec y$ est 
\medskip\centerline{$Cas (\vec y) = det\; Y= \pmatrix{y_1 & y_2 & \ldots & y_{\mu}\cr
y_1^{\sigma}& y_2^{\sigma} &\ldots & y_{\mu}^{\sigma}\cr && \ddots & \cr
y_1^{\sigma^{\mu -1}} & y_2^{\sigma^{\mu -1}} & \ldots & y_{\mu}^{\sigma^{\mu -1}}}$.} 
\medskip\noindent Si $k$ est l'anneau des constantes (suppos\'e int\`egre), $Cas (\vec y)$
s'annule si et seulement si les $y_i$ sont lin\'eairement d\'ependants sur $k$. Il v\'erifie
l'\'equation aux diff\'erences d'ordre un
\par\centerline{$\sigma(Cas (\vec y))=(-1)^{\mu}a_0.Cas (\vec y)$}
\medskip\noindent qui correspond \`a $\Lambda^{\mu}\nabla$. 

\medskip {\bf 4.4. Hom interne.} Soient ${\cal M}'=(M',\nabla')$ et ${\cal
M}"=(M",\nabla")$ deux $A$-modules \`a connexion. On suppose que $\Omega^1$ est plat sur $A$ \`a droite,
que $M'$ est de pr\'esentation finie, et que la volte de
$\nabla'$ est inversible. On peut alors identifier les $A$-$A$-bimodules $\Omega^1\otimes
Hom_A(M',M")$ et $Hom_A(M',\Omega^1\otimes M")$ (cf. [Bou85]I.2.9).
\par \noindent On d\'efinit une connexion sur
$Hom_A(M',M")$ en posant, pour tout $f\in Hom_A(M',M")$ et tout $m\in M':$ 
\medskip \centerline{$\nabla(f)(m)=\nabla"(f(m))-\phi(\nabla")(f\otimes 1_{\scriptstyle
\Omega^1})\phi(\nabla')^{-1}(\nabla'(m)).$}
\medskip \noindent V\'erifions la r\`egle de Leibniz:
 $\nabla(af)(m)=\nabla"(af(m))-\phi(\nabla")(af\otimes 1_{\scriptstyle
\Omega^1})\phi(\nabla')^{-1}(\nabla'(m)) \break =  a\nabla"(f(m))+da\otimes
f(m)-a\phi(\nabla")(f\otimes 1_{\scriptstyle
\Omega^1})\phi(\nabla')^{-1}(\nabla'(m))=a\nabla(f)(m)+da\otimes f(m)$.
\par \noindent On note ce module \`a connexion ${\cal I}hom({\cal M}',{\cal M}")$. Lorsque ${\cal
M}"=(A,d)$, on le note aussi par abus $\check{\cal M}'=(\check{M}',\check{\nabla}')$ (m\^eme si $M'$ n'est
pas projectif de type fini). Dans le cas o\`u $M'$ est projectif de type fini, l'isomorphisme canonique de
$A$-modules
$Hom_A(M',M")\rightarrow M"\otimes \check{M}'$ est horizontal, i.e. induit un isomorphisme
${\cal I}hom({\cal M}',{\cal M}")\cong{\cal M}"\otimes \check{\cal M}'$.    
\par \noindent Il d\'ecoule de 2.4.4 que tout homomorphisme horizontal ${\cal M}'\rightarrow
{\cal M}"$ d\'efinit un \'el\'ement de ${\cal I}hom({\cal M}',{\cal M}")^{\nabla}$.  

\medskip {\bf 4.5. Changement d'anneau diff\'erentiel.}
\medskip \noindent Consid\'erons un morphisme d'anneaux diff\'erentiels r\'eduits $u=(u^0, u^1): \;(A
{\buildrel{d}\over \rightarrow}\Omega^1) \rightarrow (A'{\buildrel{d'}\over
\rightarrow}\Omega'^1)$, $A$ et $A'$ \'etant suppos\'es commutatifs. 
\proclaim Lemme 4.5.1. {\it $\Bigl\{ A-\hbox{modules \`a connexion} \Bigr\}
{\buildrel{u^{\ast}}\over \longrightarrow} \Bigl\{ A'-\hbox{modules \`a connexion}\Bigr\}$} 
\hfill\break est un foncteur mono\"{\i}dal.
\par En effet, soient $(M_1,\nabla_1)$ et $(M_2,\nabla_2)$ deux $A$-modules \`a connexion. 
Soient $a'_1\otimes m_1 \in u^{\ast}M_1,\;a'_2\otimes m_2 \in u^{\ast}M_2$. Notons $\nabla$ (resp.
$\nabla"$) la connexion produit tensoriel de $\nabla_1$ et $\nabla_2$ (resp. de
$\nabla'_1=u^{\ast}(\nabla_1)$ et $\nabla'_2=u^{\ast}(\nabla_2)$), et posons $\nabla'=u^{\ast}(\nabla)$. 
Comparons $\nabla'$ et $\nabla"$ sous les identifications
\par \centerline{$\Omega'^1\otimes_{\scriptstyle A'}u^{\ast}(M_i)=(u^1\otimes
id_i)(\Omega^1\otimes_{\scriptstyle A}M_i)$ et} 
\par \centerline{$u^{\ast}(M_1)\otimes_{\scriptstyle
A'}u^{\ast}(M_2)=u^{\ast}(M_1\otimes_{\scriptstyle A}M_2)$.}
\par \noindent On a: $\;\nabla"((a'_1\otimes
m_1)\otimes_{\scriptstyle A'} (a'_2\otimes m_2))= \nabla"(\nabla)((a'_1a'_2\otimes m_1)\otimes (1\otimes
m_2)) =\nabla'_1(a'_1a'_2\otimes m_1)\otimes m_2 + a'_1a'_2(u^1\otimes
id_{12})({\phi}(\nabla'_1)\otimes id_2) (m_1\otimes \nabla'_2(1\otimes m_2))=
a'_1a'_2(u^1\otimes
id_{12})\nabla_1(m_1)\otimes m_2 +d'(a'_1a'_2)\otimes (m_1 \otimes m_2) + a'_1a'_2(u^1\otimes
id_{12})({\phi}(\nabla_1)\otimes id_2) (m_1\otimes \nabla_2(m_2))
= a'_1a'_2.(u^1\otimes id_{12})(\nabla(m_1\otimes m_2))+d(a'_1a'_2)\otimes (m_1\otimes m_2) =
\nabla'((a'_1\otimes m_1)\otimes_{\scriptstyle A'} (a'_2\otimes m_2)) $. 
\par \noindent Par ailleurs, il est clair que $u^{\ast}$ est compatible aux contraintes d'unit\'e,
d'associativit\'e et de commutativit\'e.  
\medskip Il d\'ecoule du lemme que $u^{\ast}$ transforme $A$-modules \`a connexion rigides en $A'$-modules
\`a connexion rigides. Si l'application $(A'\otimes \Omega^1 \otimes A')
\rightarrow \Omega'^1$ induite par $u^1$ est surjective (ce qui st le cas lorsque $\Omega^1$ est un
bimodule commutatif, ou, plus g\'en\'eralement, un sesquimodule), il est facile de calculer la volte de
$u^{\ast}(\nabla)$ en fonction de la volte de $\nabla$, et de voir que $\phi(u^{\ast}(\nabla))$ est
inversible si $\phi(\nabla)$ l'est; en particulier, $u^{\ast}$ commute \`a la formation de ${\cal
I}hom$.    

\bigskip {\bf 5. Localisation et rigidit\'e.} 
\medskip {\bf 5.1. Passage \`a l'anneau total de fractions.} Soit $(A,d)$ un anneau diff\'erentiel, $A$
\'etant commutatif. On note $Q(A)$ comme en 1.3.5 l'anneau total de fractions de $A$.  
\proclaim Proposition 5.1.1. On suppose que 
\hfill \break $i)$ l'anneau diff\'erentiel $(A,d)$ est simple,
\hfill \break $ii)$ $\Omega^1=dA.A$ est fid\`ele et projectif de type fini \`a droite, 
\hfill \break $iii)$ $d(Q(A))\subset dA.Q(A)$. 
\hfill \break Soit ${\cal M}'=(M',\nabla')$ un $A$-module de pr\'esentation finie muni d'une connexion de
volte $\phi(\nabla')$ inversible. Alors pour tout $A$-module \`a connexion ${\cal M}"=(M",\nabla")$, 
l'application naturelle $\;\;{\hbox {Mor}}({\cal M}',{\cal M}")\rightarrow {\hbox {Mor}}({\cal
M}'_{\scriptstyle Q(A)},{\cal M}"_{\scriptstyle Q(A)})\;\;$
 est un isomorphisme. 
\par {\it Preuve.} Commen\c cons par quelques remarques. L'hypoth\`ese $ii)$ entra\^{\i}ne
que $\Omega^1$ est sans torsion \`a droite. L'hypoth\`ese
$d(Q(A))\subset dA.Q(A)$ \'equivaut \`a dire que l'homomorphisme naturel de
$A$-$A$-bimodules $\Omega^1 \otimes_{\scriptstyle A}Q(A) \rightarrow Q(A)\otimes_{\scriptstyle A}
\Omega^1\otimes_{\scriptstyle A}Q(A)$ est surjectif (elle est trivialement satisfaite si $\Omega^1$
est un bimodule commutatif, ou, plus g\'en\'eralement, un sesquimodule). Comme $\Omega^1$ est sans
torsion \`a droite, on a finalement $d(Q(A)).Q(A) \cong \Omega^1\otimes_{\scriptstyle A}Q(A)$, i.e.
$(Q(A), d)$ est une extension diff\'erentielle de $(A,d)$.  
\par \noindent Prouvons l'injectivit\'e de $\;{\hbox {Mor}}({\cal M}',{\cal M}")\rightarrow 
{\hbox {Mor}}({\cal M}'_{\scriptstyle Q(A)},{\cal M}"_{\scriptstyle Q(A)})$. Consid\'erons pour
cela le $A$-module \`a connexion ${\cal I}hom({\cal M}',{\cal M}")$. On a ${\hbox {Mor}}({\cal M}',{\cal
M}")\subset {\cal I}hom({\cal M}',{\cal M}")^{\nabla}$. D'autre part, comme $M'$ est de pr\'esentation
finie, l'application naturelle $Q(A)\otimes_A Hom_A(M',M") \rightarrow  Hom_{\scriptstyle
Q(A)}(M'_{\scriptstyle Q(A)},M"_{\scriptstyle Q(A)})$ est bijective. Il suffit
donc de faire voir que l'application ${\cal I}hom({\cal M}',{\cal M}")^{\nabla}\rightarrow Q(A)
\otimes_A {\cal I}hom({\cal M}',{\cal M}")$ est injective, c'est-\`a-dire que 
l'annulateur $I\subset A$ de tout \'el\'ement $f\in {\cal I}hom({\cal M}',{\cal M}")^{\nabla}$
\par \noindent $ \subset {\cal I}hom({\cal M}',{\cal M}")$ est $0$ ou $A$. Soit $D\in {}^{\vee}\Omega^1$
et $i\in I$. On a $0=\nabla_D(if)=\langle D,di \rangle f + \nabla_{D.i}(f)=\langle D,di \rangle f$, ce qui
montre que $I$ est un id\'eal diff\'erentiel (1.3.3). On conclut par la simplicit\'e de $(A,d)$.  
\par \noindent Passons \`a la surjectivit\'e. Soit $f\in {\hbox {Mor}}({\cal
M}'_{\scriptstyle Q(A)},{\cal M}"_{\scriptstyle Q(A)})$. Soit $1\otimes {\cal M}'$ l'image de ${\cal
M}'$ dans ${\cal M}'_{\scriptstyle Q(A)}$. Alors
$f(1\otimes{\cal M}')$ est un sous-$A$-module \`a connexion de ${\cal M}"_{\scriptstyle Q(A)}$ (vu comme
$(Q(A),d)\otimes {\cal M}"$), de m\^eme que l'image $1\otimes {\cal M}"$ de ${\cal M}"$, et
$f$ induit un morphisme de $A$-modules \`a connexion 
\par \centerline{${\bar f}:\;{\cal M}'\rightarrow {\cal N}:=f(1\otimes{\cal M}')/(f(1\otimes{\cal M}')\cap
(1\otimes {\cal M}"))$.}
\par \noindent Or ${\cal N}_{\scriptstyle Q(A)}=0$, donc l'image de ${\bar f}$
dans ${\hbox {Mor}}({\cal
M}'_{\scriptstyle Q(A)},{\cal N}_{\scriptstyle Q(A)})$ est nulle. D'apr\`es ce qui
pr\'ec\`ede, ceci entra\^{\i}ne que ${\bar f}=0$, donc que $f\in {\hbox {Mor}}({\cal M}',{\cal
M}")$.
\par
\proclaim Corollaire 5.1.2. Sous les hypoth\`eses de 5.1.1, le foncteur de localisation
\medskip \centerline{\it $\Bigl\{ \matrix{A-\hbox{modules \`a connexion} \cr \hbox{rigides}} \Bigr\}
\longrightarrow \Bigl\{ \matrix{Q(A)-\hbox{modules \`a connexion}\cr
\hbox{rigides}}\Bigr\}$} 
\medskip\noindent est pleinement fid\`ele.
\par Le foncteur est bien d\'efini compte tenu de la remarque suivant 4.5.1. L'assertion est une
cons\'equence imm\'ediate de 5.1.1, puisque tout $A$-module projectif de type fini est de pr\'esentation
finie.
\medskip {\bf 5.2. Le th\'eor\`eme de rigidit\'e.}
\proclaim Th\'eor\`eme 5.2.1. Ajoutons aux hypoth\`eses de 5.1.1 que $Q(A)$ est semi-simple (i.e. produit
fini de corps). Soit $(M,\nabla)$ un $A$-module de pr\'esentation finie muni d'une connexion de volte
inversible. Alors $M$ est projectif (et $(M,\nabla)$ est donc rigide).  
\par {\it Preuve.} Soit $\check{\cal M}$ le dual de ${\cal M}=(M,\nabla)$ au sens de 4.4. On dispose du
morphisme d'\'evaluation $\;\epsilon :\;\; \check{\cal M} \otimes {\cal M} \rightarrow (A,d)$. 
\par \noindent D'autre part, comme $Q(A)$ est suppos\'e semi-simple, tout $Q(A)$-module est projectif.
D'apr\`es la remarque suivant 4.5.1, la volte de ${\cal M}_{\scriptstyle Q(A)}$ est inversible. On en
conclut que ${\cal M}_{\scriptstyle Q(A)}$ est rigide et que $(\check{\cal M})_{\scriptstyle
Q(A)}=\check{{\cal M}_{\scriptstyle Q(A)}}$. En particulier, on dispose de la co\'evaluation
$\eta_{\scriptstyle Q(A)} : (Q(A),d)\rightarrow {\cal M}_{\scriptstyle Q(A)} \otimes \check{\cal
M}_{\scriptstyle Q(A)}$ v\'erifiant $\;(id\otimes \epsilon_{\scriptstyle
Q(A)})(\eta_{\scriptstyle Q(A)}\otimes id)=id_{M_{\scriptstyle Q(A)}}$.
\par \noindent Selon 5.1.1, $\eta_{\scriptstyle Q(A)}$ provient d'un homomorphisme 
$\eta :\;\; (A,d)\rightarrow {\cal M} \otimes \check{\cal M}$ v\'erifiant $\;(id\otimes
\epsilon)(\eta\otimes id)=id_M$. D'apr\`es le lemme classique de la base duale, ceci entra\^{\i}ne que $M$
est projectif. 
\medskip \noindent Consid\'erer des anneaux $A$ avec $Q(A)$ semi-simple plut\^ot que des anneaux int\`egres
n'est pas un luxe: de tels anneaux (diff\'erentiels) sont in\'evitables dans la th\'eorie de
Picard-Vessiot pour les \'equations aux diff\'erences (cf. [vdPS97]1).

\par \proclaim Corollaire 5.2.2. Soit $X$ une vari\'et\'e alg\'ebrique lisse sur un corps $k$ de
caract\'eristique nulle ou une vari\'et\'e analytique lisse (r\'eelle ou complexe, voire $p$-adique
rigide). Tout ${\cal O}_X$-module coh\'erent muni d'une connexion ({\rm non n\'ecessairement
int\'egrable}) est localement libre. \par
 Il est bien connu que, r\'eciproquement, tout ${\cal O}_X$-module localement libre de type fini
peut \^etre muni d'une connexion.
\medskip {\bf 5.3. Le th\'eor\`eme tannakien.}
\proclaim Proposition 5.3.1. Ajoutons aux hypoth\`eses de 5.2.1 que $\Omega^1
\otimes Q(A) \cong Q(A)\otimes \Omega^1$. Soit ${\cal M}=(M,\nabla)$ un $A$-module \`a connexion
rigide et ${\cal M}'=(M',\nabla')$ un sous-objet. Soit ${\cal M}"=(M",\nabla")$ le sous-objet de ${\cal
M}$ d\'efini par $M"=M \cap {M'}_{\scriptstyle
Q(A)}\subset M_{\scriptstyle Q(A)}$, et supposons $M"$ de pr\'esentation finie. Alors $M'=M"$,
et ${\cal M}'$ est rigide.     
\par {\it Preuve.} Ecrivons $Q(A)=\Pi_1^r K_i$ comme produit fini de corps. La donn\'ee d'un
$Q(A)$-module de type fini \'equivaut \`a celle d'une famille d'espaces vectoriels $V_i$ de dimension finie
(un sur chaque $K_i$). On peut donc d\'efinir sa dimension comme le $r$-uplet des dimensions des
$V_i$.
\par \noindent L'hypoth\`ese que $\Omega^1 \otimes_{\scriptstyle A}Q(A) \cong Q(A)\otimes_{\scriptstyle A}
\Omega^1$ est plus forte que $d(Q(A))\subset dA.Q(A)$ (elle est trivialement satisfaite si $\Omega^1$
est un sesquimodule relatif \`a un automorphisme $\sigma$); elle implique que pour tout $A$-module de type
fini $N$, $\Omega^1\otimes_{\scriptstyle A}(Q(A)\otimes_{\scriptstyle A}N)$ et $(Q(A)\otimes_{\scriptstyle
A}N)\otimes_{\scriptstyle A}\Omega^1$ sont des $Q(A)$-modules \`a droite de m\^eme dimension
(finie).    
\par \noindent Par ailleurs, on a $M'\subset M"$ et ${\cal M}'_{\scriptstyle
Q(A)}={\cal M}"_{\scriptstyle Q(A)}$. Puisque $\phi(\nabla)$ est inversible, il en est de m\^eme de
$\phi(\nabla_{\scriptstyle Q(A)})$, donc $\phi(\nabla"_{\scriptstyle
Q(A)})$ est injective. Comme application entre $Q(A)$-modules \`a droite de m\^eme
dimension, elle est donc inversible. Puisque $\Omega^1$ est plat \`a droite, on a $\Omega^1\otimes M" \cong
(\Omega^1\otimes M) \cap (\Omega^1\otimes {M"}_{\scriptstyle Q(A)})$; puisque $M$ est fid\`ele et
projectif, $M \otimes \Omega^1\subset M_{\scriptstyle Q(A)}\otimes \Omega^1$, et
$Im(M" \otimes \Omega^1\rightarrow M \otimes \Omega^1) \cong ( M \otimes \Omega^1) \cap ({M"}_{\scriptstyle
Q(A)}\otimes\Omega^1)$. On d\'eduit alors de ce qui pr\'ec\`ede que $\phi(\nabla")$ se
factorise \`a travers un isomorphisme $\phi": \;Im(M" \otimes
\Omega^1\rightarrow M \otimes \Omega^1)\rightarrow  \Omega^1\otimes M"$. D'apr\`es 2.4.6, ceci
entra\^{\i}ne que la volte de ${\cal M}/{\cal M}"$ est inversible. Comme $M$ est projectif
de type fini et $M"$ de type fini, $M/M"$ est de pr\'esentation finie, et on conclut de 5.2.1 que
${\cal M}"$ est rigide. En particulier, $M"$ est facteur direct de $M$, et 2.4.5.$ii)$ entra\^{\i}ne
que $\phi(\nabla")$ est inversible. Comme $M"$ est suppos\'e de pr\'esentation finie, 5.2.1
implique derechef que ${\cal M}"$ est rigide. En appliquant 5.1.1 \`a l'isomorphisme ${\cal
M}"_{\scriptstyle Q(A)}\rightarrow {\cal M}'_{\scriptstyle
Q(A)}$, on conclut que ${\cal M}"={\cal M}'$.      
\par
\proclaim Th\'eor\`eme 5.3.2. On suppose que 
\hfill \break $i)$ l'anneau commutatif $\;A\;$ est noeth\'erien, d'anneau total de fractions 
$Q(A)$ semi-simple, 
\hfill \break $ii)$ $\Omega^1=dA.A$ est fid\`ele et projectif de type fini \`a droite, et $\Omega^1
\otimes_{\scriptstyle A}Q(A) \cong Q(A)\otimes_{\scriptstyle A}
\Omega^1$, 
\hfill \break $iii)$ l'anneau diff\'erentiel $(A,d)$ est simple. 
\hfill \break Alors tout sous-quotient d'un $A$-module \`a connexion rigide est rigide. La cat\'egorie des
$A$-modules \`a connexion rigides est tannakienne sur le corps des constantes $C$. 
\par {\it Preuve.} Pour la premi\`ere assertion, le point est, en vertu de 5.2.1, de faire voir
que tout sous-quotient d'un objet rigide est \`a volte inversible; 2.4.5.$i)$ ram\`ene la question au cas
d'un sous-objet. L'assertion r\'esulte alors de la proposition pr\'ec\'edente. De ceci et de 4.2.2, il
d\'ecoule que la cat\'egorie des $A$-modules \`a connexion rigides est tensorielle sur $C$ et autonome.
Pour montrer qu'elle est tannakienne, il suffit d'exhiber un foncteur fibre (i.e. un foncteur mono\"{\i}dal
$C$-lin\'eaire exact et fid\`ele) \`a valeurs dans les $B$-modules projectifs de type fini, pour une
$C$-alg\`ebre $B$ fid\`ele convenable ([De90]). Le foncteur d'oubli de la connexion convient avec $B=A$, d'apr\`es
5.2.1.    
\bigskip
\vfill \eject
 
\centerline{\big{ {\bf \S { III}} {\bf Groupes de Galois diff\'erentiels et extensions de
Picard-Vessiot}}}
\centerline{{\bf en situation semi-classique.}}

\medskip {\it Dans tout le reste de l'article, nous nous pla\c cons dans la situation
semi-classique: on travaille avec un anneau diff\'erentiel $(A {\buildrel{d}\over
\rightarrow}\Omega^1)$ o\`u $A$ est un anneau commutatif, mais o\`u le bimodule des
$1$-formes diff\'erentielles $\Omega^1=A.dA.A$ n'est pas n\'ecessairement commutatif}.

\bigskip {\bf 1. Solubilit\'e.}

\medskip {\bf 1.1. Connexions triviales.}
\medskip \noindent Rappelons qu'un $A$-module \`a connexion $(M,\nabla)$ est dit trivial
(II.2.3.3) s'il est isomorphe \`a un module \`a connexion de la forme $(A\otimes_{\scriptstyle C} N, d\otimes
id_N)$, o\`u $N$ est un module
\`a gauche quelconque sur l'anneau des constantes $C$. Remarquons que sa volte est alors induite 
par l'application $1_N\otimes \omega \rightarrow \omega\ \otimes 1_N$ (compte tenu de ce que le
$C$-$C$-bimodule sous-jacent \`a
$\Omega^1$ est commutatif, cf. II.3.2); c'est donc un isomorphisme.
\medskip \noindent La notion de connexion triviale est toutefois plus subtile qu'il ne
para\^{\i}t: 
\par \noindent $\bullet$ la fl\`eche naturelle $N\rightarrow M \cong A\otimes_{\scriptstyle C} N$ se factorise
\`a travers
$M^{\nabla}:= Ker_M(\nabla)$, mais $N\rightarrow M^{\nabla}$ n'est ni injective ni surjective en
g\'en\'eral. 
\par \noindent Elle est injective si $A$ est fid\`element plat sur $C$ (ce qui \'equivaut
\`a dire que $dA$ est un $C$-module plat, cf. [Bou85]1.3.5). Elle est surjective si $N$ est plat
sur $C$ (si $A$ est fid\`element plat
sur $C$, ceci \'equivaut \`a dire que $M$ est un $A$-module plat). L'exemple $A=\Omega^1={\bf
Z}[z]$ muni de $d/dz$, $N = {\bf Z}/p{\bf Z} \;( \; (A\otimes_{\scriptstyle C} N)^{\nabla}= {\bf Z}/p{\bf
Z}[z^p])$ montre qu'on n'a pas toujours surjectivit\'e. Plus pr\'ecis\'ement, la chasse au
diagramme de suite exactes
$$\matrix{ & & N  & \rightarrow  &
A\otimes_{\scriptstyle C} N &
\rightarrow & dA\otimes_{\scriptstyle C} N & \rightarrow  & 0 \cr 
 && \downarrow && \downarrow \cong  && \downarrow && 
\cr 0 &\rightarrow & M^{\nabla} & \rightarrow  &
M &
\rightarrow & \Omega^1 \otimes_{\scriptstyle A} M \cong \Omega^1 \otimes_{\scriptstyle C} N & & }$$
montre que $N\rightarrow M^{\nabla}$ est surjective si et seulement si $dA\otimes_{\scriptstyle C} N \rightarrow
\Omega^1 \otimes_{\scriptstyle A} M $ est injective.
\par \noindent $\bullet$ Des exemples de m\^eme farine montrent que $$Hom((A\otimes_{\scriptstyle C} N, d\otimes
id_N),(A\otimes_{\scriptstyle C} N', d\otimes
id_{N'}))$$ n'est en g\'en\'eral \'egal ni \`a $Hom_{\scriptstyle C}(N,N')$, ni \`a $Hom_{\scriptstyle C}((A\otimes_{\scriptstyle C}
N)^{\nabla},(A\otimes_{\scriptstyle C} N')^{\nabla'})$.
\par \noindent $\bullet$ Le conoyau d'un morphisme de connexions triviales n'est pas 
n\'ecessairement une connexion triviale: dans l'exemple pr\'ec\'edent, consid\'erer le morphisme ${\bf
Z}[z]\rightarrow {\bf Z}/p{\bf Z}[z]$ donn\'e par $1\mapsto z^p$.
\par \noindent $\bullet$ Pour toute connexion triviale
$(M,\nabla)$, l'homomorphisme naturel de
$A$-modules (\`a connexion, si l'on veut) 
$$A\otimes_{\scriptstyle C}\; M^{\nabla}\;\rightarrow M\;$$ est surjectif (en fait, si $A$ est
fid\`element plat sur $C$, $M$ est m\^eme un r\'etracte de $A\otimes_{\scriptstyle C}\; M^{\nabla}$), mais pas
n\'ecessairement injectif. 
\medskip Tous ces probl\`emes disparaissent lorsque $C$ est un corps. Toutefois, l'exemple de
$A=\Omega^1=k[z],\;d=zd/dz$ et de l'id\'eal diff\'erentiel $zA$ montre que m\^eme si $C$ est un
corps, une sous-connexion d'une connexion triviale n'est pas n\'ecessairement triviale.   

\proclaim Lemme 1.1.1. $i)$ La cat\'egorie des $A$-modules \`a connexion triviaux est stable par
$\oplus , \otimes$.
\hfill \break $ii)$ Supposons $A$ fid\`element plat sur $C$, et $M$ plat sur $A$. Alors
$(M,\nabla)$ est trivial si et seulement si $A\otimes_{\scriptstyle C}\; M^{\nabla}\;\cong M$. Dans ce cas,
$(M,\nabla)$ est rigide si et seulement si $M^{\nabla}$ est projectif de type fini sur $C$, et
alors $(\check M)^{\check \nabla}$ est le $C$-dual de $M^{\nabla}$. 
\par 
{\it Preuve}. On a 
\par \noindent $(A\otimes_{\scriptstyle C} N_1, d\otimes id_{N_1})\oplus (A\otimes_{\scriptstyle C} N_2,
d\otimes id_{N_2})=(A\otimes_{\scriptstyle C} (N_1 \oplus N_2), d\otimes id_{N_1\oplus N_2})),$ 
\par \noindent $(A\otimes_{\scriptstyle C} N_1, d\otimes id_{N_1})\otimes (A\otimes_{\scriptstyle C} N_2,
d\otimes id_{N_2})=(A\otimes_{\scriptstyle C} (N_1 \otimes_{\scriptstyle C} N_2), d\otimes id_{N_1\otimes N_2})),$ 
\par \noindent Supposons $(M,\nabla)$ trivial, $\cong (A\otimes_{\scriptstyle C} N, d\otimes
id_N)$. Sous les hypoth\`eses de $ii)$, $N$ est plat sur $C$, donc $N\cong M^{\nabla}$, et
$A\otimes_{\scriptstyle C}\; M^{\nabla}\;\cong M$. Le reste est ais\'e ($({}^\vee(A\otimes_{\scriptstyle C} N),{}^\vee(d\otimes
id_N))=(A\otimes {\check N}, d\otimes id_{\check N})$). 
\medskip \noindent {\bf 1.1.2. Remarque.} Supposons $A$ plat sur $C$. Alors $Hom_{\scriptstyle C}(N_1,N_2)$
s'envoie injectivement vers $Mor((A\otimes_{\scriptstyle C} N_1, d\otimes id_{N_1}), (A\otimes_{\scriptstyle C} N_2,
d\otimes id_{N_2}))$. Appelons {\it admissibles} les morphismes de connexions triviales qui
se trouvent dans l'image. Si on ne consid\`ere que les morphismes admissibles, la
cat\'egorie des connexions triviales qu'on obtient est ab\'elienne.     

\medskip {\bf 1.2. Crit\`eres d'injectivit\'e de $A\otimes_{\scriptstyle C} M^{\nabla}\rightarrow M$.}

\proclaim Proposition 1.2.1. Supposons $(A,d)$ simple, et $\Omega^1$ fid\`ele et projectif de type
fini \`a droite sur $A$. Alors pour tout module \`a connexion $(M,\nabla)$, l'homomorphisme
naturel $A\otimes_{\scriptstyle C} M^{\nabla}\rightarrow M$ est injectif.
\par {\it Preuve}. On a vu (II.1.3.4) que sous l'hypoth\`ese de simplicit\'e, $C$ est un corps.
Soit $\lbrace m_1,\ldots, m_i,\ldots\rbrace$ une famille finie d'\'el\'ements de $M^{\nabla}$
lin\'eairement ind\'ependants sur $C$. Supposons qu'il existe une combinaison $A$-lin\'eaire 
$\sum\;\alpha_i.m_i$ nulle dans $M$, avec $\alpha_1\neq 0$. Quitte \`a omettre
certains $m_i$, on peut supposer cette combinaison lin\'eaire de longueur minimale,
avec $\alpha_i \neq 0$ pour tout $i$. Consid\'erons alors toutes les combinaisons
$A$-lin\'eaires $\sum\;a_i.m_i$ nulles. L'ensemble des
coefficients $a_1$ intervenant dans une telle combinaison forment un id\'eal non nul $I\subset A$.
Si $I\neq A$, alors d'apr\`es II.1.3.3 et puisque
$(A,d)$ est simple, il existe $D\in {}^{\vee}\Omega^1$ et $a_1\in I$ tels que $\langle
D,da_1\rangle \notin I.$ On a $\nabla_{\scriptstyle D}(\sum\;a_i.m_i)=\sum\;\langle D,da_i \rangle.m_i = 0$ (cf.
II.2.1.2), d'o\`u $\langle D,da_i\rangle  \in I$: contradiction. Si $I = A$, on peut choisir
$a_1=1$, et $\nabla_{\scriptstyle D}(\sum\;a_i.m_i)=\sum_{i>1}\;\langle D,da_i \rangle.m_i = 0$ contredit la
minimalit\'e de la combinaison originale $\sum\;\alpha_i.m_i=0$. On conclut qu'une telle
combinaison $\sum\;\alpha_i.m_i=0$ n'existe pas.

\proclaim Proposition 1.2.2. Supposons $(A,d)$ simple par couches ({\rm cf. II.1.3.6}), $\Omega^1$
fid\`ele et projectif de type fini \`a droite sur $A$, et $C$
noeth\'erien. Alors pour tout module \`a connexion
$(M,\nabla)$ noeth\'erien (en tant que module \`a connexion), l'homomorphisme naturel $A\otimes_{\scriptstyle C}
M^{\nabla}\rightarrow M$ est injectif.  
\par {\it Preuve}. On peut supposer $M$ non nul. Par r\'ecurrence, on construit une filtration
croissante $F_n(M,\nabla)$ comme suit: 
\par \noindent $\bullet$ $F_0=0$, 
\par \noindent $\bullet$ $F_{n -1}(M,\nabla)$ \'etant suppos\'e construit, on choisit un id\'eal
{\goth p}$_n$ de
$C$ maximal parmi les annulateurs de sous-objets non nuls de $(M,\nabla)/F_{ n - 1}(M,\nabla)$ (un
tel id\'eal existe puisque $C$ est noeth\'erien); 
\par \noindent $\bullet$ on pose alors $F_n(M)=\lbrace m\in M, {\hbox {\goth p}}_n.m \subset
F_{n-1}(M)\rbrace$; il est clair que ce module est sous-jacent \`a un sous-objet $F_n(M,\nabla)$
de $(M,\nabla)$. Comme ce dernier est noeth\'erien, la filtration $F_.$ est finie 
et exhaustive.
\par \noindent Les id\'eaux {\goth p}$_n$ sont premiers. En effet, soit $(N,\nabla)$ un sous-objet
non nul de $(M,\nabla)/F_{ n - 1}(M,\nabla)$ d'annulateur {\goth p}$_n$. Soient $c_1, c_2 \in C$
tels que $c_1. c_2 \in$ {\goth p}$_n$ mais $c_2 \notin$ {\goth p}$_n$ (s'il en est). On a
$c_2.N\neq 0$, l'annulateur de $c_2.(N,\nabla)$ contient {\goth p}$_n$, donc co\"{\i}ncide avec
{\goth p}$_n$ par maximalit\'e. Puisque $c_1. c_2.(N,\nabla)=0$, on voit que $c_1 \in$ {\goth
p}$_n$, ce qui montre que {\goth p}$_n$ est premier (c'est l'argument assassin
familier). En outre, du fait de la maximalit\'e des {\goth p}$_n$, les gradu\'es associ\'es
$Gr_n(M)$ sont des $C/${\goth p}$_n$ sans torsion.
\par \noindent Tirons alors parti de ce que $(A,d)$ est
simple par couches, et notons (comme en II.1.3.6) $\kappa$({\goth p}$_n$) le corps de
fractions de $C/${\goth p}$_n$. Il d\'ecoule alors de la proposition pr\'ec\'edente appliqu\'ee
\`a $\kappa$({\goth p}$_n$)$\otimes Gr_n(M)$ que l'homomorphisme    
\medskip \centerline{$(\kappa({\hbox {\goth p}}_n)\otimes A)\otimes_{C/{\hbox {\goth p}}_n}
(Gr_n(M))^{\nabla}\hookrightarrow \kappa({\hbox {\goth p}}_n) \otimes Gr_n(M)$}
\medskip \noindent est injectif. D'o\`u l'injectivit\'e de 
\medskip \centerline{$(A/{\hbox {\goth p}}_nA)\otimes_{C/{\hbox {\goth p}}_n}
(Gr_n(M))^{\nabla}\hookrightarrow Gr_n(M).$} 
\medskip \noindent On a d'autre part un diagramme de suites exactes 
$$\matrix{& & (A/{\hbox {\goth p}}_nA)\otimes_{C/{\hbox {\goth p}}_n} (Gr_n(M))^{\nabla} & &  \cr 
          & & \downarrow \scriptstyle{\cong} & & \cr
0 &\rightarrow & A\otimes (Gr_n(M))^{\nabla} & \rightarrow  &
A\otimes (M/F_{n-1}(M))^{\nabla} &
\rightarrow & A\otimes (M/F_n(M))^{\nabla}  \cr 
 && \downarrow && \downarrow \scriptstyle{\iota}_{n-1}  && \downarrow \scriptstyle{\iota}_n  
\cr 0 &\rightarrow & Gr_n(M) & \rightarrow  &
M/F_{n-1}(M) & \rightarrow & M/F_n(M)  \rightarrow   0}$$
Une chasse \'el\'ementaire permet de conclure, par r\'ecurrence descendante, que $\iota_0:\;
A\otimes_{\scriptstyle C} M^{\nabla}\rightarrow M$ est injectif.

\proclaim Corollaire 1.2.3. Sous les hypoth\`eses de 1.2.1 (resp. 1.2.2), tout sous-quotient
d'un module \`a connexion trivial (resp. et de type fini sur $A$) est trivial. Sur la cat\'egorie
de ces connexions, le foncteur $C$-lin\'eaire $(M,\nabla)\mapsto M^{\nabla}$ est fid\`ele et
exact. 
\par {\it Preuve.} Sous les hypoth\`eses de 1.2.1 ou 1.2.2, $A$ est
fid\`element plat sur $C$ noeth\'erien; donc le $A$-module $M$ est noeth\'erien si et seulement si
$M^{\nabla}$ est de type fini sur $C$. Consid\'erons une suite exacte 
$$0\rightarrow (M',\nabla') \rightarrow (M,\nabla)\rightarrow (M",\nabla")\rightarrow 0.$$
On en d\'eduit un diagramme de suites exactes 
$$\matrix{0 &\rightarrow & A\otimes_{\scriptstyle C} M'^{\nabla'} & \rightarrow & A\otimes_{\scriptstyle C}
M^{\nabla} & \rightarrow & A\otimes_{\scriptstyle C} M"^{\nabla"} & & &  
\cr &  & \downarrow  & & \downarrow & & \downarrow & & \cr
0 &\rightarrow & M' & \rightarrow  &
M & \rightarrow &  M" & \rightarrow  & 0  \cr  }$$
D'apr\`es les propositions pr\'ec\'edentes, les fl\`eches verticales sont injectives. Si
$(M,\nabla)$ est trivial, celle du milieu est alors bijective, et on d\'eduit du lemme des cinq
qu'elles sont toutes trois bijectives. La seconde assertion en d\'ecoule. 

\medskip {\bf 1.3. Connexions solubles dans une extension diff\'erentielle.}
\medskip \noindent Soit $(A {\buildrel{d}\over
\rightarrow}\Omega^1) {\buildrel{u}\over
\longrightarrow} (A'{\buildrel{d'}\over
\rightarrow}\Omega'^1)$ une extension diff\'erentielle (II.1.2.5); rappelons que $\Omega'^1\cong
\Omega^1\otimes A'$ comme $A'$-module \`a droite. Soit ${\cal M} = (M,\nabla)$ un $A$-module \`a
connexion.
\proclaim Definition 1.3.1. On dit que ${\cal M}$ est soluble dans l'extension
diff\'erentielle $(A', d')$ ({\rm ou simplement}: dans $A'\;$) si
$\; u^{\ast}{\cal M}= (M_{A'},\nabla_{A'})\;$ est trivial.
\par 
\proclaim Proposition 1.3.2. Supposons que $A'$ soit fid\`element plat sur $A$ et
sur $C'$, et que $\Omega^1$ soit fid\`ele et projectif de type
fini \`a droite sur $A$. Alors
\hfill \break $i)$ si ${\cal M}$ est soluble dans $A'$, sa volte est inversible.
\hfill \break $ii)$ si $(A',d')$ est simple (resp. simple par couches, avec $C'$
noeth\'erien), alors pour tout module \`a connexion $\cal M$ soluble dans $A'$ (resp. et
noeth\'erien), on a $A'\otimes_{C'} M_{A'}^{\nabla_{A'}} {\buildrel{\simeq}\over {\rightarrow}}
M_{A'}$. Tout sous-quotient de $\cal M$ est soluble dans $A'$ (resp. et noeth\'erien). La
cat\'egorie des modules \`a connexions solubles dans $A'$ (resp. et noeth\'eriens) est
ab\'elienne mono\"{\i}dale, et le foncteur ``solutions dans $A'$"  \break $\omega_{\scriptstyle
A'}:\;{\cal M}\mapsto M_{A'}^{\nabla_{A'}}$ est fid\`ele et exact.
\hfill \break $iii)$ Si $(A',d')$ est simple, tout $A$-module $M$ de type fini muni d'une
connexion soluble dans $A'$ est projectif. 
\hfill \break $iv)$ Si $(A',d')$ est simple par couches, alors pour tout id\'eal premier {\goth
p} de $C'$ et tout module \`a connexion soluble dans $A'$ et noeth\'erien $\cal M$, on a
$M_{A'}^{\nabla_{A'}}\otimes \kappa({\hbox{\goth p}})\cong   (M_{A'\otimes \kappa({\hbox{\goth
p}})})^{\nabla_{A'\otimes \kappa({\hbox{\smgoth p}})}}$.  
\hfill \break $v)$ supposons $C'$ est r\'egulier de dimension $\leq 1$ (par exemple un corps).
Soit $(M',\nabla')\hookrightarrow (M,\nabla)$ un monomorphisme d'objets solubles. Alors si $M$
est plat sur $A$ (resp. rigide, i.e. projectif de type fini d'apr\`es $i)$), il en est de m\^eme
de $M'$.  
\par {\it Preuve.} $i)$ r\'esulte par descente fid\`element plate de $A'$ \`a $A$ du fait que
la volte de $\nabla_{A'}$, qui est inversible, est l'unique prolongement $A'$-lin\'eaire \`a
droite de la volte de $\nabla$.    
\par \noindent $ii)$ suit de 1.1.1 et 1.2.3.
\par \noindent $iii)$: Si ${\cal M}$ est soluble dans $(A',d')$
simple, $C'$ est un corps et $A'\otimes_{C'} M_{A'}^{\nabla_{A'}}$; si $M$ est de type fini sur
$A$, $M_{A'}^{\nabla_{A'}}$ est donc un espace vectoriel de dimension finie sur $C'$. Donc
$M_{A'}$ est libre de type fini sur $A'$ et on conclut par descente fid\`element plate.
\par \noindent $iv)$: comme $(A',d')\otimes \kappa({\hbox{\goth p}})$ est simple par
hypoth\`ese, on a $A'\otimes_{C'}(M_{A'\otimes \kappa({\hbox{\goth p}})})^{\nabla_{A'\otimes
\kappa({\hbox{\smgoth p}})}} \break 
\hookrightarrow M_{A'\otimes \kappa({\hbox{\goth p}})}$. Par ailleurs, $A'\otimes_{C'} M_{A'}^{\nabla_{A'}}
{\buildrel{\simeq}\over {\rightarrow}} M_{A'}$. Par cons\'equent, le compos\'e 
\par \noindent $A'\otimes_{C'} M_{A'}^{\nabla_{A'}}\otimes \kappa({\hbox{\goth p}})\rightarrow A'\otimes_{C'}  (M_{A'\otimes
\kappa({\hbox{\goth p}})})^{\nabla_{A'\otimes \kappa({\hbox{\smgoth p}})}}\rightarrow(M_{A'\otimes
\kappa({\hbox{\goth p}})})$ est un isomorphisme, et il en est donc de m\^eme de
$M_{A'}^{\nabla_{A'}}\otimes \kappa({\hbox{\goth p}})\rightarrow(M_{A'\otimes
\kappa({\hbox{\goth p}})})^{\nabla_{A'\otimes \kappa({\hbox{\smgoth p}})}}$.   
\par \noindent $v)$: $M$ est $A$-plat (resp. projectif de
type fini) si et seulement si $M_{A'}$ est $A'$-plat. Pour $(M,\nabla)$ soluble dans $A'$, cela
revient \`a dire que $(M_{A'})^{\nabla_{A'}}$ est $C'$-plat (resp. projectif de
type fini). Or, si $C'$ est r\'egulier de dimension $\leq 1$ (i.e. produit fini d'anneaux de
Dedekind et de corps), un $C'$-module $N$ quelconque est plat (resp. projectif de
type fini) si et seulement si pour tout id\'eal premier {\goth p} de $C'$, le $C'_{\hbox {\goth
p}}-$module $N_{\hbox {\goth p}}$ est sans torsion (resp. et de type fini); ceci se propage donc
\`a tout $C'$-sous-module. D'o\`u le r\'esultat. (De m\^eme, on observe que les hypoth\`eses $C'$
r\'egulier de dimension $\leq 1$ et $\Omega'^1$ fid\`ele et projectif sur $A'$ entra\^{\i}nent que
$dA'$ est $C'$-plat, donc que $A'$ est fid\`element plat sur $C'$).
  
\medskip \noindent {\bf 1.3.3. Solubilit\'e et int\'egrabilit\'e.} Une
connexion {\it soluble dans une extension diff\'erentielle n'est pas n\'ecessairement
int\'egrable}, m\^eme dans la situation classique o\`u $\Omega^1$ est un bimodule commutatif (en
d\'epit de II.2.3.2 et de l'int\'egrabilit\'e des connexions triviales). 
\par \noindent Soit par exemple
$A= k[z_1,z_2]$ consid\'er\'e comme anneau diff\'erentiel g\'en\'eralis\'e via les d\'erivations
$d/dz_1,\;d/dz_2$; nous prenons donc $\Omega^1=\Omega^1_{A/k}= dz_1.A\oplus dz_2.A$,
le module de K\"ahler usuel. 
\par \noindent Consid\'erons le module \`a connexion
$(A,\nabla)$ correspondant au syst\`eme diff\'erentiel
\medskip \centerline{$\nabla(d/dz_1)y=0,\;\;\nabla(d/dz_2)y=z_1y.$}
\medskip \noindent Ce syst\`eme n'est pas
int\'egrable. 
\par \noindent Consid\'erons d'autre part l'anneau $A'= k[z_1,z_2,z_3,{1\over z_3}]$
comme anneau diff\'erentiel via les d\'erivations $d/dz_1,\;d/dz_2 + z_1d/dz_3$, et comme
extension diff\'erentielle de $(A,d)$ (noter que $(d/dz_2 + z_1d/dz_3)_{\vert A}=d/dz_2 $); le
crochet $[d/dz_1,d/dz_2 + z_1d/dz_3]$ n'est pas nul (c'est $d/dz_3$), mais sa restriction \`a $A$
l'est. L'anneau des constantes de $A'$ est $k$. 
\par \noindent Pour calculer $\Omega'^1$, remarquons que la base
duale de $(d/dz_1,\;d/dz_2 + z_1d/dz_3, d/dz_3)$ dans le module de K\"ahler $\Omega^1_{A'/k}$ est
$(dz_1,dz_2, dz_3-z_1dz_2)$; on a donc $\Omega'^1= \Omega^1_{A'/k}/A'.(dz_3-z_1dz_2)\;$ ($\cong
\Omega^1_{A/k}\otimes_{\scriptstyle A} A'$).  
\par \noindent Il est facile de voir que $(A,\nabla)$ est soluble dans $A'$:
$y=z_3$ est une base de solutions. On a $\Omega^2=\Omega^2_{A/k}= (dz_1\wedge dz_2).A$ (du moins
si $car \;k\neq 2$, cf. II.1.1). L'image de $dz_1\wedge dz_2$ dans $\Omega^2_{A'/k}$
s'\'ecrit $d(dz_3-z_1.dz_2)$, donc s'annule dans le quotient $\Omega'^2$ (en fait, on montre
facilement que $\Omega'^2=0$). 
\par \noindent On prendra garde \`a la confusion que peut cr\'eer le mot d'int\'egrabilit\'e, un
syst\`eme non-int\'egrable pouvant \^etre ``int\'egr\'e" symboliquement. En fait, nous verrons
qu'on peut d\'evelopper la th\'eorie de Picard-Vessiot sans hypoth\`ese sur la courbure. 
\par Toutefois, on a un carr\'e commutatif 
\medskip\centerline{ $\matrix{M & \rightarrow &\Omega^2\otimes_A M \cr
         \downarrow &&\downarrow  \cr M'&\rightarrow &\Omega'^2 \otimes_{A'} M'\cong \Omega'^2
\otimes_A M}$}
\medskip\noindent o\`u les fl\`eches horizontales sont les courbures de $\nabla$ et $\nabla'$ respectivement. 
Si la fl\`eche canonique $\Omega^2\rightarrow \Omega'^2 $ est injective, et si $M$ est plat sur $A$,
la solubilit\'e de $\nabla$ dans $A'$ implique l'int\'egrabilit\'e de  $\nabla$.

\medskip \noindent {\bf 1.3.4. Remarque.} Pour travailler sans aucune hypoth\`ese de simplicit\'e
diff\'erentielle (par exemple en caract\'eristique mixte), on est amen\'e \`a ne consid\'erer que
les morphismes ``admissibles" de connexions solubles (i.e. admissibles au sens de 1.1.2 apr\`es
extension \`a $A'$). Si $A'$ est fid\`element plat sur $A$ et sur $C'$, on obtient une
cat\'egorie ab\'elienne mono\"{\i}dale. Nous ne d\'evelopperons pas ce point de vue ici. 

\medskip {\bf 1.4. Exemples fondamentaux (cas int\'egrables).}
\medskip \noindent {\bf 1.4.1.} Soit $X$ une vari\'et\'e alg\'ebrique lisse sur un corps $k$ de
caract\'eristique nulle ou une vari\'et\'e analytique lisse sur
$k={\bf R}$ ou $\bf C$ (voire une vari\'et\'e analytique $p$-adique rigide lisse). On note $m$ sa
dimension. Soient $x$ un $k$-point de $X$, ${\cal O}_{X,x}$ l'anneau local de $X$ au point (d\'efini par)
$x$, et $z_1,\ldots,z_m$ des coordonn\'ees locales autour de $x$. Alors le compl\'et\'e ${\hat A}={\hat
{\cal O}_{X,x}}\cong k[[z_1,\ldots,z_m]]$ est une extension fid\`element plate de $A$. On consid\`ere $A$
(resp. ${\hat A}$) comme anneau diff\'erentiel g\'en\'eralis\'e, gr\^ace \`a la d\'erivation
$d$ vers le module de diff\'erentielles usuelles $\Omega^1 = \Omega^1_{X/k,x}$ (resp. $\Omega^1
\otimes_{\scriptstyle A} {\hat A}$). On obtient ainsi une extension d'anneaux diff\'erentiels simples, de corps
de constantes $k$. D'apr\`es le th\'eor\`eme de Frobenius formel, tout ${\cal O}_{X,x}$-module
de pr\'esentation finie $M$ muni d'une connexion {\it int\'egrable} $\nabla$ est soluble dans
${\hat {\cal O}_{X,x}}$. En particulier, il est libre de type fini sur l'anneau local ${\cal
O}_{X,x}$. On a ${\hat M}\cong {\hat {\cal O}_{X,x}}\otimes M$, et ${\hat M}^{\nabla}$
s'identifie \`a la fibre de $M$ en $x$.   
\par \noindent Dans la situation alg\'ebrique, on peut plus g\'en\'eralement remplacer $k$ par
une $\bf Q$-alg\`ebre commutative int\`egre noeth\'erienne quelconque, pourvu que les fibres de
$X$ soient g\'eom\'etriquement connexes. La m\^eme construction fournit une extension fid\`element
plate d'anneaux diff\'erentiels simples par couches. 
\par \noindent On en d\'eduit que l'objet $(M,\nabla)$ est rigide (cf. II.4.2.1) si et seulement si
${\hat M}^{\nabla}$ est plat sur $k$. En particulier, lorsque $k$ est l'anneau de fonctions d'une
courbe affine lisse, un ${\cal O}_{X,x}$-module de type fini $M$ muni d'une connexion int\'egrable
$\nabla$ (relativement \`a $k$) s'interpr\`ete comme famille \`a un param\`etre de syst\`emes
diff\'erentiels lin\'eaires int\'egrables. Il est rigide si et seulement s'il est sans torsion sur
$k$. 
\par \noindent On en d\'eduit aussi que le foncteur ``tige en $x$" 
\medskip \centerline{\it $\Bigl\{ \matrix{ {\cal O}_X \hbox{-modules coh\'erents} \cr \hbox{\`a
connexion int\'egrable} } \Bigr\} \longrightarrow \Bigl\{ \matrix{{\cal O}_{X,x}\hbox{-module
de type fini} \cr \hbox{\`a connexion
int\'egrable} } \Bigr\} $}
\medskip \noindent est exact et pleinement fid\`ele, et que le foncteur ``fibre en $x$" 
\medskip \centerline{\it $\Bigl\{ \matrix{ {\cal O}_X \hbox{-modules coh\'erents} \cr \hbox{\`a
connexion int\'egrable}}\Bigr\} \longrightarrow
\lbrace  k\hbox{-modules de type fini}\rbrace$}
\medskip \noindent est exact et fid\`ele.
\medskip \noindent {\bf 1.4.2.} Dans la situation I.2.2, $\bar{\cal O}_x \hookrightarrow {\cal
O}^{pd}_x$ donne lieu \`a une extension fid\`element plate d'anneaux diff\'erentiels simples, de
corps des constantes le corps $k$ de caract\'eristique $p$. On a vu que toute connexion
int\'egrable $\nabla$ sur $\bar X_x$ est soluble dans l'alg\`ebre \`a puissances divis\'ees
compl\'et\'ee ${\cal O}^{pd}_x$.
\par \noindent L\`a encore, on a une variante ``relative" (remplacer $k$ par
une ${\bf F}_p$-alg\`ebre commutative int\`egre noeth\'erienne quelconque). 

\medskip \noindent {\bf 1.4.3.} {\it Equations aux diff\'erences.} On
part d'un anneau $k$ quotient de ${\bf C}[[h_1,\ldots,h_m]]$. On note
 ${\hbar}_i$ l'image de $h_i$ dans $k$. On fixe un $k$-point $a=(a_1,\ldots,a_m)$ de l'espace
affine ${\bf A}^m_{\scriptstyle k}$, et on note $\bar a$ le point ferm\'e correspondant. On consid\`ere
l'anneau local $A={\cal O}_{{\bf A}^m,{\bar a}}\;$ muni des endomorphismes $\sigma_i: z_j
\mapsto z_j\; {\hbox{si}}\; i\neq j, \;\;z_i
\mapsto z_i+{\hbar}_i$, ainsi que son compl\'et\'e ${\hat A}\cong k[[z_1-a_1,\ldots,z_m-a_m]]$.
Noter que ${\hat A}$ est fid\`element plat sur $A$ et sur $k$. On fait de $A$
et $\hat A$ des anneaux diff\'erentiels comme en I.4.4.2: $\Omega^1 = \oplus 
\;\Omega^1_{\sigma_i}\cong \oplus \; dz_i.A, \;,\;{\hat\Omega}^1=\Omega^1\otimes {\hat A}\cong
\oplus \; dz_i.{\hat A}, $ et $d$ est donn\'ee par $ df= \sum dz_i.\delta_i(f)$. Les anneaux de
constantes co\"{\i}ncident avec $k$. 
\par V\'erifions que l'anneau diff\'erentiel $\hat A$ est {\it simple par
couches}. Pour cela, on remarque que tout \'el\'ement de $\hat A$ s'\'ecrit, de mani\`ere unique,
comme s\'erie 
\medskip \centerline{$f=\Sigma_{(n_1,\ldots,n_m)}\;\alpha_{n_1,\ldots,n_m}\Pi_i\Pi_{j=0}^{j=n_i}
(z_i-a_i+j{\hbar}_i)$}
\medskip \noindent \`a coefficients $\;\alpha_{n_1,\ldots,n_m}
\in k,\;\;$ et que  
\medskip \centerline{$\delta_i[\Pi_{j=0}^{j=n_i}
(z_i-a_i+j{\hbar}_i)]=n_i\Pi_{j=0}^{j=n_i-1}
(z_i-a_i+j{\hbar}_i)\;$ pour $n>0$.}
\par \noindent (De cette derni\`ere \'egalit\'e, on d\'eduit par ailleurs le calcul des
coefficients: 
\medskip \centerline{$\alpha_{n_1,\ldots,n_m}\;
 = [(\Pi_i{1\over 
n_i!}\;\delta_i^{n_i})(f)]_{\vert a}\;$.)}
\par \noindent On conclut en consid\'erant, dans un id\'eal diff\'erentiel $I$ non nul, un
\'el\'ement $f$ de plus bas degr\'e total en $z_1-a_1,\ldots,z_m-a_m$; comme $n_i$ est inversible
dans $k$, on voit que le premier coefficient de $f$ est dans $I$.
\medskip Consid\'erons un syst\`eme
lin\'eaire aux diff\'erences
\medskip \centerline{$\sigma_i(Y)={\cal A}_{(i)}.Y$}
\par \noindent o\`u ${\cal A}_{(i)} \in GL_{\mu}(A)$. On suppose de plus: 
\par $i)$ que le syst\`eme est int\'egrable: ${{\cal A}_{(j)}}^{\sigma_i}{\cal
A}_{(i)}={{\cal A}_{(i)}}^{\sigma_j}{\cal A}_{(j)}$ (II.2.3.3);
\par $ii)$ qu'une fois r\'ecrit sous la forme
\medskip \centerline {$\delta_i(Y)={\cal
G}_{(i)}.Y$}
\par \noindent les matrices ${\cal G}_{(i)}$ sont encore \`a coefficients dans $A$;
sous cette hypoth\`ese, on a alors {\it confluence} vers un syst\`eme
diff\'erentiel sans singularit\'e en $\bar a$ lorsque $\hbar_i \rightarrow 0$. 
\medskip \noindent On peut r\'ecrire ce syst\`eme sous forme de $A$-module libre $M= \oplus
A.m_{\scriptstyle k}$ \`a connexion:  
\medskip \centerline {$\nabla(m_{\scriptstyle k})=\oplus_{i,l}\; {dz_i\over {\hbar}_i}\otimes ({\cal
A}^{-1}_{(i)lk}-\delta_{lk})\;m_l\;\;\;$ (loc.
cit.).} 
\medskip \noindent {\it Une telle connexion est toujours soluble dans $\hat A$} (c'est
l'analogue pour les \'equations aux diff\'erences du th\'eor\`eme de Frobenius formel). Pour le
voir, il suffit de trouver une solution $Y \in GL_{\mu}({\hat A})$ du syst\`eme $\delta_i(Y)={\cal
G}_{(i)}.Y, \;\;i=1,\ldots m $. En it\'erant l'action des $\delta_i$, on d\'efinit formellement
pour tout multi-indice ${\u n}=(n_1,\ldots,n_m)$ une matrice ${\cal G}_{(\u n)} \in M_{\mu}(A)$
d\'efinie par   
\medskip \centerline {$(\Pi_i{1\over 
n_i!}\;\delta_i^{n_i})(Y)={\cal G}_{(\u n)}.Y$} 
\par \noindent (l'ordre dans le produit $\Pi_i$ est sans
cons\'equence en vertu de l'hypoth\`ese d'int\'egrabilit\'e).
\par \noindent Il d\'ecoule alors des calculs pr\'ec\'edents sur les s\'eries en 
$\Pi_{j=0}^{j=n_i-1}(z_i-a_i+j{\hbar}_i)$ qu'il existe une unique solution $Y \in GL_{\mu}({\hat A})$
de 
\medskip \centerline {$\delta_i(Y)={\cal
G}_{(i)}.Y, \;(i=1,\ldots m), \;Y(a) = id $.}
\par \noindent Elle est donn\'ee par 
\medskip \centerline {$Y=\sum_{{\u n}=(n_1,\ldots,n_m)}\;{\cal G}_{(\u n)}(a).\Pi_i
[\Pi_{j=0}^{j=n_i} (z_i-a_i+j{{\hbar}}_i)]$.}
\par \noindent On a ${\hat M}\cong {\hat A}\otimes M$, et ${\hat M}^{\nabla}$
s'identifie au $k$-module fibre de $M$ en $a$. 
\par \noindent Il r\'esulte de ce qui pr\'ec\`ede et de 1.3.2 que tout sous-quotient de
$(M,\nabla)$ est soluble dans $\hat A$, et que pour tout sous-objet $(N,\nabla)$ de
$(M,\nabla)$, $N$ est libre de type fini sur $A$.      

\medskip \noindent {\bf 1.4.4.} {\it Equations aux $q$-diff\'erences.} La situation est analogue
\`a la pr\'ec\'edente, avec les changements suivants:
\par \noindent $i)$ changer $h_i$ en $q_i-1$,  
\par \noindent $ii)$ $\sigma_i: z_j
\mapsto z_j\; {\hbox{si}}\;\; i\neq j, \;\;z_i
\mapsto q_iz_i$,
\par \noindent $iii)$ tout \'el\'ement de $\hat A$ s'\'ecrit comme s\'erie 
\medskip \centerline{$f=\Sigma_{(n_1,\ldots,n_m)}\;\alpha_{n_1,\ldots,n_m}\Pi_i\Pi_{j=0}^{j=n_i}
(z_i-q_i^j a_i)$}
\medskip \noindent \`a coefficients $\alpha_{n_1,\ldots,n_m}= [(\Pi_i{1\over 
[n_i]!_{q_i}}\;\delta_i^{n_i})(f)]_{\vert a}\;\in k$, et on a  
\medskip \centerline{$\delta_i[\Pi_{j=0}^{j=n_i}
(z_i-q_i^j a_i)]=(n_i)_{q_i}\Pi_{j=0}^{j=n_i-1}
(z_i-q_i^j a_i)\;$ pour $n>0$,}
\par \noindent avec $(n_i)_{q_i}=1+q_i+\ldots
+q_i^{n_i-1},\;[n_i]!_{q_i}=\Pi_{j=1}^{j=n_i}\;(j)_{q_i}\;$,
\par \noindent $iv)$ $\nabla(m_{\scriptstyle k})=\oplus_{i,l}\; dz_i\otimes {{\cal
A}^{-1}_{(i)lk}-\delta_{lk}\over (q_i-1)z_i}\;m_l$,
\par \noindent $v)$ $(\Pi_i {1\over 
{[n_i]!_{q_i}}}\;\delta_i^{n_i})(Y)={\cal G}_{(\u n)}.Y, \;\;Y=\sum_{{\u
n}=(n_1,\ldots,n_m)}\;{\cal G}_{(\u n)}(a).\Pi_i [\Pi_{j=0}^{j=n_i} (z_i-q_i^ja_i)]$.
\medskip \noindent L'exemple typique est celui des \'equations $q$-hyperg\'eom\'etriques confluant
vers une \'equation diff\'erentielle hyperg\'eom\'etrique (on prend ${\bar a}\neq 0,1$).  
\medskip Dans tous ces exemples, on peut v\'erifier que l'int\'egrabilit\'e est
pr\'eserv\'ee par produit tensoriel et passage au dual. 
\medskip
\bigskip {\bf 2. Groupes de Galois diff\'erentiels.} 
\medskip {\it Dans toute la suite, on consid\`ere un anneau diff\'erentiel $(A,d)$
(avec $A$ commutatif, bien s\^ur), d'anneau de constantes $C=k$} noeth\'erien{\it ,
et on suppose que $\Omega^1=dA.A$ est} fid\`ele et projectif de type fini {\it comme $A$-module
\`a droite (ou, ce qui revient au m\^eme, fid\`element plat et de pr\'esentation finie).} 

\medskip {\bf 2.1. Definition et exemples.} Soit ${\cal M}=(M,\nabla)$ un $A$-module \`a 
connexion rigide (i.e. $M$ est projectif de type fini et la volte de $\nabla$ est un
isomorphisme, cf. II.4.2.1). Pour tout couple d'entiers naturels $(i,j)$, on pose $T^{i,j}({\cal M})={\cal
M}^{\otimes i}\otimes {\check{\cal M}}^{\otimes j}$. 
\par \noindent On note ${\scriptstyle <}{\cal M}{\scriptstyle >}^{\scriptstyle\otimes}$ la sous-cat\'egorie strictement pleine
($k$-lin\'eaire, ab\'elienne, mono\"{\i}dale sym\'etrique) de la cat\'egorie des $A$-modules \`a
connexion form\'ee des sous-quotients des sommes finies $\oplus \; T^{i,j}({\cal M})$.  

\proclaim Theor\`eme 2.1.1. Supposons donn\'ee une sous-cat\'egorie $\;{\cal C}\;$ pleine
mono\"{\i}dale ab\'elienne de $\;{\scriptstyle <}{\cal M}{\scriptstyle >}^{\scriptstyle\otimes}\;$ contenant $\;{\cal M}\;$. On suppose 
que dans $\;{\cal C} \;$, tout objet est quotient d'un objet rigide, et que le dual d'un objet
rigide est encore un objet de $\;\cal C$. Supposons donn\'e en outre un foncteur mono\"{\i}dal,
$k$-lin\'eaire, exact et fid\`ele $\;\;\omega:\; {\cal C} \;\longrightarrow $ ($k$-modules de
type fini)$\;\;$ qui envoie objets rigides sur $k$-modules projectifs de type fini. Alors
\hfill\break $i)$ il existe un $k$-groupe affine (fid\`element) plat $\;G({\cal C},\omega)$, muni
d'un monomorphisme naturel
$\;\iota:\;G({\cal C},\omega)\rightarrow GL(\omega({\cal M}))$, tel que $\;\omega\;$ induise une
$\otimes$-\'equivalence de cat\'egories entre $\;{\cal C} \;$ et la cat\'egorie des
repr\'esentations de type fini sur $\;k \;$ de $\;G({\cal C},\omega)$. 
\hfill \break $ii)$ Si $\;k\;$
est un corps, le monomorphisme $\;\iota\;$ est une immersion ferm\'ee et son image s'identifie au 
sous-groupe ferm\'e de
$\;GL(\omega({\cal M}))\;$ qui stabilise les $\;\omega({\cal N})$, pour tout sous-objet
$\;\cal N \;$ d'une somme finie quelconque $\;\oplus \; T^{i,j}({\cal M})$. 
\par 
\medskip \noindent {\bf 2.1.2. Notation et d\'efinition.} Si $\;{\cal C}=\;{\scriptstyle <}{\cal
M}{\scriptstyle >}^{\scriptstyle\otimes}\;$, on \'ecrira $\;Gal({\cal M},\omega)\;$ au lieu de $\;G({\cal
C},\omega)\;$. C'est le {\it groupe de Galois diff\'erentiel de $\cal M$ point\'e en $\omega$}. 

\medskip {\it Preuve de 2.1.1}. La premi\`ere assertion de 2.1.1. est une application
directe du lemme tannakien 8.1.2 de [An96]; $G({\cal C},\omega)={Aut}^{\scriptstyle\otimes}\omega$
repr\'esente le foncteur $\u{Aut}^{\scriptstyle\otimes}\omega$ sur les $k$-alg\`ebres commutatives unitaires
(dans loc. cit., on demande que $\omega$ commute \`a la dualit\'e sur les objets rigides, mais
c'est automatique, cf. [Bru94]2.2). 
\par \noindent Si $k$ est un corps, le foncteur $\u{Stab} \lbrace
\omega({\cal N}) \rbrace$ qui \`a toute $k$-alg\`ebre commutative unitaire $k'$ associe le
stabilisateur dans $GL(\omega({\cal M}))(k')$ des $\omega({\cal N})\otimes k'$ (pour tout 
sous-objet $\cal N$ dans $\cal C$ d'une somme finie quelconque $\;\oplus \; T^{i,j}({\cal M})\;$)
est repr\'esentable par un sous-sch\'ema ferm\'e $Stab \lbrace \omega({\cal N}) \rbrace$ de
$GL(\omega({\cal M}))$ ([DG70]II.1.3.6); c'est le stabilisateur dont il s'agit dans $ii)$. On a
clairement $\u{Aut}^{\scriptstyle\otimes}\omega \; \subset \u{Stab} \lbrace \omega({\cal N}) \rbrace$, ce qui
signifie que le monomorphisme $\iota$ se
factorise \`a travers $Stab \lbrace \omega({\cal N}) \rbrace$. R\'eciproquement, toute
repr\'esentation de type fini sur $k$ de $G({\cal C},\omega)$ (c'est-\`a-dire l'image par
$\omega$ de tout objet de ${\cal C}$) est quotient de l'image par $\omega$ d'un
sous-objet de $\oplus \; T^{i,j}({\cal M})$; elle d\'efinit donc
une repr\'esentation de type fini sur $k$ de $Stab \lbrace \omega({\cal N}) \rbrace$. D'apr\`es
[Saa72]II.3.3 (apr\`es passage aux cat\'egories des $ind$-objets, cf. loc. cit 2.2.3.1), cette
correspondance fonctorielle mono\"{\i}dale entre cat\'egories de repr\'esentations provient d'un
homomorphisme $Stab \lbrace \omega({\cal N}) \rbrace\rightarrow G({\cal C},\omega)$, qui est
inverse de 
\break $G({\cal C},\omega) \hookrightarrow Stab \lbrace \omega({\cal N})
\rbrace$.    
\medskip \noindent {\bf 2.1.3. Exemples g\'en\'eraux.} $i)$ Soit $(A', d')$ une extension
diff\'erentielle de $(A,d)$, avec $C'=k$ r\'egulier de dimension $\leq 1$, 
et $A'$ fid\`element plat sur $A$.  
On suppose de plus $(A',d')$ simple ou simple par couches, et ${\cal M}$ soluble dans $A'$. On
prend pour $\omega$ le foncteur $\omega_{\scriptstyle
A'}:\;{\cal N}\mapsto ({\cal N}_{A'})^{\nabla_{A'}}$ sur ${\cal
C}=\;{\scriptstyle <}{\cal M}{\scriptstyle >}^{\scriptstyle\otimes}$. 
\par \noindent Il r\'esulte de la proposition 1.3.2 que les hypoth\`eses de 2.1.1 sont
satisfaites.  
\par \noindent Ceci s'applique aux exemples de 1.4. En 1.4.3 (resp. 1.4.4), on prendra pour $k$
un quotient r\'egulier de dimension un de ${\bf C}[[h_1,\ldots,h_m]]$ (resp. ${\bf
C}[[q_1-1,\ldots,q_m-1]]$), par exemple celui d\'efini par $h_1=\ldots =h_m$ (resp. $q_1=\ldots
=q_m$).   
\medskip $ii)$ Si $(A,d)$ v\'erifie les hypoth\`eses de II.5.3.2, ${\scriptstyle <}{\cal
M}{\scriptstyle >}^{\scriptstyle\otimes}$ est tannakienne. Si elle admet un foncteur fibre $\omega$ sur
$k$, $\;Gal({\cal M},\omega)\;$ est le groupe tannakien associ\'e.  
\medskip $iii)$ Les groupes de Galois diff\'erentiels de la th\'eorie de Picard-Vessiot
(I.3), les groupes de Galois aux diff\'erences \'etudi\'es dans [vdP95], et les groupes ``mixtes"
de [Bi62] sont des exemples de groupes $Gal({\cal M},\omega)$; ce sont du reste des cas particuliers de
$i)$ ou $ii)$, o\`u $A$ est un corps.   
\medskip $iv)$ Les enveloppes alg\'ebriques des groupes d'holonomie (I.1.3) sont aussi
des incarnations de groupes $G({\cal C},\omega)$. En effet, soit $P\rightarrow X$ un fibr\'e
principal \`a droite sous $GL_n$, muni d'une connexion $\aleph$, et soit $(E,\nabla)$ le fibr\'e
vectoriel \`a connexion sur $X$ associ\'e. On peut prendre pour
$\cal C$ la cat\'egorie des fibr\'es \`a connexion de la forme $V/W$, o\`u $W$, resp. $V$, sont
des champs de $p$-plans, resp.$q$-plans, stables sous connexion dans les
$\oplus \; T^{i,j}(E)$ (pour $p\leq q$ arbitraires). Le foncteur fibre en $x$ induit une
\'equivalence $\omega$ entre $\cal C$ et la cat\'egorie des rep\'esentations r\'eelles de
dimension finie de l'enveloppe alg\'ebrique $G({\cal C},\omega)$ de $Hol_x(\nabla)$: en effet,
d'une part toute rep\'esentation r\'eelle de dimension finie de $G({\cal C},\omega)$ est
sous-quotient d'un $\oplus \; T^{i,j}(E_x)$, d'autre part $P$ provient d'un fibr\'e principal sous
$Hol_x(\nabla)$, donc aussi d'un fibr\'e principal sous $G({\cal C},\omega)$; la
construction du fibr\'e vectoriel associ\'e \`a \break toute rep\'esentation r\'eelle de dimension
finie de $G({\cal C},\omega)$ fournit un inverse de $\omega$. 
\par \noindent Si $X$ est riemannienne simplement connexe et si $E$ est le fibr\'e
tangent muni de la connexion de Levi-Civita, le groupe compact connexe $Hol_x(\nabla)$
s'identifie au groupe des points r\'eels de $G({\cal C},\omega)$.
\par \noindent A l'autre extr\^eme, si la connexion $\nabla$ sur $E$ est int\'egrable,
$G({\cal C},\omega)$ est l'enveloppe alg\'ebrique du groupe de monodromie, ce qui exprime
l'\'equivalence bien connue entre fibr\'es \`a connexion int\'egrable et syst\`emes locaux.      

\proclaim Lemme 2.1.4. Supposons que les hypoth\`eses de 2.1.1 sont satisfaites avec $\;{\cal C}=
\;{\scriptstyle <}{\cal M}{\scriptstyle >}^{\scriptstyle\otimes}\;$. Alors $\;Gal({\cal M},\omega)\;$ est le groupe trivial si et
seulement si $\cal M$ est un module \`a connexion trivial. 
\par \noindent En effet, si $Gal({\cal M},\omega)$ est trivial, on a
$M^{{\scriptstyle \nabla}}\cong Mor((A,d),{\cal M})\cong \omega({\cal M})$. En particulier,
$M^{{\scriptstyle \nabla}}$ est projectif de type fini sur $k$, donc facteur d'un $k^n$; ainsi ${\cal
M}'=(A\otimes M^{{\scriptstyle \nabla}}, d\otimes id_{M^{\scriptscriptstyle{\nabla}}})$ est facteur
direct de $(A,d)^n$, donc est un objet de ${\scriptstyle <}{\cal M}{\scriptstyle >}^{\scriptstyle\otimes}$. Le morphisme naturel ${\cal M}'
\rightarrow {\cal M}$ est un isomorphisme, puisqu'il en est ainsi de son image par $\omega$.
R\'eciproquement, si $\cal M$ est un module \`a connexion trivial, il est quotient d'un
$(A,d)^n$. Ainsi ${\scriptstyle <}{\cal M}{\scriptstyle >}^{\scriptstyle\otimes}$ n'est autre que la cat\'egorie mono\"{\i}dale des
sous-quotients des $(A,d)^n$, \'equivalente via $\omega$ \`a la cat\'egorie mono\"{\i}dale des
$k$-modules de type fini. Donc $\;Gal({\cal M},\omega)\;$ est le groupe trivial. 

\medskip \noindent {\bf 2.1.5. Question de finitude.} Consid\'erons l'exemple $k= {\bf
C}[[h]], \;A=$ le localis\'e de $k[z]$ en $z=0,\; d: A \rightarrow Adz $ la diff\'erentielle de
K\"ahler standard. Si nous munissons $A$ de la connexion $\nabla(1)=dz$, et prenons pour $\omega$
le foncteur fibre en $z=0$, le groupe de Galois diff\'erentiel est simplement ${\bf G}_m = Spec\;
k[z,{1\over z}]\cong Spec\; k[x,y]/(xy+x+y)$ (avec $x=z-1, y={1\over z}-1$). 
\par \noindent L'objet $\cal M$ obtenu en munissant $A$ de la connexion $\nabla(1)=h.dz$ est plus
d\'elicat. Remarquons que pour tout $N{\scriptstyle >}0, \;(k/h^{\scriptscriptstyle N} k) \otimes
\cal M\;$ est trivial: une base de solutions est donn\'ee par $\sum_{\scriptstyle
n=0}^{\scriptstyle n=N-1} {(-hz)^n \over n!}\;.$ Le groupe de Galois diff\'erentiel est 
\medskip \centerline{$Gal({\cal M},\omega) = \limproj \;Spec\; k[x,y]/(h^nxy+x+y)$}
\medskip \noindent (les morphismes de transition \'etant donn\'es par la multiplication de $x$ et
$y$ par $h$). Ce sch\'ema en groupes plat {\it n'est pas de type fini} sur $k$, et le
monomorphisme $\iota: \;Gal({\cal M},\omega) \rightarrow {\bf G}_m$ n'est pas une immersion
ferm\'ee. Si $k'$ est une ${\bf C}((h))$-alg\`ebre, on a $Gal({\cal
M},\omega)(k')=k'^{\ast}$, tandis que si $\;k'=\bf C\;$ est le quotient de $\;k,\;$ on a $\;Gal({\cal
M},\omega)(k')=\lbrace 1 \rbrace$.
\par L'apparition de sch\'emas en groupes non de type fini est un trait caract\'eristique de
la th\'eorie tannakienne sur des bases qui ne sont pas des corps.   

\bigskip {\bf 2.2. Torseur des solutions et groupe de Galois diff\'erentiel
intrins\`eque.}
\medskip \noindent {\bf 2.2.1.} Nous revenons \`a la situation de 2.1.2. en faisant
l'hypoth\`ese suppl\'ementaire que {\it tout objet de 
${\scriptstyle <}{\cal M}{\scriptstyle >}^{\scriptstyle\otimes}$ est rigide} et que
l'anneau $k$ des constantes est un corps.
\par \noindent Cette hypoth\`ese est v\'erifi\'ee en particulier si ${\cal M}$ est soluble dans
une extension diff\'erentielle {\it simple} $(A', d')$ de $(A,d)$, avec $A'$ fid\`element plat sur
$A$.   
\par \noindent Sous cette hypoth\`ese, on dispose du foncteur ``oubli de la connexion"
\medskip \centerline{$oubli:\;\;{\scriptstyle <}{\cal M}{\scriptstyle >}^{\scriptstyle\otimes} \; \rightarrow \; \lbrace A-{\hbox{\it
modules projectifs de type fini}}\rbrace$} 
\medskip \noindent Le foncteur
$\u{Isom}^{\scriptstyle\otimes}(\omega \otimes 1_{\scriptstyle A},\;\hbox{\smit{oubli}})$ sur les
$A$-alg\`ebres commutatives unitaires est alors repr\'esentable par un $A$-sch\'ema affine
$\;\Sigma({\cal M},\omega)$. C'est un torseur sous 
\break $\;Gal({\cal M},\omega)\otimes_{\scriptstyle k}
A$ (agissant \`a droite par composition); en particulier, il est fid\`element plat sur $A$
[Saa72]II.4.2.
\par \noindent Dans la situation \'evoqu\'ee ci-dessus o\`u est ${\cal M}$ est soluble dans $(A',
d')$ simple, et o\`u $\omega$ est le foncteur $\omega_{\scriptstyle
A'}:\;{\cal N}\mapsto N_{A'}^{\nabla_{A'}}$, ce torseur
est appel\'e {\it torseur des solutions} de ${\cal M}$ dans $A'$; en effet, on verra plus loin que
lorsque $M$ est libre, l'alg\`ebre des fonctions sur $\Sigma({\cal M},\omega)$ est engendr\'ee
par les coefficients d'une ``matrice fondamentale de solutions" (4.2).   
\medskip \noindent {\bf 2.2.2.} De m\^eme, le foncteur 
$\u{Aut}^{\scriptstyle\otimes}(\hbox{\smit{oubli}})$
sur les $A$-alg\`ebres commutatives unitaires est repr\'esentable par un $A$-sch\'ema en groupes
affine plat $\;Gal({\cal M},\hbox{\smit{oubli}})$, appel\'e {\it groupe de Galois
diff\'erentiel intrins\`eque} (ou encore {\it g\'en\'erique} si $A$ est un corps) de ${\cal
M}$.  
\par \noindent Ce sch\'ema en groupes agit \`a gauche
par composition sur $\Sigma({\cal M},\omega)$, qui est en fait un
bitorseur sous $\;Gal({\cal
M},\hbox{\smit{oubli}})$ et $\;Gal({\cal M},\omega)\otimes_{\scriptstyle k} A\;$.  
\par \noindent Le sch\'ema en groupes $\;Gal({\cal M},\hbox{\smit{oubli}})$ n'est autre que
le sous-groupe ferm\'e de $\;GL( M)\;$ qui stabilise les sous-objets ${\cal N}=(N,\nabla)$ des
sommes finies $\;\oplus \; T^{i,j}({\cal M})$ (noter que les 
sous-$A$-modules sous-jacents sont par hypoth\`ese facteurs directs, de sorte que le
stabilisateur $Stab \lbrace {\cal N} \rbrace$ en question est bien d\'efini: il repr\'esente le
foncteur associant \`a toute $A$-alg\`ebre commutative unitaire $R$ le stabilisateur dans $\;GL(
M_R)\;$ des $N_R$). Ceci peut se voir en utilisant la propri\'et\'e correspondante pour
$\;Gal({\cal M},\omega)\otimes_{\scriptstyle k} A$ (cf. 2.1.1.$ii)$), et la formule 
$\;Gal({\cal M},\hbox{\smit{oubli}})=Aut_{(Gal({\cal
M},\omega)\otimes_{\scriptstyle k} A)}\Sigma({\cal M},\omega)$. 
\par \noindent Il faut toutefois prendre garde \`a la r\'eciproque: il n'est pas vrai en
g\'en\'eral que tout sous-$A$-module de $M$ stable sous $\;Gal({\cal
M},\hbox{\smit{oubli}}) = Stab \lbrace {\cal N} \rbrace$ soit stable sous la connexion.
\medskip \noindent {\bf 2.2.3.} Il r\'esulte de ce qui pr\'ec\`ede (en prenant
$R=A[\epsilon]/(\epsilon^2)$) que l'alg\`ebre de Lie de
$\;Gal({\cal M},\hbox{\smit{oubli}})$ est la sous-$A$-alg\`ebre de Lie 
$\;LieStab \lbrace
{\cal N}\rbrace$ de {\goth gl}$(M)$ qui stabilise les sous-objets des sommes finies
$\;\oplus \; T^{i,j}({\cal M})$, pour l'action de Lie naturelle sur les puissances
tensorielles mixtes. 
\par \noindent Remarquons que {\goth gl}$(M)$ s'identifie au $A$-module sous-jacent au module \`a
connexion  $T^{1,1}({\cal M})= {\cal M}\otimes {\check{\cal M}}$. 
\proclaim Proposition 2.2.4. $LieGal({\cal M},\hbox{\smit{oubli}})$ est (sous-jacent \`a) un
sous-module \`a connexion de $T^{1,1}({\cal M})$. Tout sous-module \`a connexion de
$LieGal({\cal M},\hbox{\smit{oubli}})$ est un id\'eal de Lie. \par  {\it Preuve.} Pour
prouver la premi\`ere assertion, il s'agit de faire voir que pour tout
$\ell \in LieStab \lbrace
{\cal N}\rbrace$, et tout sous-objet ${\cal N}=(N,\nabla)$ d'une 
somme finie $\;{\cal M}':= \oplus \; T^{i,j}({\cal M}), \;(\nabla \ell)N \subset \Omega^1\otimes
N$.  Pour tout $n\in N, \;(\nabla \ell)n$ se calcule via l'action de Lie naturelle   
\medskip \centerline{{\goth L}: $T^{1,1}({\cal M}) \otimes {\cal M}' \rightarrow
{\cal M}'$}
\medskip \noindent qui est un morphisme de modules \`a connexion:
\medskip \centerline{$\nabla(\ell)n=(1_{ \Omega^1}\otimes\hbox{\goth L} )[\nabla(\ell)\otimes
n] $.}
\medskip \noindent On a {\goth L}$(\ell\otimes n)=\ell.n\in N$, d'o\`u $\nabla (\ell.n)\in
\Omega^1\otimes N$. Calculons ce dernier en faisant intervenir la contrainte de commutativit\'e
$c=c(T^{1,1}({\cal M}),\;{\cal M}')$ (l'\'echange des facteurs, qui est un morphisme de modules
\`a connexion): 
\medskip \centerline{$\nabla (\ell.n)=\nabla (\hbox{\goth L}(c(n\otimes \ell )))=(1_{
\Omega^1}\otimes(\hbox{\goth L}\circ c))\nabla (n\otimes \ell) $}
\medskip \centerline{$=
(1_{ \Omega^1}\otimes(\hbox{\goth L}\circ c))[\nabla(n)\otimes \ell + (\phi(\nabla_{{\cal
M}'})\otimes 1_{T^{1,1}({\cal M})})(n\otimes
\nabla(\ell))]$}
\medskip \centerline{$=
(1_{ \Omega^1}\otimes(\hbox{\goth L}\circ c))[\nabla(n)\otimes \ell + (\phi(\nabla_{\cal
N})\otimes 1_{T^{1,1}({\cal M})})(n\otimes
\nabla(\ell))]$.}
\medskip \noindent On en d\'eduit que $(1_{ \Omega^1}\otimes(\hbox{\goth L}\circ c))(\phi(\nabla_{\cal N})\otimes
1_{T^{1,1}({\cal M})})(n\otimes
\nabla(\ell))\in \Omega^1\otimes N $. Comme $\Omega^1$ est suppos\'e fid\`element plat \`a droite
sur $A$, et comme $n$ est arbitraire, ceci entra\^{\i}ne que $(\phi(\nabla_{\cal N})\otimes
1_{T^{1,1}({\cal M})})(N\otimes 
\nabla(\ell))\subset \Omega^1\otimes N \otimes LieStab_{T^{1,1}({\cal M})}(N)$, o\`u
$LieStab_{T^{1,1}({\cal M})}(N)$ d\'esigne la plus grande sous-$A$-alg\`ebre de Lie de
$T^{1,1}(M)$ telle que {\goth L}$(LieStab_{T^{1,1}({\cal M})}(N)\otimes N )\subset N$.
\par \noindent Or $\phi(\nabla_{\cal
N})$ est un isomorphisme (puisque $\cal N$ est rigide), donc $\;N\otimes 
\nabla(\ell)\subset N \otimes  \Omega^1 \otimes LieStab_{T^{1,1}({\cal M})}(N)$, d'o\`u finalement
$\nabla(\ell)\in \Omega^1 \otimes LieStab_{T^{1,1}({\cal M})}(N)$ et $(\nabla \ell)N \subset
\Omega^1\otimes N$.  
\par \noindent La seconde assertion en d\'ecoule, car tout sous-module \`a connexion de $LieStab
\lbrace {\cal N}\rbrace \break\subset T^{1,1}({\cal M})$ est stable sous l'action adjointe de
$LieStab \lbrace {\cal N}\rbrace$. 
\medskip \noindent {\bf 2.2.5. Le cas o\`u $A$ est de caract\'eristique $p>0$, et o\`u $\Omega^1$
est un bimodule commutatif.} Dans ce cas, l'anneau des constantes $k$ - qu'on suppose \^etre un
corps - contient $A^p$ (la commutativit\'e de $\Omega^1$ est essentielle ici). De plus,
${}^{\vee}\Omega^1$ est muni d'une $p$-structure. On peut d\'efinir l'op\'erateur de $p$-courbure
$R_{{\scriptstyle \nabla},p}$: pour tout $D\in {}^{\vee}\Omega^1,\;R_{{\scriptstyle \nabla}, p}(D) = (\nabla_{\scriptstyle
D})^p -\nabla_{(D^p)}$, un endomorphisme $A$-lin\'eaire de $M$.  
\par \noindent Supposons $\nabla$ {\it int\'egrable}. Alors
$R_{{\scriptstyle \nabla},p}$ est additif et $p$-lin\'eaire en $D$. Pla\c cons-nous 
dans la situation de I.2.1.1
ou I.2.2.1: plus pr\'ecis\'ement, $A$ est le corps de fonctions de $X$ lisse sur $k_{\scriptstyle
0}$ parfait (auquel cas le corps des constantes est $A^p$), ou bien $A$ est de la forme
${\bar{\cal O}}_x$ (auquel cas le corps des constantes est $k=k_{\scriptstyle 0}$). Les calculs
de [Kat70]5.2 montrent alors qu'on a $R_{{\scriptstyle \nabla}, p}(D)\in \;$End${\cal M}$, et que la
$p$-alg\`ebre de Lie engendr\'ee sur $A$ par les
$R_{{\scriptstyle \nabla}, p}(D)$ et leurs puissances $p^n$-i\`emes est ab\'elienne. 
\par Le cas de la
situation I.2.1.1 \'etant trait\'e dans [vdP95] (voir loc. cit. 6.6 pour le cas de dimension
sup\'erieure), nous nous limiterons {\it au cas $A={\bar{\cal O}}_x$}. Un foncteur fibre $\omega$
est alors donn\'e par  ${\cal N}\mapsto ({\cal N}_{A'})^{\nabla_{A'}}$ sur ${\scriptstyle <}{\cal
M}{\scriptstyle >}^{\scriptstyle\otimes}$, avec $A'={\cal
O}^{pd}_x$ (cf.1.4.2.), et tout objet de ${\scriptstyle <}{\cal
M}{\scriptstyle >}^{\scriptstyle\otimes}$ est rigide (cf. 1.3.2.$i),iii)$). La $p$-alg\`ebre de
Lie ab\'elienne ${\cal L}$ engendr\'ee sur $k$ par les endomorphismes horizontaux
$R_{{\scriptstyle \nabla}, p}(D)$ (et leurs puissances $p^n$-i\`emes) peut \^etre vue comme
sous-alg\`ebre de Lie de $LieGal({\cal M},\omega)$; elle est invariante sous l'action adjointe de
$Gal({\cal M},\omega)$. Notons $G({\cal L})$ le sous-sch\'ema en groupes ferm\'e infinit\'esimal
de hauteur $1$ (i.e. annul\'e par Frobenius relatif \`a $k$) de $Gal({\cal M},\omega)$ de
$p$-alg\`ebre de Lie $\cal L$ (cf [DG70]7, $n^{os}$ 3 et 4). On peut aussi voir ${\cal L}\otimes_k
A$ comme sous-alg\`ebre de Lie de $LieGal({\cal M},\hbox{\smit{oubli}})$. 
\proclaim Th\'eor\`eme 2.2.6. $i)$ $Gal({\cal M},\omega)= G({\cal L})$; en particulier,
$Gal({\cal M},\omega)$ est ab\'elien infinit\'esimal de hauteur $1$.
\hfill \break $ii)$ $LieGal({\cal M},\hbox{\smit{oubli}})$ est la $A$-alg\`ebre de Lie
engendr\'ee par les $R_{{\scriptstyle \nabla}, p}(D)$. 
\hfill \break $iii)$ Le torseur des solutions $\Sigma({\cal M},\omega)$ est un torseur trivial; en
particulier, \break $Gal({\cal M},\hbox{\smit{oubli}})\cong Gal({\cal
M},\omega)\otimes_{\scriptstyle k} A$.  
\par {\it Preuve.} D'apr\`es le crit\`ere de Chevalley, il suffit de montrer que toute droite
$k.\ell$ dans une somme finie $\oplus \; \omega (T^{i,j}({\cal M}))$ stable sous $G({\cal L})$
l'est aussi sous $Gal({\cal M},\omega)$. On peut du reste remplacer $k.\ell$ par une
quelconque puissance tensorielle non nulle; en particulier, par sa puissance $p$-i\`eme. Comme
$G({\cal L})$ est de hauteur $1$, on peut donc supposer $\ell$ fixe sous $G({\cal L})$. Alors
$\ell$ est annul\'e par $\cal L$. Il en est de m\^eme du $Gal({\cal
M},\omega)$-sous-espace de $\oplus \; \omega (T^{i,j}({\cal M}))$ engendr\'e par $\ell$,
puisque $\cal L$ est invariante sous l'action adjointe de $Gal({\cal
M},\omega)$. Le sous-module \`a connexion ${\cal N}\subset \oplus \; T^{i,j}({\cal
M})$ correspondant est donc annul\'e par les
$R_{{\scriptstyle \nabla}, p}(D)$. D'apr\`es I.2.2.1, ${\cal N}$ est par suite trivial. On
conclut que $\ell$ est fix\'e par $Gal({\cal
M},\omega)$, ce qui \'etablit $i)$. Les points $ii)$ et $iii)$ en r\'esultent
ais\'ement, et sont laiss\'es au lecteur.           
\medskip \noindent On peut voir 2.2.6.$ii)$ comme un analogue du th\'eor\`eme d'Ambrose-Singer
(I.1.2) en caract\'eristique $p$ (I.2.2).

\bigskip {\bf 2.3. Fonctorialit\'es.} 
\medskip \noindent {\bf 2.3.1.} Commen\c cons par quelques pr\'eliminaires sur la cat\'egorie
$Ind{\scriptstyle <}{\cal M}{\scriptstyle >}^{\scriptstyle\otimes}$ des $ind$-objets de ${\scriptstyle <}{\cal M}{\scriptstyle >}^{\scriptstyle\otimes}$. Elle admet les
descriptions \'equivalentes suivantes: 
\par \noindent $i)$ c'est la cat\'egorie des (petits) syst\`emes inductifs filtrants $({\cal
M}_{\alpha})$ de ${\scriptstyle <}{\cal M}{\scriptstyle >}^{\scriptstyle\otimes}$, avec $Mor(({\cal
M}_{\alpha}),({\cal
N}_{\beta}))={\displaystyle
\limproj_{\scriptstyle\alpha}\limind_{\scriptstyle\beta}}\; Mor({\cal M}_{\alpha},{\cal
N}_{\beta})$,
\par \noindent $ii)$ c'est la cat\'egorie des foncteurs contravariants ${\scriptstyle <}{\cal
M}{\scriptstyle >}^{\scriptstyle\otimes}\rightarrow Ens$ qui sont la (petite) limite inductive filtrante de foncteurs
repr\'esentables $h_{{\cal
M}_{\alpha}}$.
\par \noindent Le passage de $i)$ \`a $ii)$ est $({\cal
M}_{\alpha})\mapsto \limind h_{{\cal M}_{\alpha}}\;$ (cf.[De89] 4).  
\par \noindent Comme ${\scriptstyle <}{\cal M}{\scriptstyle >}^{\scriptstyle\otimes}$ est $k$-lin\'eaire ab\'elienne mono\"{\i}dale, il en
est de m\^eme de \break $Ind{\scriptstyle <}{\cal M}{\scriptstyle >}^{\scriptstyle\otimes}$, et ${\scriptstyle <}{\cal M}{\scriptstyle >}^{\scriptstyle\otimes}$ s'identifie \`a une
sous-cat\'egorie strictement pleine de $Ind{\scriptstyle <}{\cal M}{\scriptstyle >}^{\scriptstyle\otimes}$. 
\par \noindent Par ailleurs, la cat\'egorie des $A$-modules \`a connexion poss\`ede de mani\`ere
\'evidente des (petites) limites inductives filtrantes, et on a $h_{\limind {\cal
M}_{\alpha} }= \limind h_{{\cal M}_{\alpha}}$. On en d\'eduit que le foncteur 
\par \centerline{${\limind }:\; Ind{\scriptstyle <}{\cal M}{\scriptstyle >}^{\scriptstyle\otimes} \rightarrow \lbrace A-modules\; {\hbox
{\it \`a }}\; connexion \rbrace$}
\par \noindent est pleinement fid\`ele (et d'ailleurs mono\"{\i}dal), ce qui nous permet
d'identifier dor\'enavant les $ind$-objets de ${\scriptstyle <}{\cal M}{\scriptstyle >}^{\scriptstyle\otimes}$ \`a des $A$-modules \`a 
connexion. 
\par \noindent Dans la situation 2.1.1-2.1.2, tout objet de ${\cal C}={\scriptstyle <}{\cal M}{\scriptstyle >}^{\scriptstyle\otimes}$ est
noeth\'erien, et on en d\'eduit ([De89] 4.2.1) que tout $ind$-objet de ${\scriptstyle <}{\cal M}{\scriptstyle >}^{\scriptstyle\otimes}$ est
(petite) limite ordonn\'ee filtrante de sous-objets qui sont objets de ${\scriptstyle <}{\cal M}{\scriptstyle >}^{\scriptstyle\otimes}$.
Les objets de ${\scriptstyle <}{\cal M}{\scriptstyle >}^{\scriptstyle\otimes}$ sont les objets noeth\'eriens de $Ind{\scriptstyle <}{\cal M}{\scriptstyle >}^{\scriptstyle\otimes}$.  
\medskip \noindent {\bf 2.3.2. Fonctorialit\'e en $\cal M$.} Dans la situation 2.1.2,
consid\'erons un objet rigide $\cal N$ de ${\scriptstyle <}{\cal M}{\scriptstyle >}^{\scriptstyle\otimes}$. Supposons que tout objet de
${\scriptstyle <}{\cal N}{\scriptstyle >}^{\scriptstyle\otimes}$ soit quotient d'un objet rigide. On peut appliquer 2.1.1.$i)$ \`a la
sous-cat\'egorie pleine ${\scriptstyle <}{\cal N}{\scriptstyle >}^{\scriptstyle\otimes}$ de ${\scriptstyle <}{\cal M}{\scriptstyle >}^{\scriptstyle\otimes}$. On obtient alors un
homomorphisme canonique 
\medskip \centerline{$Gal({\cal M},\omega)\rightarrow Gal({\cal N},\omega)$.}

\proclaim Lemme 2.3.3. Si $k$ est r\'egulier de dimension $\leq 1$, l'homomorphisme de big\`ebres
${\cal O}(Gal({\cal N},\omega))\rightarrow {\cal O}(Gal({\cal M},\omega))$ est
injectif. Si $k$ est un corps, $Gal({\cal M},\omega)\rightarrow Gal({\cal N},\omega)$
est fid\`element plat.    
\par {\it Preuve.} On \'etend $\omega$ \`a la cat\'egorie des $ind$-objets de ${\scriptstyle <}{\cal
M}\; {\rm ou} \;{\cal N}{\scriptstyle >}^{\scriptstyle\otimes}$, en consid\'erant des repr\'esentations non n\'ecessairement
de type fini de $Gal({\cal M}\; {\rm ou} \;{\cal N},\omega)$, cf. [Saa72]II.2.3.4. Alors
le foncteur $X \in Ind{\scriptstyle <}{\cal M}{\scriptstyle >}^{\scriptstyle\otimes} \mapsto Hom_{\scriptstyle k}(\omega(X),k)$ (resp. $X \in Ind{\scriptstyle <}{\cal
N}{\scriptstyle >}^{\scriptstyle\otimes} \mapsto Hom_{\scriptstyle k}(\omega(X),k)$) est repr\'esentable par un objet $B_{\cal M}$ (resp.
$B_{\cal N}$), et on a un morphisme naturel $B_{\cal N}\rightarrow B_{\cal M}$ dans $Ind{\scriptstyle <}{\cal
M}{\scriptstyle >}^{\scriptstyle\otimes}$; l'image de ce dernier par $\omega$ n'est autre que ${\cal O}(Gal({\cal
N},\omega))\rightarrow {\cal O}(Gal({\cal M},\omega))$, cf. loc. cit. 2.3.2.1. 
\par \noindent Soit $X$ un sous-objet (dans $Ind{\scriptstyle <}{\cal N}{\scriptstyle >}^{\scriptstyle\otimes}$) du noyau de $B_{\cal
N}\rightarrow B_{\cal M}$, tel que $\omega(X)$ soit de type fini sur $k$; on peut voir $X$
comme un objet de ${\scriptstyle <}{\cal N}{\scriptstyle >}^{\scriptstyle\otimes}$. Les identifications $Hom_{\scriptstyle k}(\omega(X),k)=Mor_{Ind{\scriptstyle <}{\cal
N}{\scriptstyle >}^{\scriptstyle\otimes}}(X,B_{\cal N})=Mor_{Ind{\scriptstyle <}{\cal M}{\scriptstyle >}^{\scriptstyle\otimes}}(X,B_{\cal M})$ montrent que
$Hom_{\scriptstyle k}(\omega(X),k)=0$. Comme $B_{\cal N}$ est $k$-plat, on en d\'eduit que $X=0$ si $k$ est
r\'egulier de dimension $\leq 1$. Or $\omega(Ker(B_{\cal N}\rightarrow B_{\cal M}))$ est limite
inductive de tels $\omega(X)$ (cf. [Se93]1.4), donc $B_{\cal N}\hookrightarrow B_{\cal M}$.
\par \noindent La seconde assertion du lemme d\'ecoule de la premi\`ere (qui, dans le cas o\`u
$k$ est un corps, s'obtient directement \`a partir de 2.1.1. $ii)$; voir aussi [Saa72]II.4.3.2.). 

\medskip \noindent {\bf 2.3.4. Changement d'anneau diff\'erentiel.} Soit $(A,d)\rightarrow
(\tilde{A}, \tilde{d})$ un morphisme d'anneaux diff\'erentiels (cf.II.4.5.). On ne suppose {\it pas} que la
$k$-alg\`ebre $\tilde{C}$ des constantes de $\tilde{A}$ soit \'egale \`a $k$. 
\par \noindent Nous aurons besoin de la cat\'egorie $\tilde{C}$-lin\'eaire ab\'elienne
mono\"{\i}dale
$Ind{\scriptstyle <}{\cal M}{\scriptstyle >}^{\scriptstyle\otimes}_{(\tilde{C})}$ form\'ee des 
$ind$-objets de ${\scriptstyle <}{\cal M}{\scriptstyle >}^{\scriptstyle\otimes}$ munis
d'une action de $\tilde{C}$; sous les hypoth\`eses de 2.1.1 (avec ${\cal C} ={\scriptstyle <}{\cal
M}{\scriptstyle >}^{\scriptstyle\otimes}$), cette cat\'egorie est $\otimes$-\'equivalente \`a la cat\'egorie des
repr\'esentations sur $\tilde{C}$ de
$Gal({\cal M},\omega)\otimes_{\scriptstyle k} \tilde{C}$, cf. [Saa72]II.1.5, 2.0, III.1.1. 
\medskip \noindent Pour simplifier, nous supposerons que $\tilde{C}\otimes_{\scriptstyle k} \tilde{C} \cong
\tilde{C}$; c'est le cas si
$\tilde{C}$ est une {\it localisation} ou bien un {\it quotient} de $k$. On peut alors
consid\'erer
$Ind{\scriptstyle <}{\cal M}{\scriptstyle >}^{\scriptstyle\otimes}_{(\tilde{C})}$ comme une 
sous-cat\'egorie pleine de $Ind{\scriptstyle <}{\cal
M}{\scriptstyle >}^{\scriptstyle\otimes}$, et m\^eme comme une sous-cat\'egorie mono\"{\i}dale en vertu des isomorphismes
canoniques $M'\otimes_{A\otimes 
\tilde{C}} M"\cong M'\otimes_{\scriptstyle A}
M"$. On dispose d'autre part d'un foncteur mono\"{\i}dal
\par \centerline{ $u:\;Ind{\scriptstyle <}{\cal
M}{\scriptstyle >}^{\scriptstyle\otimes}_{(\tilde{C})}\;\longrightarrow \; Ind{\scriptstyle <}{\cal
M}_{\tilde{A}}{\scriptstyle >}^{\scriptstyle\otimes}, \;\; N\mapsto {\tilde{A}}\otimes_{A\otimes_{\scriptstyle k} \tilde{C}}N \cong
N_{\tilde{A}}$.}
\medskip \noindent Supposons de plus $\tilde{A}$ fid\`element plat sur $A\otimes_{\scriptstyle k}
\tilde{C}$, de sorte que $u$ est fid\`ele et exact.
\par \noindent Notons $u^{-1}{\scriptstyle <}{\cal
M}_{\tilde{A}}{\scriptstyle >}^{\scriptstyle\otimes}$ la sous-cat\'egorie pleine de 
$Ind{\scriptstyle <}{\cal
M}{\scriptstyle >}^{\scriptstyle\otimes}_{(\tilde{C})}$ form\'ee des objets dont l'image par 
$\; u \;$ est dans ${\scriptstyle <}{\cal M}_{\tilde{A}}{\scriptstyle >}^{\scriptstyle\otimes}$.
C'est la sous-cat\'egorie pleine des objets noeth\'eriens; elle est
$\otimes$-\'equivalente \`a la cat\'egorie des repr\'esentations de type fini sur $\tilde{C}$ de
$Gal({\cal M},\omega)\otimes_{\scriptstyle k} \tilde{C}$.   
\medskip \noindent Supposons donn\'e en outre un foncteur mono\"{\i}dal
$\tilde{\omega}$ sur ${\scriptstyle <}{\cal M}_{\tilde{A}}{\scriptstyle >}^{\scriptstyle\otimes}$ comme en 2.1.1, tel que 
$\;\tilde{\omega} \circ u \;$ soit isomorphe \`a la restriction de $\omega \otimes_{\scriptstyle
k} id_{\tilde{C}}$ \`a
$u^{-1}{\scriptstyle <}{\cal M}_{\tilde{A}}{\scriptstyle >}^{\scriptstyle\otimes}$. Le foncteur mono\"{\i}dal
\par \centerline{ $u:\;u^{-1}{\scriptstyle <}{\cal M}_{\tilde{A}}{\scriptstyle >}^{\scriptstyle\otimes}\;\longrightarrow \;{\scriptstyle <}{\cal
M}_{\tilde{A}}{\scriptstyle >}^{\scriptstyle\otimes}, \;\; N\mapsto N_{\tilde{A}}$}
\par \noindent induit un homomorphisme 
\medskip \centerline{$Gal({\cal M}_{\tilde{A}},\tilde{\omega}) \rightarrow Gal({\cal
M},\omega)\otimes_{\scriptstyle k} \tilde{C}$.}
\par \proclaim Lemme 2.3.5. Cet homomorphisme est un monomorphisme (en particulier, c'est une
immersion ferm\'ee si $\;\tilde{C}\;$ est un corps). C'est un isomorphisme si $(\tilde{A},
\tilde{d})=(A,d)\otimes_{\scriptstyle k} \tilde{C}$ (localisation ou sp\'ecialisation des constantes).   
\par {\it Preuve.} On a clairement  
$\u{Aut}^{\scriptstyle\otimes}((\omega \otimes_{\scriptstyle k} id_{\tilde{C}})_{\mid u^{-1}{\scriptstyle <}{\cal
M}_{\tilde{A}}{\scriptstyle >}^{\scriptstyle\otimes}}) = \u{Aut}^{\scriptstyle\otimes}
(\tilde{\omega}\circ u)\hookleftarrow 
\u{Aut}^{\scriptstyle\otimes}(\tilde{\omega})_{\mid {\scriptstyle <}{\cal M}_{\tilde{A}}{\scriptstyle >}^{\scriptstyle\otimes}}$, puisque tout objet de
${\scriptstyle <}{\cal M}_{\tilde{A}}{\scriptstyle >}^{\scriptstyle\otimes}$ est sous-quotient d'un objet de $u(u^{-1}{\scriptstyle <}{\cal
M}_{\tilde{A}}{\scriptstyle >}^{\scriptstyle\otimes})$; d'o\`u la premi\`ere assertion. La seconde provient de ce que si
$(\tilde{A}, \tilde{d})=(A,d)\otimes_{\scriptstyle k} \tilde{C}$, alors $\;u: \;N\mapsto N_{\tilde{A}}\cong N$
est une \'equivalence de cat\'egories mono\"{\i}dales. 

\proclaim Lemme 2.3.6. Supposons $A$ noeth\'erien d'anneau total de fractions $Q(A)$ semi-simple, $\Omega^1$
fid\`ele et projectif de type fini \`a droite et tel que $\;\Omega^1
\otimes_{\scriptstyle A}Q(A) \cong Q(A)\otimes_{\scriptstyle A}
\Omega^1$, et $(A,d)$ simple. Alors le foncteur de localisation
$\;{\scriptstyle <}{\cal M}{\scriptstyle >}^{\scriptstyle\otimes}\;\rightarrow
\;{\scriptstyle <}{\cal M}_{\scriptstyle Q(A)}{\scriptstyle >}^{\scriptstyle\otimes},\;\;\;{\cal N}
\mapsto {\cal N}_{\scriptstyle Q(A)}$ 
est une \'equivalence de cat\'egories mono\"{\i}dales. Si l'on dispose d'un foncteur
$\omega$ comme en 2.1.1, cela donne lieu \`a un isomorphisme 
\hfill\break \centerline{$Gal({\cal M}_{\scriptstyle Q(A)},\omega)
\cong Gal({\cal M},\omega)$.}
\par {\it Preuve.} Sous ces hypoth\`eses, tout objet de ${\scriptstyle <}{\cal M}{\scriptstyle
>}^{\scriptstyle \otimes}$ est rigide (II.5.3.2), et la pleine fid\'elit\'e r\'esulte de II.5.2.2. Pour la
surjectivit\'e essentielle, il suffit de montrer que tout sous-objet $\tilde{\cal N}$ d'une somme finie
$\oplus \; T^{i,j}({\cal M})_{\scriptstyle Q(A)}$ est de la forme
${\cal N}_{\scriptstyle Q(A)}$, pour un sous-objet ${\cal N}$ convenable de $\oplus \; T^{i,j}({\cal
M})$. Il suffit de prendre ${\cal N}=\tilde{\cal N}\cap (\oplus \; T^{i,j}({\cal
M}))$.
\medskip \noindent {\bf 2.3.7. Changement de constantes.} Soit $k \rightarrow
\tilde{k}$ un homomorphisme d'anneaux. Un $A\otimes {\tilde{k}}$-module
\`a connexion (relativement \`a l'anneau diff\'erentiel $(A,d)\otimes_{\scriptstyle k} {\tilde{k}}\;$) n'est rien
d'autre qu'un $A$-module \`a connexion $\cal N$ muni de ${\tilde{k}}\rightarrow \;$End$(\cal N)$,
ce qui fournit une \'equivalence tautologique de cat\'egories mono\"{\i}dales
${\tilde{k}}$-lin\'eaires 
$Ind{\scriptstyle <}{\cal M}{\scriptstyle >}^{\scriptstyle\otimes}_{(\tilde{k})}\;\cong \;Ind{\scriptstyle <}{\cal
M}\otimes_{\scriptstyle k}{\tilde{k}}{\scriptstyle >}^{\scriptstyle\otimes}$.
\par \noindent Le foncteur $\omega$ s'\'etend en un foncteur
mono\"{\i}dal $k$-lin\'eaire, fid\`ele et exact, encore not\'e $\omega:\;Ind {\scriptstyle <}{\cal
M}{\scriptstyle >}^{\scriptstyle\otimes}\rightarrow \lbrace k-modules \rbrace$. Il induit un foncteur mono\"{\i}dal
$\tilde{k}$-lin\'eaire $\tilde{\omega}:\; Ind{\scriptstyle <}{\cal
M}{\scriptstyle >}^{\scriptstyle\otimes}_{(\tilde{k})}\;\cong \;Ind{\scriptstyle <}{\cal
M}\otimes_{\scriptstyle k}{\tilde{k}}{\scriptstyle >}^{\scriptstyle\otimes}\rightarrow \lbrace {\tilde{k}}-modules \rbrace$, encore exact
et fid\`ele. On en d\'eduit:
\proclaim Lemme 2.3.8. Sous ces hypoth\`eses, $Gal({\cal
M}\otimes_{\scriptstyle k}{\tilde{k}},\tilde{\omega}) \cong Gal({\cal M},\omega)\otimes_{\scriptstyle k} \tilde{k}$.
\par (Voir l'exemple 2.1.5.)
\par \noindent Par ailleurs le diagramme suivant de foncteurs est commutatif
\medskip \centerline{$\matrix{{\scriptstyle <}{\cal
M}{\scriptstyle >}^{\scriptstyle\otimes} & \;{\buildrel{\omega}\over
\longrightarrow}\; & {\lbrace k-modules \rbrace } \cr 
{\scriptstyle{\otimes_{\scriptstyle k} {\tilde k}}}\downarrow \hphantom{\scriptstyle \otimes_{\scriptstyle k} {\tilde k}} & &
\downarrow {\scriptstyle \otimes_{\scriptstyle k} {\tilde k}} 
\cr  {Ind{\scriptstyle <}{\cal M}{\scriptstyle >}^{\scriptstyle\otimes}_{(\tilde{k})}} & \;
{\buildrel{\tilde{\omega}}\over
\longrightarrow}\; & \lbrace {\tilde k}-modules \rbrace }$}

\bigskip {\bf 3. Le th\'eor\`eme de sp\'ecialisation.}

\medskip {\bf 3.1.} Dans ce paragraphe, nous d\'etaillons un cas particulier des fonctorialit\'es
pr\'ec\'edentes. Comme dans l'exemple 2.1.3, consid\'erons une
extension diff\'erentielle $(A', d')$ de $(A,d)$, simple par couches. On suppose $C'= C=k$ {\it
de Dedekind} (par exemple un anneau de valuation discr\`ete),
et $A'$ fid\`element plat sur $A$. 
\par \noindent On se donne un $A$-module \`a connexion ${\cal M}$ projectif de type fini sur $A$,
soluble dans $A'$. On consid\`ere le foncteur $ \omega_{\scriptstyle A'}: \;{\cal N}\mapsto
({\cal N}_{A'})^{\nabla_{A'}}$ sur $\;{\scriptstyle <}{\cal M}{\scriptstyle >}^{\scriptstyle\otimes}$, et le groupe de Galois
diff\'erentiel $Gal({\cal M},\; \omega_{\scriptstyle A'})$ (un $k$-sch\'ema en groupes plat, non
n\'ecessairement de type fini). 
\par \noindent Si $\tilde{k}$ est une localisation de $k$ (par exemple son
corps de fractions $\kappa(k)$) ou bien le quotient de $k$ par un id\'eal maximal, alors les
hypoth\`eses de 2.1.1. sont encore v\'erifi\'ees si l'on remplace
$(A,d)$ par $(A,d)\otimes \tilde{k}$, $(A', d')$ par
$(A', d')\otimes \tilde{k}$, etc..., et $\omega_{\scriptstyle
A'}$ par $\omega_{\scriptstyle
A'}\otimes 1\cong \omega_{\scriptstyle
A'\otimes {\tilde k}}$ (cf. 1.3.2.$iv)$). On dispose donc du
$\tilde{k}$-sch\'ema en groupes de Galois diff\'erentiel $Gal({\cal M}\otimes_{\scriptstyle k}
\tilde{k},\;\omega_{\scriptstyle
A'\otimes {\tilde k}})$.  

\proclaim Th\'eor\`eme 3.1.1. On a un isomorphisme canonique 
\medskip $Gal({\cal M}\otimes_{\scriptstyle k} \kappa(k),\;\omega_{\scriptstyle
A'\otimes \kappa(k)})\; {\buildrel{\sim}\over
\rightarrow}\; Gal({\cal M},\omega_{\scriptstyle
A'})\otimes_{\scriptstyle k} \kappa(k).$
\medskip \noindent De plus, pour tout id\'eal maximal  {\goth p} 
de $k$, on a un isomorphisme canonique 
\medskip $Gal({\cal M}\otimes_{\scriptstyle k} k/{\hbox{\goth p}},\;\omega_{\scriptstyle
A'\otimes k/{\hbox{\smgoth p}}})\; {\buildrel{\sim}\over
\rightarrow}\; Gal({\cal M},\omega_{\scriptstyle
A'})\otimes_{\scriptstyle k} k/{\hbox{\goth p}}.$
\par \noindent Cela d\'ecoule tant de 2.3.5 (seconde assertion) que de
2.3.8, avec $\tilde{C}=\tilde{k}=\kappa(k)$, resp. $\tilde{C}=\tilde{k}=k/${\goth p}. 

\medskip {\bf 3.2.} Pour \'eviter les $k$-sch\'emas en groupes non
de type fini et rendre l'\'enonc\'e plus tangible, introduisons le sous-groupe ferm\'e 
de $GL(\omega_{\scriptstyle A'}({\cal M}))$ adh\'erence de Zariski de $Gal({\cal
M}\otimes_{\scriptstyle k}
\kappa(k),\omega_{\scriptstyle
A'\otimes \kappa(k)})$.
C'est un sch\'ema en groupes plat {\it de type fini} sur
$k$ \'egal au stabilisateur $Stab \lbrace 
\omega_{\scriptstyle A'}({\cal N}) \rbrace$ des $\;\omega_{\scriptstyle
A'}({\cal N})$, pour
tout sous-objet $\;\cal N \;$ d'une somme finie quelconque $\;\oplus \; T^{i,j}({\cal M})$ tels
que $\;\omega_{\scriptstyle A'}({\cal N})\;$ soit facteur direct dans le $k$-module $\;\oplus
\;\omega_{\scriptstyle A'}( T^{i,j}({\cal
M}))$ (cf. [DG70]II.1.3.6). Sa fibre g\'en\'erique est $Gal({\cal M}\otimes_{\scriptstyle k}
\kappa(k),\omega_{\scriptstyle
A'\otimes\kappa(k)})$. On a alors

\proclaim Corollaire 3.2.1. Pour tout id\'eal maximal  {\goth p}  de $k$, on a un immersion
ferm\'ee canonique 
$\;Gal({\cal M}\otimes_{\scriptstyle k} k/{\hbox{\goth p}},\;\omega_{\scriptstyle
A'\otimes k/{\hbox{\smgoth p}}})\; \hookrightarrow\; Stab \lbrace \omega_{\scriptstyle
A'}({\cal N}) \rbrace\otimes_{\scriptstyle k} k/{\hbox{\goth p}}.$ En
particulier, 
\break $dim\;Gal({\cal M}\otimes_{\scriptstyle k} k/{\hbox{\goth p}},\;\omega_{\scriptstyle
A'\otimes k/{\hbox{\smgoth p}}}) \;\leq \;dim\;Gal({\cal M}\otimes_{\scriptstyle k}
\kappa(k),\;\omega_{\scriptstyle
A'\otimes \kappa(k)})$.
\par \noindent {\bf 3.2.2. Remarque.} Cet \'enonc\'e ne fournit pas exactement un th\'eor\`eme
de semicontinuit\'e du fait que la fonction {\goth p}$\;\mapsto \; dim\;Gal({\cal M}
\otimes_{\scriptstyle k} k/{\hbox{\goth p}},\;\omega_{\scriptstyle
A'\otimes k/{\hbox{\smgoth p}}})$ n'est pas alg\'ebriquement constructible en
g\'en\'eral. On peut n\'eanmoins montrer la constructibilit\'e - et partant la semicontinuit\'e -
pour la topologie {\it classique}, dans la situation analytique de 1.4.1. 
\medskip Dans la situation alg\'ebrique de 1.4.1, on obtient l'\'enonc\'e suivant,
d\^u \`a O. Gabber (cf. [Kat90]2.4), par une autre voie:
\proclaim Corollaire 3.2.3. Soient $X$ un sch\'ema lisse s\'epar\'e de type fini
 \`a fibres \break
g\'eom\'etriquement connexes sur ${\bf C}[[h]]$, $x$ un $k$-point de $X$, $E$ un ${\cal
O}_X$-module localement libre de rang fini, muni d'une connexion int\'egrable relative $
E \rightarrow \Omega^1_{X/{\bf C}[[h]]}\otimes_{{\cal O}_X} E$. Alors le groupe de Galois
diff\'erentiel de $E/hE$ (point\'e en $x$) est contenu dans la fibre sp\'eciale de la ${\bf
C}[[h]]$-adh\'erence de Zariski du groupe de Galois diff\'erentiel de $E \otimes {\bf C}((h))$
(point\'e en $x$).
\par
Le corollaire 3.2.1 s'applique aussi bien \`a la situation de
{\it confluence d'\'equations aux ($q$-)diff\'erences} \'etudi\'ee en 1.4.3, 1.4.4, lorsque
$k$ est de dimension un: 
\proclaim Corollaire 3.2.4. Dans cette situation, le groupe de Galois
diff\'erentiel du syst\`eme diff\'erentiel obtenu par confluence est contenu dans la fibre
sp\'eciale de la $k$-adh\'erence de
Zariski du groupe de Galois diff\'erentiel du syst\`eme aux \break ($q$-)diff\'erences sur le
corps de fractions de $k$.
\par D'apr\`es 2.3.6, on peut aussi bien consid\'erer l'\'equation aux ($q$-)diff\'erences comme
d\'efinie sur $A$ ou sur le corps de fractions de $A$. 

\medskip {\bf 3.3. Une application.} Ce dernier r\'esultat permet de ``calculer" le
groupe de Galois diff\'erentiel de certaines \'equations $q$-hyperg\'eom\'etriques. A titre
d'illustration, consid\'erons l'\'equation aux $q$-diff\'erences d'ordre $r > 1$ 
\medskip \centerline{$(\ast)_{{\u\alpha},{\u\beta},q}\;\;\;\;\; \Pi_{j=1}^{j=r}
(q^{\beta_j-1}z\delta_q + {q^{\beta_j-1}-1\over q-1})y= z\Pi_{i=1}^{i=r} (q^{\alpha_i}z\delta_q +
{q^{\alpha_i}-1\over q-1})y\;,
\;\;(\alpha_i, \beta_j \in {\bf C},\beta_r=1)$}
\medskip \noindent satisfaite par la s\'erie
$q$-hyperg\'eom\'etrique (cf. [GR90], p.27)
\medskip \centerline{$_r\phi_{r-1}(\matrix{q^{\alpha_1}&\ldots, &q^{\alpha_r} \cr q^{\beta_1}
&\ldots & q^{\beta_{r-1}}}; q ; z)= \sum_{n\geq 0}\;{(q^{\alpha_1};q)_n \ldots
(q^{\alpha_r};q)_n\over (q^{\beta_1};q)_n\ldots (q^{\beta_{r-1}};q)_n.(q;q)_n}\;z^n\;.$}
\medskip \noindent Ici, une expression comme $q^{\alpha}$ est \`a interpr\'eter comme
l'\'el\'ement $(1+(q-1))^{\alpha}$ de ${\bf C}[[q-1]]$, et
$(q^{\alpha};q)_n=(1-q^{\alpha})(1-q^{\alpha +1})\ldots (1-q^{\alpha +n-1})$.
\par \noindent Cette \'equation aux $q$-diff\'erences conflue vers
l'\'equation diff\'erentielle ordinaire 
\medskip \centerline{$(\ast)_{{\u\alpha},{\u\beta}}\;\;\;\;\;\;\;\;\;\;\Pi_{j=1}^{j=r} (z{d\over
dz} +
\beta_j-1)y= z\Pi_{i=1}^{i=r} (z{d\over dz} + \alpha_i)y$}
\medskip \noindent satisfaite par la s\'erie
hyperg\'eom\'etrique $_rF_{r-1}(\matrix{\alpha_1&\ldots, &\alpha_r \cr \beta_1 &\ldots &
\beta_{r-1}}; z).$
\medskip \noindent On conna\^{\i}t pr\'ecis\'ement les conditions sur les param\`etres
$\alpha_i,\beta_j$ qui entra\^{\i}nent que le groupe de Galois diff\'erentiel de
$(\ast)_{{\u\alpha},{\u\beta}}$ {\it contient} $SL_{r,{\bf C}}$ (cf. [Beu92]4). C'est par exemple
le cas lorsque les conditions suivantes sont simultan\'ement satisfaites (toutes sauf la
derni\`ere sont d'ailleurs n\'ecessaires): 
\proclaim Conditions. $\bullet$ Pour tout $i \neq j$, $\alpha_i \not\equiv \beta_j$ mod.$\bf Z$,
\medskip \noindent $\bullet$ pour aucun entier naturel $d > 1$, on n'a \hfill \break
$\lbrace \alpha_1 + {1 \over d}, \ldots ,
\alpha_r + {1 \over d}\rbrace \equiv \lbrace \alpha_1 , \ldots ,
\alpha_r \rbrace,\; \lbrace \beta_1 +{1 \over d}, \ldots ,
\beta_r + {1 \over d}\rbrace \equiv \lbrace \beta_1 , \ldots ,
\beta_r \rbrace $ mod.$\bf Z$,
\medskip \noindent $\bullet$ pour aucun entier naturel $d < r$, et aucun couple $(x,y)\in
{\bf C}^2$, on n'a \hfill \break $\lbrace x, x + {1 \over d},\ldots , x + {d-1 \over d}, \; y,\; y
+ {1 \over r-d}, \ldots , y + {r-d-1 \over r-d} \rbrace \equiv \lbrace ... \alpha_i,...
\rbrace\;$ (resp. $\lbrace ... \beta_i,... \rbrace$),
\hfill \break $\lbrace {ax+by \over a+b}, {ax+by \over a+b} + {1 \over r}\;,\ldots , {ax+by \over
a+b} + {r-1 \over r} \rbrace \equiv \lbrace ... \beta_i,...\rbrace\;$ (resp. $\lbrace ...
\alpha_i,... \rbrace\;$) mod.$\bf Z$,
\medskip \noindent $\bullet$ s'il existe $x \in \bf C$ tel que les $\;\alpha_i + x \;$ et
$\;\beta_i + x\;$ soient tous rationnels, alors il existe un entier $\; n \;$ premier au
d\'enominateur commun de ces nombres, tel que les ensembles $\lbrace e^{2\pi in(\alpha_i + x)}
\rbrace$ et $\lbrace e^{2\pi in(\beta_i + x)}\rbrace$ ne s'entrelacent pas sur le cercle unit\'e,
\par \noindent $\bullet$ $\Sigma (\alpha_i - \beta_i)$ n'est pas demi-entier.
\par 
\proclaim Th\'eor\`eme 3.3.1. Sous ces conditions, le groupe de Galois diff\'erentiel de
l'\'equation
$q$-hyperg\'eom\'etrique $\;(\ast)_{{\u\alpha},{\u\beta},q}\;$ sur
$\;{\bf C}((q-1))(z)\;$ est $\;GL_{r,{\bf C}((q-1))}$.

\par En effet, d'apr\`es 3.2.4 (compl\'et\'e par 2.3.6), ce groupe de Galois diff\'erentiel
contient $\;SL_{r,{\bf C}((q-1))}$. Il reste donc \`a examiner l'\'equation aux $q$-diff\'erences
d'ordre un satisfaite par le d\'eterminant de Casorati (cf. II.4.3). C'est 
\medskip \centerline{$Cas (\vec y)(qz)={(-1)^r q^{\Sigma (1-\beta_j)}(1-z)\over 1-q^{\Sigma
(1+\alpha_j-\beta_j)}z}\;Cas (\vec y)(z)$}
\par \noindent En testant la r\'egularit\'e en $0$ et \`a l'infini (cf. [vdPS97]12.19), on
v\'erifie qu'elle correspond \`a une connexion triviale si et seulement si
$r$ est pair et ${\Sigma (1-\beta_j)}={\Sigma \alpha_j}=0$; or ces \'egalit\'es sont
exclues puisque $\Sigma (\alpha_i - \beta_i)$ n'est pas demi-entier.
 \medskip Lorsque les param\`etres ${\u \alpha}, {\u \beta}$ sont tous rationnels, de
d\'enominateur commun $N$, on peut voir $\;(\ast)_{{\u\alpha},{\u\beta},q}\;$ comme d\'efinie sur
${\bf Q}(q^{\scriptstyle 1/N},z)$, qu'on peut plonger dans $\bf C$ en donnant \`a $\;q\;$ une
valeur complexe transcendante arbitraire. On d\'eduit a fortiori du th\'eor\`eme, joint \`a
2.3.8, que le groupe de Galois diff\'erentiel de
$\;(\ast)_{{\u\alpha},{\u\beta},q}\;$ sur ${\bf C}(z)$ est $\;GL_{r,{\bf C}}$. Je ne connais pas de
d\'emonstration directe de ce r\'esultat: l'outil essentiel dans la d\'etermination du groupe de Galois
diff\'erentiel hyperg\'eom\'etrique, \`a savoir la pseudo-r\'eflexion donn\'ee par la monodromie locale au
point $1$, ne se transporte pas au cas $q$-hyperg\'eom\'etrique. Pour une approche analytique de ce type
de confluence, voir [Sau99].

\medskip\bigskip {\bf 4. Extensions de Picard-Vessiot.}
\medskip {\it Nous supposons d\'esormais que l'anneau des constantes $k$ est un corps.}
\medskip {\bf 4.1. D\'efinition.}
\par \noindent On consid\`ere un $A$-module \`a connexion
${\cal M}=(M,\nabla)$ rigide au sens de II.4.2: $M$ est projectif de type fini et la volte
$\phi(\nabla)$ est un isomorphisme. Il admet un dual $\check{\cal M}$. On rappelle que ${\scriptstyle <}{\cal
M}{\scriptstyle >}^{\scriptstyle\otimes}$ d\'esigne la sous-cat\'egorie strictement pleine de la cat\'egorie des $A$-modules
\`a connexion form\'ee des sous-quotients des sommes finies $\oplus \; {\cal M}^{\otimes
i}\otimes {\check{\cal M}}^{\otimes j}$.  
\par \noindent Soit $A'{\buildrel{d'}\over \rightarrow}\Omega'^1\cong \Omega^1\otimes A'$ une
extension diff\'erentielle. On dispose des $A'$-modules \`a connexion ${\cal M}_{A'}$ et
$\check{\cal M}_{A'}$. On note $\omega_{\scriptstyle
A'}$ le foncteur ${\cal N} \mapsto N_{A'}^{\nabla_{A'}}=
Ker_{N_{A'}} (\nabla_{A'})$ sur ${\scriptstyle <}{\cal M}{\scriptstyle >}^{\scriptstyle\otimes}$. On dispose des sous-$A$-modules de type
fini $\langle M,\omega_{\scriptstyle
A'}(\check{\cal M})\rangle$ et $\langle \check{M},\omega_{\scriptstyle
A'}({\cal M})\rangle$ de $A'$. 
\proclaim D\'efinition 4.1.1. On dit que $(A',d')$ est une
extension \hbox{\rm{de Picard-Vessiot}} (enti\`ere) de $(A,d)$
pour ${\cal M}$ si
\hfill \break $i)$ $A'$ est fid\`element plate sur $A$,
\hfill \break $ii)$ $(A',d')$ est simple,
\hfill \break $iii)$ l'anneau des constantes de $(A',d')$ est le corps $\;k$,
\hfill \break $iv)$ ${\cal M}$ est soluble dans $(A',d')$,
\hfill \break $v)$ $A'$ est engendr\'ee \hbox{\rm{comme $A$-alg\`ebre}} par $\langle M,
\omega_{\scriptstyle A'}(\check{\cal M})\rangle$ et $\langle \check{M},\omega_{\scriptstyle
A'}({\cal M})\rangle$. \par
\medskip
\proclaim D\'efinition 4.1.2. Supposons que $A$ soit semi-simple (i.e. produit fini de corps). On dit que
$(A',d')$ est une extension \hbox{\rm{de Picard-Vessiot fractionnaire}} de $(A,d)$ si 
\hfill \break $i)$ $A'$ est semi-simple,
\hfill \break $ii)$ $(A',d')$ est simple,
\hfill \break $iii)$ l'anneau des constantes de $(A',d')$ est le corps $\;k$,
\hfill \break $iv)$ ${\cal M}$ est soluble dans $(A',d')$,
\hfill \break $v)$ $A'$ est \hbox{\rm une localisation} de la $A$-alg\`ebre engendr\'ee par
$\langle M, \omega_{\scriptstyle A'}(\check{\cal M})\rangle$ et $\langle \check{M},\omega_{\scriptstyle
A'}({\cal M})\rangle$. 
\par \medskip \noindent {\bf 4.1.3. Remarques.} $a)$ Il semble impossible de trouver une terminologie
compatible \`a toutes celles de la litt\'erature. Nous avons privil\'egi\'e, comme dans
[vdPS97] les extensions de type fini par rapport \`a leurs contreparties ``birationnelles"; les
extensions de la th\'eorie de Picard-Vessiot classique (cf. I.3.1) sont, dans la terminologie
ci-dessus, des extensions de Picard-Vessiot fractionnaires avec $A=\Omega^1$ et $A'=\Omega'^1$.
Nous verrons ult\'erieurement que toute extension de Picard-Vessiot
fractionnaire appara\^{\i}t comme anneau total de fractions d'une extension de Picard-Vessiot enti\`ere. 
\par \noindent $b)$ Toute extension de Picard-Vessiot est de type fini en
tant qu'extension d'anneaux.    
\par \noindent $c)$ Toute extension de Picard-Vessiot fractionnaire $A'/A$ est fid\`element plate (en
tant que module sur le produit fini de corps $\Pi K_i$, $A'\cong \Pi V_i$, o\`u $V_i$ est un espace
vectoriel non nul sur $K_i$).
\par \noindent $d)$ Supposons qu'il existe une extension de Picard-Vessiot fractionnaire $A'/A$ pour
${\cal M}$. Soit $J/A$ une extension diff\'erentielle interm\'ediaire, avec $J$ semi-simple. Alors
$A"/J$ est une extension de Picard-Vessiot fractionnaire pour
${\cal M}_{J}$ (c'est clair).

\proclaim Lemme 4.1.4. S'il existe une extension de Picard-Vessiot $A'/A$ pour $\cal
M$ (resp. de Picard-Vessiot fractionnaire), alors
$\;{\scriptstyle <}{\cal M}{\scriptstyle >}^{\scriptstyle\otimes}\;$ est une cat\'egorie
tannakienne neutre sur $k$: tout objet est rigide et le foncteur ``solutions dans $A'\;$" est un
foncteur fibre (i.e. un foncteur mono\"{\i}dal
$k$-lin\'eaire fid\`ele et exact) \`a valeurs dans les $k$-espaces vectoriels de dimension
finie.
\par Cela r\'esulte de 1.3.2 (resp. et de 4.1.3.$c)$ ).

\medskip {\bf 4.2. Foncteurs fibres, torseurs de solutions et extensions de Picard-Vessiot.}
\par \noindent Ce paragraphe est inspir\'e de [De90]9. On suppose que
${\scriptstyle <}{\cal M}{\scriptstyle >}^{\scriptstyle\otimes}$ est une cat\'egorie tannakienne
neutre sur $k$. On fixe un foncteur fibre $\omega$. Ce foncteur fibre rend ${\scriptstyle <}{\cal
M}{\scriptstyle >}^{\scriptstyle\otimes}$ \'equivalente \`a la cat\'egorie mono\"{\i}dale des
repr\'esentations de dimension finie de $Gal({\cal M},\omega)$.
\par \noindent Pour tout sous-objet $\cal N$ d'un $\oplus \; {\cal M}^{\otimes
i}\otimes {\check{\cal M}}^{\otimes j}$, notons $\;{\scriptstyle <}{\cal N}{\scriptstyle
>}^{\scriptstyle\oplus}\;$ la sous-cat\'egorie ab\'elienne strictement pleine de 
$\;{\scriptstyle <}{\cal M}{\scriptstyle >}^{\scriptstyle\otimes}$ form\'ee des sous-quotients
des sommes directes finies de copies de $\cal N$. Ordonnons les objets $\cal N$ par la relation
``\^etre facteur direct".
\par \noindent Soit $\;{\cal O}(Gal({\cal M},\omega))\;$ la $k$-alg\`ebre des fonctions sur
$Gal({\cal M},\omega)$, vue comme repr\'esentation de $Gal({\cal M},\omega)$ 
(repr\'esentation r\'eguli\`ere gauche: $f(x)\mapsto f(g^{-1}x)$). On a la formule 
\par \centerline{${\cal O}(Gal({\cal M},\omega))={\displaystyle\limind_{\scriptstyle\cal
N}\;}({\rm End}(\omega
\mid {\scriptstyle <}{\cal N}{\scriptstyle >}^{\scriptstyle\oplus}))^{\vee}$}
\par \noindent o\`u ${}^{\vee}$ d\'esigne un $k$-dual (cf. [Saa72]II.3.2.6, [DeMi82]2.14).
\par \noindent Observons que ${\rm End}(\omega \mid {\scriptstyle <}{\cal
N}{\scriptstyle >}^{\scriptstyle\oplus})$ est un sous-$k$-espace (et m\^eme une
sous-$k$-alg\`ebre) de ${\rm End}(\omega({\cal N}))$. Si l'on fait agir $Gal({\cal M},\omega)$
par composition \`a gauche sur ${\rm End}(\omega({\cal N}))$, ce sous-espace est stable. 
En consid\`erant ${\rm End}(\omega \mid {\scriptstyle <}{\cal
N}{\scriptstyle >}^{\scriptstyle\oplus})^{\vee}$ comme $k$-repr\'esentation quotient de $({\rm
End}(\omega({\cal N})))^{\vee}$, on voit que la formule ${\cal O}(Gal({\cal
M},\omega)) \break ={\displaystyle\limind\;}({\rm End}(\omega
\mid {\scriptstyle <}{\cal N}{\scriptstyle >}^{\scriptstyle\oplus}))^{\vee}$ est compatible
\`a l'action de $Gal({\cal M},\omega)$. 
\medskip Consid\'erons d'autre part la $A$-alg\`ebre ${\cal O}(\Sigma({\cal M},\omega))$ des
fonctions sur le torseur $\Sigma({\cal M},\omega)$ (cf. 2.2.1). On a la formule 
\par \centerline{${\cal O}(\Sigma({\cal
M},\omega))={\displaystyle\limind_{\scriptstyle\cal
N}\;}({\rm Hom}(\omega \otimes 1_{\scriptstyle A}, \hbox{\smit{oubli}} \mid
{\scriptstyle <}{\cal N}{\scriptstyle >}^{\scriptstyle\oplus}))^{\vee}$}
\par \noindent o\`u ${}^{\vee}$ d\'esigne ici un $A$-dual. 
\par \noindent Observons que  $\;{\rm Hom}(\omega \otimes 1_{\scriptstyle A}, \hbox{\smit{oubli}}
\mid {\scriptstyle <}{\cal N}{\scriptstyle >}^{\scriptstyle\oplus}))\;$ est un
sous-$A$-module de ${\cal I}hom(A\otimes_k \omega({\cal N}),{\cal N})$. La connexion sur ce
dernier (o\`u $A\otimes_k \omega({\cal N})$ est muni de la connexion triviale) laisse stable $\;{\rm Hom}(\omega \otimes 1_{\scriptstyle A}, \hbox{\smit{oubli}}
\mid {\scriptstyle <}{\cal N}{\scriptstyle >}^{\scriptstyle\oplus}))$. Ainsi, la connexion
sur ${\cal I}hom((A\otimes_k \omega({\cal N}),{\cal N})^{\vee} \cong \check{\cal N}\otimes_k
\omega({\cal N})$ passe au quotient $({\rm Hom}(\omega \otimes 1_{\scriptstyle A},
\hbox{\smit{oubli}} \mid {\scriptstyle <}{\cal N}{\scriptstyle >}^{\scriptstyle\oplus}))^{\vee}$
et induit une connexion $\;d'\;$ sur ${\cal O}(\Sigma({\cal M},\omega))$. 
\par \noindent On peut donc
consid\'erer $\;({\cal O}(\Sigma({\cal M},\omega)),d')\;$ comme une extension diff\'erentielle de
$\;(A,d)$.    
\proclaim Lemme 4.2.1. $i)$ Le prolongement de $\;\omega \;$ \`a $\;Ind{\scriptstyle <}{\cal
M}{\scriptstyle >}^{\scriptstyle\otimes}\;$ envoie $\;({\cal O}(\Sigma({\cal
M},\omega)),d')\;$ sur $\;{\cal O}(Gal({\cal M},\omega))\;$ muni de la repr\'esentation
r\'eguli\`ere gauche. 
\hfill\break $ii)$ $({\cal O}(\Sigma({\cal M},\omega)),d')$ est une extension de Picard-Vessiot
de $(A,d)$ pour $\cal M$. 
\par {\it Preuve.} $i)$ $\;\omega \;$ envoie ${\cal I}hom(A\otimes_k \omega({\cal
N}),{\cal N})^{\vee}$ est $({\rm End}(\omega({\cal N})))^{\vee}$ muni de l'action \`a gauche de
$Gal({\cal M},\omega)$, et l'objet quotient $({\rm Hom}(\omega \otimes 1_{\scriptstyle A},
\hbox{\smit{oubli}} \mid {\scriptstyle <}{\cal N}{\scriptstyle >}^{\scriptstyle\oplus}))^{\vee}$
sur la repr\'esentation quotient $({\rm End}(\omega \mid {\scriptstyle <}{\cal N}{\scriptstyle
>}^{\scriptstyle\oplus}))^{\vee}$, d'o\`u le r\'esultat en passant \`a la limite. 
\par \noindent $ii)$ $A'={\cal O}(\Sigma({\cal M},\omega))$ est fid\`element plate sur $A$, cf.
2.2.1. Tout id\'eal diff\'erentiel d\'efinissant un sous-objet de $(A',d')$ dans
$Ind{\scriptstyle <}{\cal M}{\scriptstyle >}^{\scriptstyle\otimes}$, la simplicit\'e de $(A',d')$
\'equivaut au fait que les seuls sous-sch\'emas ferm\'es de
$G=Gal({\cal M},\omega)$ invariants par translation \`a gauche sont $\emptyset$ et $G$ (c'est
ici qu'intervient le fait que l'on travaille avec ${\scriptstyle <}{\cal M}{\scriptstyle
>}^{\scriptstyle\otimes}$ et non une quelconque sous-cat\'egorie ab\'elienne mono\"{\i}dale
pleine $\cal C$ comme en 2.1.1). Les constantes de $(A',d')$ s'identifient aux \'elements de
$\;{\cal O}(G)\;$ invariants par translation, c'est-\`a-dire aux \'el\'ements de $k$. 
\par \noindent Montrons que tout objet $\cal N$ de ${\scriptstyle <}{\cal
M}{\scriptstyle >}^{\scriptstyle\otimes}$ est soluble dans $A'$:  $A'\otimes_{k}
({\cal N}_{A'})^{\nabla_{A'}} {\buildrel{\simeq}\over {\rightarrow}} {\cal N}_{A'}$. Via $\omega$,
il s'agit d'\'etablir que la fl\`eche naturelle de $G$-comodules 
\par \centerline{${\cal O}(G)\otimes_k ({\cal O}(G)\otimes_k \omega({\cal N}))^G \rightarrow
{\cal O}(G)\otimes_k \omega({\cal N}) $}
\par \noindent est un isomorphisme. Si $E$ est le fibr\'e vectoriel trivial de fibre $\omega({\cal
N})$, cela traduit le fait que le morphisme $G$-\'equivariant 
\par \centerline{$G \times E \rightarrow G \times ((G \times E) /G), \;\;\;(g,e)\mapsto
(g,(g,e))$}
\par \noindent est un isomorphisme. Cet argument fournit de m\^eme un
isomorphisme naturel        
\par \centerline{$\omega({\cal N}) \cong \omega_{\scriptstyle
A'}({\cal N})=({\cal N}_{A'})^{\nabla_{A'}}$}
\par \noindent fonctoriel en $\cal N$.
\par \noindent Enfin, soit $A"$ la sous-$A$-alg\`ebre de $A'$ engendr\'ee par
$\langle M,
\omega_{\scriptstyle A'}(\check{\cal M})\rangle$ et $\langle
\check{M},\omega_{\scriptstyle 
A'}({\cal M})\rangle$. C'est aussi la r\'eunion des $\langle N, \omega_{\scriptstyle
A'}(\check{\cal N})\rangle$ pour tout objet rigide $\cal N$ de ${\scriptstyle <}{\cal
M}{\scriptstyle >}^{\scriptstyle\otimes}$. Il est clair que $A"$ est stable sous tout $D\in
{}^{\vee}\Omega^1$, donc fournit une extension diff\'erentielle interm\'ediaire entre $A$ et
$A'$. De plus, comme tout objet de ${\scriptstyle <}{\cal
M}{\scriptstyle >}^{\scriptstyle\otimes}$ est soluble dans $A"$, on a un $A"$-point
canonique de $\Sigma({\cal M},\omega)$, i.e. un homomorphisme d'alg\`ebres $A'\rightarrow A"$,
compatible aux connexions, et dont l'inclusion $A" \subset A'$ est une section. Comme $(A',d')$
est simple, on a donc $A'=A"$, ce qui ach\`eve la d\'emonstration.    
\medskip \noindent {\bf 4.2.2. Remarque.} Ce lemme traduit alg\'ebriquement l'id\'ee
g\'eom\'etrique suivante. La consid\'eration du fibr\'e principal \`a droite $P=\Sigma({\cal M},\omega)$,
de groupe $G=Gal({\cal M},\omega)$, de base $X=Spec\; A$, \'etablit une analogie avec la situation
I.1.2, I.1.3. Pr\'ecisons en nous pla\c cant dans cette situation: on a un fibr\'e \`a connexion $E$ sur
$X$, associ\'e \`a un $GL_n$-fibr\'e principal \`a connexion sur $X$. Ce dernier se r\'eduit \`a un
fibr\'e principal $P{\buildrel{\pi}\over\rightarrow} X$ sous le groupe d'holonomie $G$ muni d'une
connexion $\aleph$. 
\par \noindent L'analogue de l'extension de Picard-Vessiot (${\cal O}(\Sigma({\cal M},\omega)),d')$ est
l'anneau $C^{\infty}(P)$ des fonctions de $P$, muni de la {\it diff\'erentielle longitudinale} $\;d':
C^{\infty}(P)\rightarrow \Gamma(T^{\vee}(X))$ (en consid\'erant $\aleph$ comme un
``feuilletage non-int\'egrable" $\cal F$ de $P$). La connexion ``le long des feuilles" sur l'image inverse
de $E$ sur $P$: $\Gamma(E_P)\rightarrow \Gamma (T^{\vee}({\cal F})\otimes E_P)$ est alors triviale. Ceci
exprime le fait suivant. Consid\'erons l'image inverse du fibr\'e $\pi : P\rightarrow X$ sur $P$
lui-m\^eme. L'application $\theta:P\times G \rightarrow P\times_X P,\;(p,g)\mapsto (p,t_g(p))\;$ est un
isomorphisme $G$-\'equivariant (action \`a droite sur les facteurs de droite). On a:
$\;(\theta_{\ast})_{p,t_g(p)}(H_p(P),0_{\hbox{\goth g}}) = Im((id,t_g)_{\ast} H_p(P)) \subset
T_{p,t_g(p)}(P\times_X P)$, qui s'envoie isomorphiquement sur $H_{t_g(p)}(P)$ par la seconde projection. 
\proclaim Th\'eor\`eme 4.2.3. Sous l'hypoth\`ese que ${\scriptstyle <}{\cal M}{\scriptstyle
>}^{\scriptstyle\otimes}$ soit une cat\'egorie tannakienne neutre sur $k$, on a des \'equivalences
de cat\'egories quasi-inverses
\medskip \centerline{$\Bigl\{ \matrix{ \hbox{foncteurs fibres} \cr \hbox{sur}\;{\scriptstyle
<}{\cal M}{\scriptstyle >}^{\scriptstyle\otimes}}\Bigr\} 
\longleftrightarrow \Bigl\{\matrix{ \hbox{extensions de} \cr
\hbox{Picard-Vessiot pour}\;{\cal M}}\Bigr\} $}
\medskip \noindent donn\'ees par $\;\;\omega \mapsto ({\cal O}(\Sigma({\cal
M},\omega)),d'),\;\;\;\; (A',d') \mapsto \omega_{A'}$.
\hfill \break En outre, si $\omega$ est un quelconque foncteur fibre, ces cat\'egories sont
\'equivalentes \`a la cat\'egorie des $Gal({\cal M},\omega)$-torseurs. 
\par {\it Preuve.} On remarque que les morphismes de ces cat\'egories sont des isomorphismes
(du c\^ot\'e extensions de Picard-Vessiot, cela vient de leur simplicit\'e). On voit ainsi que    
$\;\omega \mapsto ({\cal O}(\Sigma({\cal
M},\omega)),d'),\;$ et $\; (A',d') \mapsto \omega_{A'}$ d\'efinissent bien des foncteurs. Qu'ils
soient quasi-inverses r\'esulte de l'isomorphisme $\omega \cong \omega_{\scriptstyle
A'}$ observ\'e ci-dessus, avec $A'={\cal O}(\Sigma({\cal M},\omega))$. La seconde assertion est
une g\'en\'eralit\'e dans les cat\'egories tannakiennes neutres, cf.[Saa72]II.3.2.3.3.  
\proclaim Corollaire 4.2.4. Il existe une extension de Picard-Vessiot de $(A,d)$
pour ${\cal M}$ si et seulement si ${\scriptstyle <}{\cal M}{\scriptstyle >}^{\scriptstyle\otimes}$ est une cat\'egorie tannakienne
neutre sur $k$. \par 
Cela r\'esulte de 4.1.4 et 4.2.3.
\proclaim Corollaire 4.2.5. Soit $\bar k$ une cl\^oture alg\'ebrique de
$k$. Le $\bar k$-groupe alg\'ebrique $Gal({\cal M},\omega)\otimes {\bar k}$ ne d\'epend pas, \`a
isomorphisme pr\`es, du foncteur fibre $\;\omega$ (suppos\'e exister). Les classes d'isomorphie
d'extensions de Picard-Vessiot sont en bijection avec $H^1({\bar k}/k, Gal({\cal
M},\omega)\otimes {\bar k})$. Deux extensions de Picard-Vessiot pour
$\cal M$ quelconques deviennent isomorphes apr\`es extension finie du corps des constantes $k$.  
\par \noindent {\bf 4.2.6. Exemple.} Pour la
connexion correspondant \`a l'\'equation diff\'erentielle ${d \over dz}\;y={1\over 2z}\;y$ sur
${\bf R}(z)$, on a deux extensions de Picard-Vessiot non-isomorphes: ${\bf R}(z,{\root\of{z}})$
et ${\bf R}(z,{\root\of{-z}})$; elles donnent lieu \`a des torseurs de
solutions non-triviaux.
 
\medskip {\bf 4.3. Existence d'extensions de Picard-Vessiot.} 
\proclaim Th\'eor\`eme 4.3.1. On suppose que 
\hfill \break $i)$ l'anneau commutatif $\;A\;$ est noeth\'erien, et son anneau total de fractions 
$Q(A)$ est semi-simple, 
\hfill \break $ii)$ $\Omega^1=dA.A$ est fid\`ele et projectif de type fini \`a droite, et $\Omega^1
\otimes_{\scriptstyle A}Q(A) \cong Q(A)\otimes_{\scriptstyle A}
\Omega^1$, 
\hfill \break $iii)$ l'anneau diff\'erentiel $(A,d)$ est simple de corps de constantes $\;k$. 
\medskip \noindent Alors quitte \`a remplacer $\;k\;$ par une extension finie, il
existe une extension de Picard-Vessiot de $(A,d)$ pour ${\cal M}$. 
\hfill \break En particulier, si $k$ est alg\'ebriquement clos, il existe une extension de Picard-Vessiot
de $(A,d)$ pour ${\cal M}$, qui est unique \`a isomorphisme non unique pr\`es. 
\par {\it Preuve.} Il suit de II.5.3.2 que ${\scriptstyle <}{\cal M}{\scriptstyle
>}^{\scriptstyle\otimes}$ est tannakienne. Comme elle admet un g\'en\'erateur tensoriel (par exemple ${\cal
M}\oplus \check{\cal M}$), il existe une extension finie $k'/k'$ et un foncteur fibre $\omega'$ \`a
valeurs dans les $k'$-espaces vectoriels de dimension finie ([De90]6.17,[Saa72]III.3.3). Ce foncteur fibre fait de
${\scriptstyle <}{\cal M}{\scriptstyle >}^{\scriptstyle\otimes}_{(k')}$ une cat\'egorie tannakienne neutre
([DeMi82]3.11). D'autre part, un $A\otimes k'$-module
\`a connexion (relativement \`a l'anneau diff\'erentiel $(A,d)\otimes_{\scriptstyle k} k'\;$) n'est rien
d'autre qu'un $A$-module \`a connexion $\cal N$ muni de $k'\rightarrow \;$End$(\cal N)$,
ce qui, du fait que $k'$ est fini sur $k$, fournit une \'equivalence tautologique de cat\'egories
mono\"{\i}dales $k'$-lin\'eaires ${\scriptstyle <}{\cal M}{\scriptstyle >}^{\scriptstyle\otimes}_{(k')}\cong
{\scriptstyle <}{\cal M}_{k'}{\scriptstyle >}^{\scriptstyle\otimes}$. D'apr\`es 4.2.4, il existe donc une
extension de Picard-Vessiot de $(A,d)\otimes_{\scriptstyle k} k'$ pour ${\cal M}_{k'}$. L'unicit\'e
lorsque $k=k'=\bar k$ suit de 4.2.5. 
\medskip \noindent {\bf 4.3.2.} Remarquons encore qu'il n'a \'et\'e fait aucune hypoth\`ese sur la
courbure de $\nabla$. En termes figur\'es, 4.3.1 est un th\'eor\`eme d'int\'egrabilit\'e symbolique des
syst\`emes diff\'erentiels (ou aux diff\'erences) non n\'ecessairement int\'egrables. Pour un exemple
concret, voir 1.3.3; pour une interpr\'etation g\'eom\'etrique, voir 4.2.2.  
\medskip \noindent {\bf 4.3.3.} La condition suppl\'ementaire que le corps des constantes $k$ est alg\'ebriquement clos, 
qui garantit l'unicit\'e de l'extension de Picard-Vessiot, est innocente en caract\'eristique $0$. Elle n'est en revanche jamais
remplie en caract\'eristique $p>0$ sous les hypoth\`eses g\'en\'erales du th\'eor\`eme: en effet, elle entra\^{\i}ne
$Q(A)^p\subset k$, et comme $Q(A)$ est suppos\'e semi-simple, on voit $Q(A)=k$;
 mais alors $\Omega^1=AdA=0$ n'est pas fid\`ele.

\medskip {\bf 4.4. Quelques propri\'et\'es des extensions de Picard-Vessiot (fractionnaires).}
\proclaim Proposition 4.4.1. Soit $(A',d')$ une extension de Picard-Vessiot pour
$\cal M$. Soit $\cal N$ un objet de ${\scriptstyle
<}{\cal M}{\scriptstyle >}^{\scriptstyle\otimes}$. Alors la sous-$A$-alg\`ebre $A"$ de
$A'$ engendr\'ee par $\langle N,
\omega_{\scriptstyle A'}(\check{\cal N})\rangle$ et $\langle
\check{N},\omega_{\scriptstyle
A'}({\cal N})\rangle$ fournit une extension diff\'erentielle interm\'ediaire entre $A$ et
$A'$. C'est une extension de Picard-Vessiot pour $\cal N$.  
\par {\it Preuve.} D'apr\`es 4.2.3, on peut supposer $(A',d')= ({\cal O}(\Sigma({\cal
M},\omega_{A"})),d')$. Consid\'erons l'homomorphisme fid\`element plat
$\phi:\;Gal({\cal N},\omega_{\scriptstyle A'})\rightarrow Gal({\cal M},\omega_{\scriptstyle A'})$
(2.3.2). Il induit un $\phi$-morphisme fid\`element plat $\phi'$ du $Gal({\cal
M},\omega_{\scriptstyle A'})_A$-torseur $\Sigma({\cal M},\omega_{\scriptstyle A'})$ vers le
$Gal({\cal N},\omega_{\scriptstyle A'})_A$-torseur $\Sigma({\cal N},\omega_{\scriptstyle A'})$,
d'o\`u une inclusion d'alg\`ebres compatible aux connexions ${\cal O}(\Sigma({\cal
N},\omega_{\scriptstyle A'}))\hookrightarrow{\cal O}(\Sigma({\cal M},\omega_{\scriptstyle A'}))$. On
conclut en appliquant 4.2.1.$ii)$ \`a $({\cal O}(\Sigma({\cal N},\omega_{A"})),d')$.    

\proclaim Proposition 4.4.2. Supposons que $A$ soit semi-simple et qu'il existe une extension de
Picard-Vessiot fractionnaire $(A",d")$ de $(A,d)$ pour ${\cal M}$. Alors il existe une unique extension
diff\'erentielle interm\'ediaire $(A',d')$ qui est de Picard-Vessiot pour ${\cal M}$; on a $A"=Q(A')$.
\par {\it Preuve.} $A"$ est l'anneau total de fractions de la
sous-$A$-alg\`ebre $A'$ de $A"$ engendr\'ee par $\langle M, \omega_{A"}(\check{\cal M})\rangle$ et $\langle
\check{M},\omega_{A"}({\cal M})\rangle$. D'apr\`es 4.2.3, il existe une extension de Picard-Vessiot
${\tilde A}/A$ (qu'on peut prendre \'egale \`a ${\cal O}(\Sigma({\cal M},\omega_{A"})$ telle que
$\omega_{A"}$ s'identifie \`a $\omega_{\tilde A}$. L'homomorphisme canonique
$({\cal O}(\Sigma({\cal M},\omega_{A"})),d')={\tilde A}\rightarrow A"$ est compatible aux
diff\'erentielles, donc injectif puisque $({\tilde A},{\tilde d})$ est simple. On en d\'eduit que
$({\tilde A},{\tilde d})$ s'identifie \`a $(A',d')$. 
\medskip Prendre garde \`a la r\'eciproque: si $A$ est semi-simple, et si $(A',d')$ est une extension de
Picard-Vessiot avec $Q(A')$ semi-simple, il n'est pas clair en g\'en\'eral que
$d'(Q(A')).Q(A')\cong d'A'.Q(A')$, et donc que la structure naturelle d'anneau
diff\'erentiel sur $Q(A')$ en fasse une extension diff\'erentielle de $(A,d)$ au sens de II.1.2.5 (voir
aussi II.1.3.5). C'est toutefois le cas si $d'A'.A'$ est un sesquimodule.   
\proclaim Corollaire 4.4.3. S'il existe une extension de
Picard-Vessiot fractionnaire $(A",d")$ pour ${\cal M}$, alors $\Sigma({\cal M},\omega_{A"})$ est lisse,
et $A"$ est le produit des corps de fonctions $K_i$ de ses composantes connexes. 
\par En effet, on a vu que l'anneau semi-simple $A"$ s'identifie \`a l'anneau total de fractions de ${\cal
O}(\Sigma({\cal M},\omega_{A"})$, donc $\Sigma({\cal M},\omega_{A"})$ est r\'eduit. Etant un torseur sur
un anneau semi-simple, il est lisse.  

\proclaim Proposition 4.4.4. Si $A$ est int\`egre de caract\'eristique nulle, toute extension de
Picard-Vessiot $A'$ est lisse. Elle est int\`egre si $\;\Omega^1$ est un bimodule commutatif. En
particulier, si $A$ est un corps, si le corps de constantes est alg\'ebriquement clos de
caract\'eristique nulle, et si $\;\Omega^1$ est un bimodule commutatif, il existe un corps extension
de Picard-Vessiot fractionnaire pour tout $A$-espace vectoriel de dimension finie \`a connexion (non
n\'ecessairement int\'egrable). 
\par {\it Preuve.} En effet, si $A$ est int\`egre de caract\'eristique nulle, tout $A$-groupe
alg\'ebrique, et tout torseur sous un tel groupe, est lisse. Or on a vu que $A'\cong {\cal
O}(\Sigma({\cal M},\omega_{A'}))$. Si en outre $\;\Omega^1$ est un bimodule commutatif, il en
est de m\^eme de $\;\Omega'^1$; tout idempotent de $A'$ est n\'ecessairement une constante,
\'egale \`a  $0$ ou $1$, d'o\`u l'assertion. Si $A$ est un corps, on obtient une extension
de Picard-Vessiot fractionnaire en passant au corps de fractions de $A'$.
\medskip \noindent Dans le cas d'une connexion attach\'ee \`a une \'equation aux diff\'erences
sur un corps de caract\'eristique nulle, les extensions de Picard-Vessiot sont lisses mais non
n\'ecessairement int\`egres (consid\'erer l'\'equation $y^{\sigma}=-y$).

\proclaim Proposition 4.4.5. Soit $(A',d')$ un corps extension de Picard-Vessiot fractionnaire de
$(A,d)$ pour ${\cal M}$. Alors le degr\'e de transcendance de $A'$ sur $A$ est \'egal \`a la dimension
de $Gal({\cal M},\omega_{A'})$. 
\par {\it Preuve.} D'apr\`es 4.4.2 et 4.2.3, $A'$ est isomorphe au corps de fonctions de la
$A$-vari\'et\'e $\Sigma({\cal M},\omega_{A'})$ qui est un torseur sous $Gal({\cal
M},\omega_{A'})\otimes_k A$. D'o\`u le r\'esultat.   

\bigskip {\bf 5. Correspondance galoisienne.}
\medskip {\bf 5.1. Le groupe alg\'ebrique affine ${\bf Aut}((A',d')/(A,d))$.}
\medskip \noindent Soit $(A',d')$ une extension diff\'erentielle de $(A,d)$. On note 
$\u{Aut}((A',d')/(A,d))$ le foncteur qui associe \`a toute $k$-alg\`ebre $k'$ (commutative
unitaire) le groupe des automorphismes de l'anneau diff\'erentiel $(A',d')\otimes_k k'$ induisant
l'identit\'e sur $A\otimes_k k'$ (et donc aussi sur $\Omega^1\otimes_k k'$). 
\proclaim Th\'eor\`eme 5.1.1. Soit $(A',d')$ est une extension de Picard-Vessiot ou une extension
de Picard-Vessiot fractionnaire pour $\cal M$. Alors le foncteur $\u{Aut}((A',d')/(A,d))$ est
repr\'esent\'e par $Gal({\cal M},\omega_{A'})$.
\par Dans ce r\^ole, ce $k$-groupe alg\'ebrique affine est aussi not\'e ${\bf Aut} ((A',d')/(A,d))$ ou
simplement ${\bf Aut}_d(A'/A)$, et appel\'e {\it groupe de Galois diff\'erentiel de $A'$ sur $A$}.      
\medskip {\it Preuve de 5.1.1.} Traitons d'abord le cas d'une extension de Picard-Vessiot.
D'apr\`es 4.2.3, on peut supposer $(A',d')= ({\cal O}(\Sigma({\cal M},\omega_{A'})),d')$. Via
$\omega_{A'}$, le foncteur $\;\u{Aut}((A',d')/(A,d))\;$ s'identifie au foncteur
$\;\u{Aut}_{\scriptstyle Gal({\cal M},\omega_{A'})}({\cal O}(Gal({\cal M},\omega_{A'})))$, qui
n'est autre que le foncteur que repr\'esente $Gal({\cal M},\omega_{A'})$. Le cas d'une extension
de Picard-Vessiot fractionnaire s'en d\'eduit gr\^ace \`a 4.4.2 (Dans la situation de 4.4.2, il est clair
que le groupe de Galois diff\'erentiel de l'extension de Picard-Vessiot fractionnaire $A"/A$ s'identifie
au groupe de Galois diff\'erentiel de l'extension de Picard-Vessiot associ\'ee $A'/A$).

\medskip {\bf 5.2. La correspondance galoisienne pour les extensions de Picard-Vessiot fractionnaires.} 
\medskip \noindent Dans la suite, on se limite au cas o\`u le corps des constantes $k$ est
alg\'ebriquement clos de caract\'eristique nulle. La donn\'ee du sch\'ema en groupes ${\bf
Aut}_d(A'/A)$ \'equivaut alors \`a celle du groupe de ses $k$-points vu comme sous-groupe de
$GL(\omega_{A'}({\cal M}))$; c'est le groupe $G=Aut_d(A'/A)$ des {\it $A$-automorphismes de $A'$ commutant
\`a $d$} (en revanche, il n'y a pas d'action alg\'ebrique de ${\bf Aut}_d(A'/A)$ sur $A'$ en
g\'en\'eral).   

\proclaim Lemme 5.2.1. Supposons que $A=K=\Pi K_i$ soit semi-simple et qu'il existe une extension de
Picard-Vessiot fractionnaire $L$ pour $\cal M$. Soit $H$ un sous-groupe Zariski-ferm\'e de $G=Aut_d(L/K)$.
Alors 
\hfill \break $i)$ l'anneau des invariants $L^G$ est $K$,
\hfill \break $ii)$ si $L^H=K$, alors $H=G$.
\par {\it Preuve.} D'apr\`es 4.4.2, $L=Q(K')$ pour une extension de Picard-Vessiot ``enti\`ere"
$K'$ de $K$. Posons ${\bf G}=Gal({\cal M},\omega_{K'}) = {\bf Aut}_d(L/K)$ (5.1.1), et notons ${\bf
H}\subset {\bf G}$ le sous-sch\'ema en groupes de groupe de $k$-points $H$. On peut identifier $K'$ \`a
${\cal O}(\Sigma({\cal M},\omega_{K'}))$ (4.2.3), donc tout \'el\'ement de $L$ \`a un $K$-morphisme
$f:\Sigma({\cal M},\omega_{K'})\rightarrow {\bf P}^1_K$ (i.e. \`a une famille de $K_i$-morphismes
$\Sigma({\cal M},\omega_{K'})\otimes_K K_i\rightarrow {\bf P}^1_{K_i}$). D'apr\`es
2.2.1, $\Sigma({\cal M},\omega_{K'})$ est un torseur \`a droite sous ${\bf G}_K={\bf G}\otimes_k K$. Il est
\'equivalent de dire qu'un \'element de $L$ est invariant sous $G$ (resp. sous $H$) ou que le morphisme
correspondant $f$ est invariant sous ${\bf G}_K$ (resp. sous ${\bf H}_K$). Comme ${\bf G}_K$ agit
transitivement sur $\Sigma({\cal M},\omega_{K'})$, un tel $f$ invariant sous ${\bf G}_K$ est constant, ce
qui implique $i)$.  Pour $ii)$, il s'agit de montrer que si $H\neq G$, il existe un fonction rationnelle
${\bf H}_K$-invariante non constante. C'est une cons\'equence de l'existence du $K$-sch\'ema $\Sigma({\cal
M},\omega_{K'})/{\bf H}_K$ (qui se d\'eduit par descente \'etale de l'existence de la $k$-vari\'et\'e
alg\'ebrique ${\bf G}/{\bf H}$). 
\medskip Lorsque $G$ et $H$ sont r\'eductifs, ${\bf G}/{\bf H}$ est affine, de m\^eme que $\Sigma({\cal
M},\omega_{K'})/{\bf H}_K$, et on peut alors remplacer ${\bf P}^1$ par ${\bf A}^1$, et extensions de
Picard-Vessiot fractionnaires par extensions de Picard-Vessiot ``enti\`eres".
\proclaim Th\'eor\`eme 5.2.2. Supposons que $A=K$ soit semi-simple et qu'il existe une extension de
Picard-Vessiot fractionnaire $L$ pour $\cal M$. Alors $H\mapsto (L^H,d)$ et $(J,d)\mapsto Aut_d(L/J)$ sont
des bijections r\'eciproques entre l'ensemble des sous-groupes Zariski-ferm\'es de $G=Aut_d(L/K)$ et
l'ensemble des extensions diff\'erentielles interm\'ediaires $(J,d)$ entre $K$ et
$L$, avec $J$ semi-simple. En outre, si $H$ est normal dans $G$, $L^H$ est une extension de
Picard-Vessiot fractionnaire $L$ pour un objet convenable $\cal N$ de 
${\scriptstyle <}{\cal M}{\scriptstyle >}^{\scriptstyle \otimes}$. 
\par {\it Preuve.} $L/J$ \'etant une extension de Picard-Vessiot fractionnaire (4.1.3.$d$), il
d\'ecoule de 5.1.1 que $Aut_d(L/J)$ est un sous-groupe Zariski-ferm\'e de
$G=Aut_d(L/K)$. Montrons, inversement, que $L^H$ d\'efinit une extension diff\'erentielle
interm\'ediaire. Il s'agit de voir que pour tout $D\in {}^{\vee}\Omega^1, \;\langle D , d(L^H)\rangle
\subset L^H$. Cela vient de ce que l'action de $G$ sur $L$ commute \`a $d$ donc \`a l'endomorphisme
$k$-lin\'eaire $x\mapsto \langle D , dx\rangle $ de $L$. Qu'on ait des bijections r\'eciproques
r\'esulte alors du lemme pr\'ec\'edent.
\par \noindent Enfin, si $H$ est normal dans $G$, ${\bf H}$ est normal dans ${\bf G}=Gal({\cal
M},\omega_L)$, et le groupe quotient ${\bf G}/{\bf H}$ est $k$-groupe associ\'e \`a
une sous-cat\'egorie tannakienne de ${\scriptstyle <}{\cal M}{\scriptstyle >}^{\scriptstyle \otimes}$.
Une telle sous-cat\'egorie tannakienne admet n\'ecessairement un g\'en\'erateur tensoriel $\cal N$, et on
a alors ${\bf G}/{\bf H}=Gal({\cal N},\omega_L)$. Comme ${\bf H}=Gal({\cal
M}_J,\omega_L)$ agit trivialement sur $\omega_L({\cal N})$, ${\cal N}$ est soluble dans $J$. La
$K$-alg\`ebre engendr\'ee par $\langle M, \omega_{J}(\check{\cal M})\rangle$ et $\langle
\check{M},\omega_{J}({\cal M})\rangle$ fournit donc une extension de Picard-Vessiot $K'$ de $K$ pour
${\cal N}$ (4.4.1). Son anneau total de fractions $Q(K')$ est contenu dans l'anneau semi-simple. Il est
donc semi-simple (et diff\'erentiellement simple d'apr\`es II.1.3.5). On a ${\bf G}/{\bf H}=Gal({\cal
N},\omega_L)= {\bf Aut}_d(Q(K')/K)$ (5.1.1), et on conclut de la correspondance galoisienne que $Q(K')=J$.
Ceci montre que $J/K$ est une extension de Picard-Vessiot fractionnaire pour $\cal N$.

\medskip Ce r\'esultat contient et unifie la correspondance de Picard-Vessiot classique (cf. [Le90]3), son
analogue aux diff\'erences (cf. [vdPS97]1.29,1.30), ainsi que les situations mixtes \'etudi\'ees par
Bialynicki-Birula [Bi62]. Il englobe le cas des syst\`emes \`a plusieurs variables, \'eventuellement
non int\'egrables. 

\vfill \eject
\centerline{\it Bibliographie.}

\bigskip \item {[An87]} Y. Andr\'e, {\it Quatre descriptions des groupes de Galois
diff\'erentiels}, S\'em. d'alg\`ebre de Paris 86/87, Springer L. N. M. 1296 (1987),
28-41.
\medskip \item {[An89[} Y. Andr\'e, {\it Notes sur la th\'eorie de Galois
diff\'erentielle}, Pr\'epubl. IHES/M/89/49 (1989).
\medskip \item {[An96]} Y. Andr\'e, {\it Pour une th\'eorie inconditionnelle des motifs}, Publ.
math. I.H.E.S. 83 (1996) 5-49.
\medskip \item {[Ao90]} K. Aomoto, {\it $q$-analogue of De Rham cohomology associated with
Jackson integrals}, Proc. japan Acad. 66 (1990) 161-164.
\medskip \item {[BeO78]} P. Berthelot, A. Ogus, {\it Notes on Crystalline Cohomology}, Math.
Notes N. 21,  Princeton Univ. Press (1978).
\medskip \item {[Ber92]} D. Bertrand, {\it Groupes alg\'ebriques et \'equations diff\'erentielles
lin\'eaires}, S\'em. Bourbaki, exp. 750, F\'ev. 1992. 
\medskip \item {[Beu92]} F. Beukers, {\it Differential Galois theory}, from Number Theory to Physics, M.
Waldschmidt, P. Moussa, J.M. Luck, C. Itzykson eds., Springer (1992).
\medskip \item {[Bi62]} A. Bialynicki-Birula, {\it On Galois theory of fields with operators},
Amer. J. of Math. 84 (1962) 89-109.
\medskip \item {[Bou81]} N. Bourbaki, {\it Alg\`ebre}, chapitres III, IV et X, Masson (1981).
\medskip \item {[Bou85]} N. Bourbaki, {\it Alg\`ebre commutative}, chapitre I, Masson
(1985).
\medskip \item {[Bru94]} A. Brugui\`eres, {\it Th\'eorie tannakienne non commutative}, Comm. in
Algebra 22 (14) (1994) 5817-5860.
\medskip \item {[Bry99]} R. Bryant, {\it Recent advances in the theory of holonomy}, S\'em.
Bourbaki, exp. 861, Juin 1999.
\medskip \item {[Coh71]} P. Cohn, {\it Free rings and their relations}, London Math. Soc. monogr.
2, (1971) Acad. Press.
\medskip \item {[Coh77]} P. Cohn, {\it Skew field constructions}, Cambridge Univ. press  (1977).
\medskip \item {[Con90]} A. Connes, {\it G\'eom\'etrie non commutative}, Inter\'editions (1990).
\medskip \item {[Con94]} A. Connes, {\it Non-commutative Geometry}, Acad. Press
(1994).
\medskip \item {[DG70]} M. Demazure, P. Gabriel, {\it Groupes alg\'ebriques 1}, North Holland
(1970).  
\medskip \item {[De89]} P. Deligne, {\it Le groupe fondamental de la droite projective moins
trois points}, Galois groups over $\bar Q$, M.S.R.I. Publ. 16 (1989), 79-297.
\medskip \item {[De90]} P. Deligne, {\it Cat\'egories tannakiennes}, Grothendieck
Festschrift, vol. 2, Birkh\"auser P.M. 87 (1990), 111-198.
\medskip \item {[DeMi82]} P. Deligne, J. Milne, {Tannakian categories}, Springer Lecture Notes 900
(1982)101-228. 
\medskip \item {[DubM96]} M. Dubois-Violette, T. Masson, {\it On the first order operators in
bimodules}, Lett. Math. Phys. 37 (1996) 464-474.
\medskip \item {[Dub97]} M. Dubois-Violette, {\it Some aspects of noncommutative differential
geometry}, Contemporary Math. vol. 203 (1997) 145-157.
\medskip \item {[Duv83]} A. Duval, {\it Lemmes de Hensel et factorisation formelle pour les
op\'erateurs aux diff\'erences}, Funkcial. Ekvac. 26 (1983) 349-368. 
\medskip \item {[E87]} T. Ekedahl, {\it Foliations and inseparable
morphisms}, Proc. Symp. Pure Math. 46 (1987) 139-149.
\medskip \item {[F92]} A. Fahim, {\it Extensions galoisiennes d'alg\`ebres diff\'erentielles}, C.R. Acad.
Sci. Paris, t. 314 (1992) 1-4.
\medskip \item {[GR90]} G. Gasper, M. Rahman, {\it Basic hypergeometric
series}, Cambridge Univ. press (1990).
\medskip \item {[He92]} Y. Hellegouarch, {\it Galois calculus and Carlitz exponentials},
The arithmetic of function fields, D. Goss, D. Hayes, M. Rosen eds., de Gruyter (1992) 33-50.
\medskip \item {[Kar87]} M. Karoubi, {\it Homologie cyclique et K-th\'eorie}, Ast\'erisque 149,
S.M.F. (1987).
\medskip \item {[Kar95]} M. Karoubi, {\it Formes diff\'erentielles non commutatives et
cohomologie \`a coefficients arbitraires}, Trans. A.M.S. 347 n. 11 (1995) 4277-4299.
\medskip \item{[Kat70]} N. Katz, {\it Nilpotent connections and the monodromy theorem.
Applications of a result of Turrittin},
Publ. Math. IHES 39 (1970), 175-232.
\medskip \item {[Kat90]} N. Katz, {\it Exponential sums and differential equations}, Annals of
Math. Studies 124 (1990) Princeton.
\medskip \item {[KN63]} S. Kobayashi, K. Nomizu, {\it Foundations of differential geometry},
Interscience I (1963).
\medskip \item {[Kol48]} E. Kolchin, {\it Algebraic matrix groups and the Picard-Vessiot theory of linear
ordinary differential equations}, Ann. of Math. 49 (1948) 1-42. 
\medskip \item {[LT99]} D. Laksov, A. Thorup, {\it These are the differentials of order $n$},
Trans. A.M.S. 351, n. 4 (1999) 1293-1353.
\medskip \item {[Le90]} A. Levelt, {\it Differential Galois
theory and tensor products}, Indag. Mathem. N.S.1 (1990) 439-450.
\medskip \item {[Li55]} A. Lichnerowicz, {\it Th\'eorie globale des connexions et des groupes
d'holonomie}, Cremona (1955). 
\medskip \item {[Mat99]} O. Mathieu, {\it Classification des alg\`ebres de
Lie simples}, S\'em. Bourbaki, exp. 858, Mars 1999.
\medskip \item {[Mou95]} J. Mourad, {\it Linear connections in noncommutative geometry}, Class. and Quant.
Grav. 12 (1995) 965-974.
\medskip \item {[N97]} P. Nuss, {\it Noncommutative descent and non-abelian cohomology}, K-theory
12 (1997) 23-74.
\medskip \item {[Pr83]} C. Praagman, {\it The formal classification of
linear difference operators}, Proc. Kon. Ned. Ac. Wet. Ser. A, 86 (1983) 249-261.
\medskip \item {[Saa72]} N. Saavedra Rivano, {\it Cat\'egories tannakiennes},
Springer L. N. M. 265 (1972). 
\medskip \item {[Sau99]} J. Sauloy, {\it Matrice de connexion d'un syst\`eme aux
$q$-diff\'erences confluant vers un syst\`eme diff\'erentiel de matrices de monodromie},
C.R.Acad. Sci. Paris, t. 328, S\'er. I (1999) 155-160. 
\medskip \item {[Se93]} J. P. Serre, {\it G\`ebres}, l'Ens. Math. 39 (1993) 33-85. 
\medskip \item {[Su99]} M. Suzuki, {\it General formulation of quantum analysis}, Reviews in
Math. Physics 11, n.2 (1999) 243-265. 
\medskip \item {[Ta39]} T. Tannaka, {\it \"Uber der Dualit\"atssatz der nichtkommutative
topologischen Gruppen}, T\^ohoku Math. J. 45 (1939) 1-12. 
\medskip \item {[TV97]} V. Tarasov, A. Varchenko, {\it Geometry of
$q$-hypergeometric functions, quantum affine algebras and elliptic quantum groups}, Ast\'erisque
246, S.M.F. (1997).  
\medskip \item {[vdP95]} M. van der Put, {\it Differential equations in
characteristic $p$}, Compos. Math. 97 (1995) 227-251.
\medskip \item {[vdPS97]} M. van der Put, M. Singer, {\it Galois theory of
difference equations}, Springer Lect. notes Math. 1666 (1997). 
\medskip \item {[W71]} F. Warner, {\it Foundations of differentiable manifolds and Lie groups},
Springer G.T.M. 94 (1971).

\end